\documentclass[11pt,a4paper,twoside]{article}
\usepackage{latexsym,amsfonts,amsmath, amsthm, amssymb}
\usepackage{fullpage}
\usepackage{xcolor,bm}
\usepackage[utf8x]{inputenc}
\usepackage{array}

\usepackage[normalem]{ulem}

\newtheorem{theorem}{Theorem}
\newtheorem{proposition}[theorem]{Proposition}
\newtheorem{lemma}[theorem]{Lemma}
\newtheorem{corollary}[theorem]{Corollary}

\theoremstyle{definition}

\newtheorem{definition}[theorem]{Definition}
\newtheorem{remark}[theorem]{Remark}
\newtheorem{remarks}[theorem]{Remarks}
\newtheorem{reminder}[theorem]{Reminder}

\newtheorem{notation}[theorem]{Notation}

\newcommand{\ls}{\leqslant}
\newcommand{\gs}{\geqslant}

\newcommand{\R}{\mathbb{R}}

\newcommand{\B}{\mathbf{B}}
\newcommand{\I}{\mathcal{I}}
\newcommand{\J}{\mathcal{J}}
\newcommand{\cZ}{\mathcal{Z}}

\newcommand{\II}{\mathfrak{I}}

\newcommand{\U}{\mathcal{U}}

\newcommand{\norm}[1]{\left\lVert{#1}\right\rVert}
\newcommand{\abs}[1]{\left|{#1}\right|}

\DeclareMathOperator{\sign}{\mathrm{sign}}
\DeclareMathOperator{\intr}{\mathrm{int}}

\begin{document}

\title{\bf Illuminating 1-unconditional convex bodies in $\R^3$ and $\R^4$, and certain cases in higher dimensions}

\medskip

\author{Wen Rui Sun and Beatrice-Helen Vritsiou}

\date{}

\maketitle

\begin{abstract}
\small
We settle the Hadwiger-Boltyanski Illumination Conjecture for all 1-unconditional convex bodies in $\R^3$ and in $\R^4$. Moreover, we settle the conjecture for those higher-dimensional 1-unconditional convex bodies which have at least one coordinate hyperplane projection equal to the corresponding projection of the circumscribing rectangular box. Finally, we confirm the conjectured equality cases of the Illumination Conjecture within the subclass of 1-unconditional bodies which, just like the cube $[-1,1]^n$, have no extreme points on coordinate subspaces.

Our methods are combinatorial, and the illuminating sets that we use consist primarily of small perturbations of the standard basis vectors. In particular, we build on ideas and constructions from \cite{Sun-Vritsiou-sym}, and mainly on the notion of \emph{deep illumination} introduced there.
\end{abstract}

\vspace{-0.5cm}

\section{Introduction}\label{sec:intro}

\vspace{-0.2cm}

This paper is a direct continuation of \cite{Sun-Vritsiou-sym}: building on the approach there, which allowed us to come up with a uniform way of illuminating 1-symmetric convex bodies of all dimensions in accordance to what the Illumination Conjecture stipulates, we extend this to certain cases of 1-unconditional convex bodies. 

Let $K$ be a \emph{convex body} in the Euclidean space $\R^n$, that is, a convex, compact set with non-empty interior. 
Given a boundary point $x$ of $K$ and a non-zero vector $d\in \R^n$ (a direction), we say that $d$ illuminates $x$ if there exists $\varepsilon > 0$ such that $x+\varepsilon d\in \intr K$. A set of directions ${\cal D}= \{d_1, d_2,\ldots, d_M\}$ such that,
\emph{for each boundary point $x$ of $K$, there is at least one $d_i\in {\cal D}$ which illuminates $x$,}
will be called an illuminating set for $K$. The smallest cardinality of an illuminating set for $K$ is called the illumination number of $K$, and we denote it by $\II(K)$. 

\smallskip

This definition of illumination is due to Boltyanski \cite{Boltyanski-1960}. There is also an equivalent definition by Hadwiger \cite{Hadwiger-1960}, where we illuminate using point `light sources' placed outside $K$ (and all rays emanating from them towards the boundary of $K$, which meet the boundary and \underline{then} cross into the interior of $K$); it can be shown that both definitions lead to the same number for a fixed body $K$. 

Moreover, we have that, for any convex body $K$, $\II(K) = N(K,\intr K)$, where the latter stands for the \emph{covering number of $K$ by $\intr K$} (that is, the smallest number of translates of $\intr K$ whose union contains $K$). Thus the Illumination Conjecture, which we formulate right below, is equivalent to Hadwiger's Covering Problem. Finally, there is yet another equivalent formulation by Gohberg and Markus \cite{Gohberg-Markus-1960}, where we cover $K$ by smaller homothetic copies of it.

\let\thefootnote\relax\footnote{\emph{Keywords:} illumination, symmetries of the cube, covering number, X-ray number, deep illumination}
\let\thefootnote\relax\footnote{\emph{2020 Mathematics Subject Classification:} 52A40, 52A37 (Primary); 52A20, 52C07 (Secondary)}

\vspace{-0.2cm}

\noindent {\bf Hadwiger's Covering Problem/The Hadwiger-Boltyanski Illumination Conjecture.} For every convex body $K$ in $\R^n$, we should have $\II(K) = N(K,\intr K)  \leq 2^n$.

\smallskip

Furthermore, the inequality should be strict, except in the case of the cube and of its affine images (parallelepipeds) in $\R^n$.

\medskip

An excellent reference on the history of these conjectures, and of related problems, and on progress up to recent years is the survey \cite{Bezdek-Khan-survey}. We also refer to the monographs \cite{Bezdek-CMSmonograph, BMS-book} and the surveys \cite{Bezdek-1992, Boltyanski-Gohberg-1995, Martini-Soltan-1999}. 

Levi in 1955 \cite{Levi-1955} fully settled the problem of bounding $N(K,\intr K)$ for planar convex bodies (showing that $N(K,\intr K) = 3$ for $K\subset \R^2$, except if $K$ is a parallelogram, in which case $N(K,\intr K) = 4$). Motivated by that, in 1957 Hadwiger \cite{Hadwiger-1957} posed the analogous question in higher dimensions. Still, aside from Levi's solution in $\R^2$, in all other dimensions the general problem is still open. In dimension 3 Lassak \cite{Lassak-1984} has shown that, if $K$ is centrally symmetric (that is, $K-x=x-K$ for some $x\in \R^3$), then $\II(K) \leq 8$. In other words, short of the equality cases, the conjecture in $\R^3$ is settled for symmetric convex bodies, but it remains open for the not-necessarily symmetric case, with the best bound being $14$ (due to Prymak \cite{Prymak-2023}). We also refer to a very recent paper by Arman, Bondarenko and Prymak \cite{ABP-2024}, where the reader can find all the progress to date and the most recent improvements on the bounds for other low dimensions.

A longstanding general upper bound (which remains the best known when specialised to the symmetric case) was already given in 1964 by Erd\"{o}s and Rogers \cite{Erdos-Rogers-1964}:
\begin{equation*}
\II(K) = N(K,\intr K) \leq \frac{{\rm vol}(K-K)}{{\rm vol}(K)}\theta(K) \leq \frac{{\rm vol}(K-K)}{{\rm vol}(K)}n\bigl(\ln n + \ln\ln n + 5\bigr)
\end{equation*}
where $\theta(K)$ is the asymptotic lower density of the most economical covering of $\R^n$ by copies (translates) of $K$. Erd\"{o}s and Rogers adapted an earlier proof by Rogers \cite{Rogers-1957} which was giving the first polynomial-order, and essentially best known to date, bound on $\theta(K)$. Combining this with the Rogers-Shephard inequality \cite{Rogers-Shephard-1957}, one is led to the bound $\II(K) \leq C4^n \sqrt{n}\ln n$ for every convex body $K\subset \R^n$, where $C$ is an absolute constant (moreover, in the symmetric case one gets $\II(K) \leq C^\prime 2^n n\ln n$). More recently, subexponential improvements to this general upper bound were given in \cite{HSTV-2021}, \cite{CHMT-2022} and \cite{Galicer-Singer-2023}, with the latter two papers attaining almost exponential improvements. The main novelty in these three papers is the use of results from Asymptotic Convex Geometry on the concentration of volume in high-dimensional convex bodies. Note however that neither the initial approach in \cite{HSTV-2021}, nor the more recent refinements, can contribute anything to the bound in the symmetric case, which would be the most relevant one for both this paper and \cite{Sun-Vritsiou-sym}.

\bigskip

The Illumination Conjecture has been fully settled for certain special classes of convex bodies. Again, we refer the reader to the survey \cite{Bezdek-Khan-survey} for a comprehensive list of references up to 2016. Just as examples, we mention that:
\begin{itemize}
\item Levi also showed in \cite{Levi-1955} that $\II(Q)=n+1$ for all smooth convex bodies $Q$ in $\R^n$. 
\item Martini \cite{Martini-1985} settled the conjecture for the class of belt polytopes (which contains the zonotopes). This was later extended by Boltyanski and Soltan \cite{Boltyanski-Soltan-1990, Boltyasnki-Soltan-1992} to zonoids, and by Boltyanski \cite{Boltyanski-1996} to belt bodies (see also \cite{Boltyanski-Martini-2001}).
\item The conjecture is fully settled for convex bodies of constant width. For dimensions $n\geq 16$, this is due to Schramm \cite{Schramm-1988}. For the remaining dimensions we have: \cite{Lassak-1997}, \cite{Weissbach-1996} (see also \cite[Section 11]{BLNP-2007}) dealing with $n=3$, \cite{Bezdek-Kiss-2009} dealing with $n=4$, and \cite{BPR-2022} dealing with $5\leq n\leq 15$.
\item Tikhomirov \cite{Tikhomirov-2017} settled the conjecture for 1-symmetric convex bodies of sufficiently large dimension (1-symmetric means that the body is invariant under reflections about coordinate subspaces and under any permutation of the coordinates). His result was the main motivation for \cite{Sun-Vritsiou-sym}, where the authors gave an alternative approach which also allows one to deal with 1-symmetric bodies in low dimensions.
\item Bezdek, Ivanov and Strachan \cite{BIS-2023} confirmed the conjecture for centrally symmetric cap bodies in dimensions $n=3$ (see also \cite{Ivanov-Strachan-2021}), $n=4$, and $n\geq 20$. They further showed that, if the cap body is 1-unconditional (we recall the definition below), then the Illumination Conjecture holds in all dimensions (and in that case $\II(K) \leq 4n$ once $n\geq 5$).
\item Gao, Martini, Wu and Zhang \cite{GMWZ-2024} verified the conjecture for polytopes which arise as the convex hull of the Minkowski sum of a finite subset of the lattice ${\mathbb Z}^n$ and of the unit-volume cube $\bigl[-\tfrac{1}{2},\,\tfrac{1}{2}\bigr]^n$.
\item Finally, Livshyts and Tikhomirov \cite{Livshyts-Tikhomirov-2020} settled the conjecture for convex bodies in sufficiently small neighbourhoods of the cube (with respect to either the geometric or the Hausdorff distance). Given that $\II(K) = N(K,\intr K)$ is an upper semicontinuous quantity (see e.g. \cite{Naszodi-2009}), the bound $2^n$ can already be deduced for bodies sufficiently close to $[-1,1]^n$, so their result is about settling the equality cases (and indeed they show that, if ${\rm dist}(K, [-1,1]^n)$ is small enough, and $K$ is not a parallelepiped, then $2^n -1$ is a \underline{sharp} upper bound for $\II(K)$).
\end{itemize}

Recall that a convex body $K$ in $\R^n$ is called \emph{$1$-unconditional} if it is invariant under reflections about coordinate subspaces. Equivalently if
\begin{center}
$x=(x_1,x_2,\ldots,x_n)\in K\ $ implies that
$\ (\epsilon_1 x_1, \epsilon_2 x_2,\ldots, \epsilon_n x_n)\in K$ 
\\
for any choice of signs $\epsilon_i\in \{\pm 1\}$, $1\leq i\leq n$.
\end{center} 
Relevant results that would apply to this class are the following, which however only deal with 3-dimensional convex bodies. Lassak \cite{Lassak-1984} showed that $\II(K)\leq 8$ for every centrally symmetric convex body $K$ in $\R^3$ (equivalently, for every origin-symmetric $K$, that is, such that $K=-K$). Moreover, he showed this while using illuminating sets formed by 4 pairs of opposite directions (and posed the question whether this is possible to do in higher dimensions as well, if $K=-K$). 

Bezdek \cite{Bezdek-1991} showed that $\II(P)\leq 8$ for any polytope in $\R^3$ which has a non-trivial affine symmetry. Finally Dekster \cite{Dekster-2000} obtained the same bound for any convex body $K$ in $\R^3$ which is symmetric about a plane. 

\smallskip

The main results of this paper are the following. Note that, given a convex body $K$, we will denote by ${\rm dim}(K)$ the dimension of the ambient Euclidean space. Moreover, assuming the dimension $n$ is clear from the context, we write ${\bm 1} = e_1+e_2+\cdots+e_n$, where $e_i$, $1\leq i\leq n$ are the standard basis vectors in $\R^n$.

\begin{theorem}\label{thm:main-dim3-dim4}
Let $K$ be a 1-unconditional convex body in $\R^3$ or $\R^4$, and assume that $K$ is not a parallelepiped. Then $\II(K) \leq 2^{{\rm dim}(K)} - 2$.

Moreover, we can use illuminating sets of this cardinality which consist of pairs of opposite directions.
\end{theorem}

Observe that, because of Lassak's and Dekster's results, the part of the above theorem which concerns dimension 3 is only novel in that we also settle the equality cases.

\begin{theorem}\label{thm:main-unit-subcubes}
Let $n\geq 3$, and let $K$ be a 1-unconditional convex body in $\R^n$. Assume without loss of generality that $e_i\in \partial K$ for all $1\leq i\leq n$ (see the next section on why this is WLOG).

In addition, suppose that there exists at least one $i_0\in \{1,2,\ldots,n\}$ such that the vector ${\bm 1}-e_{i_0}\in \partial K$ (in other words, $K$ contains at least one unit subcube of dimension $n-1$). Then, if $K$ is not a parallelepiped, we will have that $\II(K)\leq 2^n-2$.

Moreover, we can use illuminating sets of this cardinality which consist of pairs of opposite directions.
\end{theorem}

Similar to Theorem \ref{thm:main-unit-subcubes}, we also have the following

\begin{proposition}\label{prop:main-all-unit-subcubes-dim-n-2}
Let $n\geq 4$, and let $K$ be a 1-unconditional convex body in $\R^n$ such that $e_i\in \partial K$ for all $1\leq i\leq n$.

Assume that \underline{for all} $1\leq i < j\leq n$ we have that ${\bm 1}-e_i-e_j\in \partial K$ (in other words, $K$ contains \underline{all possible} unit subcubes of dimension $n-2$). Then, if $K$ is not a parallelepiped, we will have that $\II(K)\leq 2^n-2$.

Moreover, we can use illuminating sets of this cardinality which consist of pairs of opposite directions.
\end{proposition}

\begin{corollary}\label{cor:main-X-ray-conjecture}
Let $K$ be a 1-unconditional convex body satisfying any of the assumptions of Theorems \ref{thm:main-dim3-dim4} and \ref{thm:main-unit-subcubes} or of Proposition \ref{prop:main-all-unit-subcubes-dim-n-2}.

Then, given that we can illuminate $K$ using no more than $2^{{\rm dim}(K)-1}$ pairs of opposite directions, we can conclude that $K$ also satisfies the Bezdek-Zamfirescu $X$-ray conjecture, and that its $X$-ray number $X(K)\leq 2^{{\rm dim}(K)}$.
\end{corollary}

For details on this latter conjecture, see e.g. \cite{Bezdek-Zamfirescu-1994} and \cite{Bezdek-Kiss-2009}.

\smallskip

The final main result of this paper is the following

\begin{theorem}\label{thm:main-cubelike-bodies}
Let $n\geq 3$, and let $K$ be a 1-unconditional convex body in $\R^n$ with the following property:
\begin{equation*}
\hbox{if $x$ is an extreme point of $K$, then $x_i\neq 0$ for all $1\leq i\leq n$.} \tag{$\dagger$}
\end{equation*}
Then, if $K$ is not a parallelepiped, we will have that $\II(K) \leq 2^n-2$.

Moreover, we can use illuminating sets of this cardinality which consist of pairs of opposite directions.
\end{theorem}

\begin{remark}
ln this paper, 1-unconditional convex bodies which have Property ($\dagger$) will be called \emph{cubelike}. 

As we will recall in the next section, such convex bodies in $\R^n$ can be illuminated by $2^n$ directions (and in fact, they can be illuminated by any illuminating set of the cube in $\R^n$). Hence, this last theorem is about settling equality cases in this subclass of bodies.

We will also see that its proof relies on an inductive process, which however on its own can only recover the bound $2^n$; it's a combination of this inductive process and Theorem \ref{thm:main-unit-subcubes} that finally allows us to obtain the claimed result.
\end{remark}

\smallskip

The rest of the paper is organised as follows. For most of the 3-dimensional cases of Theorem \ref{thm:main-dim3-dim4}, a proof (or a proof sketch) is given in Section \ref{sec:3-dim-bodies}. The remaining cases are also special cases of Theorem \ref{thm:main-unit-subcubes}: all cases of this theorem, broken down into separate propositions, are proved in Section \ref{sec:maximal-unit-subcubes}. The proof of Proposition \ref{prop:main-all-unit-subcubes-dim-n-2} is also found at the end of this section. In Section \ref{sec:cubelike} we establish Theorem \ref{thm:main-cubelike-bodies}. Finally, the still unsettled 4-dimensional cases of Theorem \ref{thm:main-dim3-dim4}, which do not already follow as special cases of Theorem \ref{thm:main-unit-subcubes} and of Proposition \ref{prop:main-all-unit-subcubes-dim-n-2}, are handled in Section \ref{sec:4-dim-bodies}.

\bigskip

\noindent {\bf Acknowledgements.} Part of writing up the final version of this paper was done while the two authors were in residence at the Hausdorff Research Institute for Mathematics for the programme ``Synergies between modern probability, geometric analysis and stochastic geometry''. The authors are grateful to the institute and the organisers for the hospitality and the excellent working conditions. The second-named author is partially supported by an NSERC Discovery Grant.

\section{Preliminary results}\label{sec:prelims}

We write $[n]$ for the set $\{1,2,\ldots,n\}$, and $e_1, e_2,\ldots, e_n$ for the standard basis vectors of $\R^n$. For any vector $x\in \R^n$, we will denote by $\cZ_x$ the set $\{i\in [n]: x_i=0\}$. Also, we will write $\big|\vec{x}\big|$ for the vector $\sum_{i\in [n]}\abs{x_i}e_i$, namely the coordinate reflection of $x$ which has only non-negative coordinates.

Given a subset $A$ of $\R^n$, we will denote its interior and its boundary by $\intr A$ and by ${\rm bd} A$ or $\partial A$ respectively. Recall that if $A$ is a non-empty convex set, then its affine hull
\begin{equation*}
{\rm aff}(A) := \bigl\{\mu_1a_1+\mu_2a_2+\cdots+\mu_\ell a_\ell\,:\,  \ell \gs 1,\, a_i\in A \ \hbox{and}\ \mu_i\in \R \ \hbox{with}\ \mu_1 + \mu_2 +\cdots + \mu_\ell = 1\bigr\}
\end{equation*}
coincides with the smallest affine subspace of $\R^n$ which contains $A$, and that, in the subspace topology on ${\rm aff}(A)$, $A$ has non-empty interior. We call this the relative interior of $A$ and denote it by ${\rm relint} A$. Moreover, we call $A\setminus {\rm relint} A$ the relative boundary of $A$, and denote it by ${\rm relbd}A$.

\medskip

As already mentioned, a convex body $K$ in $\R^n$ is a convex set of $\R^n$ which is compact and has non-empty interior. If $K$ is also origin-symmetric, that is, $K=-K$, then $K$ is the unit ball of a certain norm on $\R^n$, which is given by $x\in \R^n \, \mapsto\, \|x\|_K:= \inf\{t>0: x\in tK\}$.

\smallskip

Recall that the illumination number of any convex body is an affine invariant: namely $\II(K) = \II(TK+z)$ for any invertible linear transformation $T\in {\rm GL}(n)$ and any (translation) vector $z$. 

Therefore, without loss of generality, we can assume that all the 1-unconditional convex bodies $\B\subset \R^n$ which we consider satisfy $e_i\in \partial \B$ for all $i\in [n]$, or equivalently that $\|e_i\|_\B = 1$ for all $i\in [n]$ (this is possible because, if $\B$ does not already have this property, it suffices to multiply it by the diagonal matrix ${\rm diag}\bigl(\|e_1\|_\B^{-1},\|e_2\|_\B^{-1},\ldots,\|e_n\|_\B^{-1}\bigr)$). We will denote this subclass of $n$-dimensional 1-unconditional convex bodies by $\U^n$.

\smallskip

Finally, we should mention the following fact about illumination (which we will often use in the sequel).

\noindent {\bf Fact A.} If a set ${\cal D}$ of directions illuminates all extreme points of a convex body $K$, then ${\cal D}$ illuminates $K$. 

\medskip

We recall \cite[Lemma 1, Corollary 2 and Remark 3]{Sun-Vritsiou-sym} and \cite[Lemma 6 and Corollary 7]{Sun-Vritsiou-sym} (they are now Lemma \ref{lem:smaller-coordinates}, Corollary \ref{cor:uncond-illum}, Remark \ref{rem:uncond-in-cube}, Lemma \ref{lem:affine-set} and Corollary \ref{cor:affine-set} respectively). Their proofs are standard, and are already given in \cite{Sun-Vritsiou-sym}, so we will not repeat them here.

\begin{lemma}\label{lem:smaller-coordinates}
Let $\B$ be a 1-unconditional convex body in $\R^n$. Suppose that $x$ is a point in $\B$, and that $y\in \R^n$ satisfies:
\begin{equation*}
\hbox{for all $i\in [n]$},\ \ |y_i|\leq |x_i|.
\end{equation*}
Then $y\in \B$ as well.

Moreover, if we have that
\begin{equation*}
\hbox{for all $i\in [n]$},\ \ |y_i| < |x_i| \ \hbox{or}\ |y_i|=|x_i|=0,
\end{equation*}
then $y\in \intr\B$.
\end{lemma}

\begin{corollary}\label{cor:uncond-illum}
Let $\B$ be a 1-unconditional convex body in $\R^n$, and let $x\in \partial \B$. Then $x$ is illuminated by any direction $d\in \R^n$ which satisfies
\begin{center}
$\cZ_d = \cZ_x$, and $d_i\cdot x_i < 0$ for all $i\in [n]\setminus \cZ_x$
\end{center}
(recall that $\cZ_x$ is the set of indices in $[n]$ which correspond to the zero coordinates of $x$). 

\smallskip

In particular, $\B$ is illuminated by the set $\{-1,0,1\}^n\setminus \{\vec{0}\}$. Furthermore, if $\B$ is cubelike (namely if it has Property $(\dagger)$ from Theorem \ref{thm:main-cubelike-bodies}), then $\B$ can be illuminated by the set $\{-1,1\}^n$ (here we also rely on Fact A).
\end{corollary}

\begin{remark}\label{rem:uncond-in-cube}
If $\B$ is a 1-unconditional convex body in $\R^n$, and $x\in \B$, then, by Lemma \ref{lem:smaller-coordinates}, we also have that $|x_i|e_i\in \B$ for all $i\in [n]$. 

Thus, if $\B\in \U^n$, then $\|x\|_\infty:=\max\limits_{i\in [n]}|x_i|\leq 1$. In other words, $\B\subseteq [-1,1]^n$.
 \end{remark}
 
 \begin{lemma}\label{lem:affine-set}
Let $K$ be a convex body in $\R^n$, and let $H$ be an affine subspace of $\R^n$. Suppose that $(\intr K) \cap H \neq \emptyset$. Then 
\begin{center}
${\rm relint}(K\cap H) = (\intr K) \cap H\,$ and $\,{\rm relbd}(K\cap H) = (\partial K)\cap H$.
\end{center}
\end{lemma}

\begin{corollary}\label{cor:affine-set}
Given the same general assumptions as in Lemma \ref{lem:affine-set}, consider $p\in {\rm relbd}(K\cap H)$, and a non-zero vector $d^\prime$ in the linear subspace $H-p\ls \R^n$ such that $p+\varepsilon d^\prime \in {\rm relint}(K\cap H)$ for some $\varepsilon > 0$. Then $p+\varepsilon d^\prime\in \intr K$.

In other words, if $p$ is $(K\cap H)$-illuminated by $d^\prime$, within $H= {\rm aff}(K\cap H)$, then it is also $K$-illuminated by $d^\prime$, viewed within $\R^n$ now.
\end{corollary}

In the sequel, we will also need the following

\begin{lemma}\label{lem:perturb}
Let $K$ be a convex body in $\R^n$, let $x_0\in \partial K$, and let $d_0$ be a direction in $\R^n$ which illuminates $x_0$. 
\begin{itemize}
\item[(a)] We can find $\rho>0$ such that, if $d^\prime\in \R^n\setminus\{\vec{0}\}$ satisfies $\|d_0-d^\prime\|_\infty < \rho$, then $d^\prime$ also illuminates $x$.
\smallskip
\item[(b)] We can find $\tau>0$ such that, for every $y\in \partial K$ which satisfies $\|x-y\|_\infty < \tau$, we will have that the direction $d_0$ illuminates $y$ as well.
\end{itemize}
\end{lemma}
\begin{proof} Fix $\varepsilon_0 > 0$ such that $x_0+\varepsilon_0d_0\in \intr K$. Then we can find $\eta_0 >0$ such that
\begin{equation*}
\bigl\{z\in \R^n: \|(x_0+\varepsilon_0d_0)-z\|_\infty < \eta_0\bigr\}\subseteq \intr K.
\end{equation*}
Now, set $\rho= \frac{1}{\varepsilon_0}\eta_0$, and consider $d^\prime\in \R^n\setminus\{0\}$ such that $\|d_0-d^\prime\|_\infty < \rho$. For $z_0 = x_0 + \varepsilon_0d^\prime$, we will have
\begin{equation*}
\big\|(x_0+\varepsilon_0d_0) - z_0\big\|_\infty = \|\varepsilon_0(d_0-d^\prime)\|_\infty = \varepsilon_0 \|d_0-d^\prime\|_\infty < \eta_0,
\end{equation*}
which shows that $z_0=x_0 + \varepsilon_0d^\prime\in \intr K$. In other words, $d^\prime$ illuminates $x_0$ too, which completes the proof of part (a).

\medskip

Similarly, set $\tau= \eta_0$. Suppose that $y\in \partial K$ and satisfies $\|x_0-y\|_\infty < \tau$. Then, for $z_1 = y+\varepsilon_0d_0$, we have
\begin{equation*}
\big\|(x_0+\varepsilon_0d_0) - z_1\big\|_\infty = \|x_0-y\|_\infty < \eta_0,
\end{equation*}
and hence $z_1 = y+\varepsilon_0d_0\in \intr K$. In other words, $d_0$ illuminates the boundary point $y$ too, which shows part (b).
\end{proof}

\subsection{A brief review of tools from \cite{Sun-Vritsiou-sym}: illuminating 1-symmetric convex bodies in all dimensions}\label{subsec:sym-results}

Recall that 1-symmetric convex bodies are a subclass of 1-unconditional convex bodies: a body $\B\subset \R^n$ is called 1-symmetric if
\begin{center}
$x=(x_1,x_2,\ldots,x_n)\in K\ $ implies that
$\ (\epsilon_1 x_{\sigma(1)}, \epsilon_2 x_{\sigma(2)},\ldots, \epsilon_n x_{\sigma(n)})\in K$ 
\\
for any choice of signs $\epsilon_i\in \{\pm 1\}$, $1\ls i\ls n$, and {\bf any permutation $\sigma$ on $[n]$}.
\end{center} In \cite{Sun-Vritsiou-sym} we dealt with illuminating 1-symmetric convex bodies in all dimensions, and for this purpose we introduced the notion of \textit{deep illumination}: given $\delta\in (0,1)$, consider the set
\begin{equation}\label{def:big-illum-set}
    G^n(\delta) := \left\{d\in \R^n: \exists\, i\in [n]\ \hbox{such that}\ d=\pm e_i\,+\!\sum_{j \in [n] \backslash \{i\}} 
   \! \pm \delta e_j \right\},
\end{equation}
of directions in ${\mathbb R}^n$. A direction $d\in G^n(\delta)$ is said to \textit{deep illuminate} a non-zero vector $x\in {\mathbb R}^n$ if (i) $\sign(x_i) = -\sign(d_i)$ for every $i \in [n]\backslash\cZ_x$ and (ii) $1=\norm{d}_\infty = \abs{d_j}$ for an index $j \in [n]\backslash \cZ_x$. A subset $S$ of $G^n(\delta)$ is said to deep illuminate a subset $A$ of ${\mathbb R}^n\setminus \{0\}$ if every $y\in A$ is deep illuminated by some $d_y\in S$.

\medskip

Subsequently we showed that there exists a smaller subset $\I^n(\delta)$ of $G^n(\delta)$ with cardinality $2^n$ which still deep illuminates every non-zero vector in ${\mathbb R}^n$. We gave two explicit constructions, a more geometric one, and a purely combinatorial and recursive one. 

In this paper, when we write $\I^n(\delta)$, we will exclusively refer to the 2nd type of construction, which we recall below. This is partly because it will be easier to define/describe variations of this set, which, as we can then show, we can use as illuminating sets in different settings.

\begin{reminder}\label{reminder:In(delta)} \emph{(Construction of $\I^n(\delta)$ from \cite{Sun-Vritsiou-sym})}

\noindent (i) Check that $\I^2(\delta):= \bigl\{\pm (1,\delta),\,\pm (\delta, -1)\bigr\}$ deep illuminates every non-zero vector in ${\mathbb R}^2$.
\medskip\\
(ii) Construct $\I^{n+1}(\delta)$ from $\I^n(\delta)$ as follows:
out of the $2^{n+1}$ directions that $\I^{n+1}(\delta)$ will have in the end, the first $2^n$ are formed by appending to each direction $\mathbf{d}^n_s$ of $\I^n(\delta)$ one more `small' coordinate at the end, so that this new coordinate will have the same sign as the last coordinate of $\mathbf{d}^n_s$; that is,
\begin{equation*}
\mathbf{d}^{n+1}_s : = \bigl(\mathbf{d}^n_s,\,\sign(\mathbf{d}^n_{s,n})\delta\bigr).
\end{equation*} 
At the same time, this direction $\mathbf{d}^n_s$ allows us to also define one of the remaining $2^n$ directions for $\I^{n+1}(\delta)$, which we will denote by $\mathbf{d}^{n+1}_{2^n+s}$: the sign of each of the first $n$ coordinates of $\mathbf{d}^{n+1}_{2^n+s}$ will be the same as for the respective coordinate of $\mathbf{d}^n_s$, while the last coordinate of $\mathbf{d}^{n+1}_{2^n+s}$ will be equal to 1 in absolute value and will have opposite sign to the previous coordinate, the $n$-th one. That is,
\begin{equation*}
\mathbf{d}^{n+1}_{2^n+s} := \Bigl(\delta\cdot\Bigl(\sign(\mathbf{d}^n_{s,1}),\, \sign(\mathbf{d}^n_{s,2}),\ldots,\sign(\mathbf{d}^n_{s,n-1}),\,\sign(\mathbf{d}^n_{s,n})\Bigr),\ \,-\!\!\!\!-\sign(\mathbf{d}^n_{s,n})\Bigr).
\end{equation*}
Given this construction, we can inductively check that $\I^n(\delta)$ deep illuminates $\R^n$.
\end{reminder}

If we consider a 1-symmetric convex body $\B$ in ${\mathbb R}^n$ and $x\in \partial \B$, it is not hard to check that, if a direction $d\in G^n(\frac{1}{n+1})$ deep illuminates the vector $x$, then $x+\varepsilon d\in \intr\B$ for some $\varepsilon>0$ (see \cite[Lemma 11]{Sun-Vritsiou-sym}). Thus, $\B$ is illuminated by both the set $G^n(\frac{1}{n+1})$ and by its smaller subset $\I^n(\frac{1}{n+1})$ (since each of these sets deep illuminates all the boundary points of $\B$).

\medskip

Moreover, by examining more carefully 
\begin{enumerate}
\item which directions are included in the set $\I^n(\frac{1}{n+1})$ when it is constructed as above, and
\item in which other cases deep illumination is guaranteed to imply illumination,
\end{enumerate}
we could establish the following

\begin{theorem}\emph{(\cite[Theorem 23]{Sun-Vritsiou-sym})}
Let $n\geq 3$, and let $\B$ be a 1-symmetric convex body in $\R^n$ which is not affinely equivalent to the cube. WLOG assume that $\B\in \U^n$. Then we can find a minimal $\alpha_\B > 1$ such that $\B\subset [-1,1]^n\subset \alpha_\B\B$.

We can illuminate $\B$ using the set
\begin{gather*}
\Bigl[\I^n\bigl(\tfrac{1}{n+1}\bigr) \backslash\Bigl\{\pm\bigl(+1,+\tfrac{1}{n+1},+\tfrac{1}{n+1},\ldots,+\tfrac{1}{n+1},+\tfrac{1}{n+1},+\tfrac{1}{n+1}\bigr), \  \phantom{\I^n\bigl(\frac{1}{n+1}\bigr)}
\\
\phantom{\pm\bigl(+1,+\tfrac{1}{n+1},\ldots,+\tfrac{1}{n+1},+\tfrac{1}{n+1}\bigr)} \pm\Bigl(+\tfrac{1}{n+1},+\tfrac{1}{n+1},+\tfrac{1}{n+1},\ldots,+\tfrac{1}{n+1},{\bf -\!\!\!\!\!-\frac{1}{n+1}},+1\Bigr)\Bigr\}\Bigr]  \phantom{\I^n\bigl(\frac{1}{n+1}\bigr)}
\\[0.3em]
\bigcup \ \,\Bigl\{\pm\bigl(+1,+\tfrac{1}{n+1},+\tfrac{1}{n+1},\ldots,+\tfrac{1}{n+1},\,{\bm \eta},\, +\tfrac{1}{n+1}\bigr)\Bigr\}
\end{gather*}
for some $\eta\in (0,\frac{1}{n+1})$ (which will only depend on how close $\alpha_\B$ is to $1$).
\end{theorem}

One of the main ingredients in the proof of this theorem is the following

\noindent {\bf Fact B.} \emph{(\cite[Lemma 24]{Sun-Vritsiou-sym})} The following subset of $\I^n(\delta)$ ($\delta\in (0,1)$):
\begin{equation}
\I_{-2}^n(\delta) := \I^n(\delta)\setminus \bigl\{\pm(+\delta,+\delta,\ldots,+\delta, +\delta, \,{\bf -\!\!\!\!\!-}{\bm \delta},\, + 1)\bigr\}
\end{equation}
deep illuminates every vector $x\in \R^n\setminus\{\vec{0}\}$ which has at least one zero coordinate, that is, every vector $x$ with $1\leq \abs{\cZ_x}\leq n-1$.

\medskip

A further observation that we will need here is the following: if $\B$ is 1-symmetric and in $\U^n$, and moreover 
\begin{equation*}
{\bm 1}-e_i\in \B
\end{equation*}
for all $i\in [n]$, while ${\bm 1}\notin \B$, then $\B$ can more simply be illuminated by the set $\I_{-2}^n(\delta)$ directly, 
with $\delta\in (0,1)$ only depending on $\|{\bm 1}\|_\B$.

As we will see (Proposition \ref{prop:all-n-1-tuples}), this can be extended to 1-unconditional convex bodies with the same properties (note that these don't necessarily have to be 1-symmetric too, because for example we could have $\B$ contain some points of the form ${\bm 1}-\frac{1}{2}{\bm e}_i$, say, but not for all $i\in [n]$, and this wouldn't violate 1-unconditionality, but it would break 1-symmetry).

\section{1-unconditional convex bodies in $\R^3$} \label{sec:3-dim-bodies}

Observe that if $\B\in \U^3$ is NOT affinely equivalent to the cube, then, by our convention, $e_i\in \partial \B$ for all $i\in [3]$, while ${\bm 1} = e_1+e_2+e_3\notin \B$. Thus we can separate the different 3-dimensional cases of Theorem \ref{thm:main-dim3-dim4} into four groups, based on whether any coordinate permutations of $e_1+e_2$ are contained in $\B$, and if yes, how many (equivalently, based on whether $\B$ contains any \underline{unit} squares (2-dimensional subcubes), and how many).

We summarise the conclusions for each of these groups of cases in the following theorem (the numerical subscripts correspond to the numbering of the proposition(s) where each illuminating set appears).

\begin{theorem}\label{thm:R^3-summary}
Let $\B\in \U^3$ which is not a parallelepiped. Then $\B$ is illuminated by a coordinate permutation of one of the following sets:
\begin{align*}
&\phantom{\hbox{or}}\  {\cal F}_{\ref{prop:R^3-no-pairs},\ref{prop:R^3-exactly-one-pair},\delta}:= \bigl\{\pm(1,\delta,0),\ \pm(-\delta,1,0),\ \pm(0,0,1)\bigr\} &&\substack{\hbox{\small cases with no unit squares} \\ \hbox{\small or one unit square}}
\\
\intertext{or}
&\substack{\hbox{\normalsize ${\cal F}^1_{\ref{prop:R^3-exactly-two-pairs},\delta_1}:=\bigl\{\pm (\delta_1,\delta_1, 1), \pm (\delta_1,\delta_1, -1), \pm(-\delta_1, 1, 0)\bigr\}$} \\[0.2em] \hbox{\normalsize or ${\cal F}^2_{\ref{prop:R^3-exactly-two-pairs},\epsilon_2,\delta_{\epsilon_2}}: = \big\{\pm(\epsilon_2,1,1),\pm(-\delta_{\epsilon_2},1,\delta_{\epsilon_2}), \pm(-\delta_{\epsilon_2},-\delta_{\epsilon_2}, 1)\bigr\}$}} && \hbox{\small cases with two unit squares}
\\
\intertext{or}
&\phantom{\hbox{or}}\  {\cal F}_{\ref{prop:all-n-1-tuples},\delta} : = \bigl\{\pm(1,\delta,\delta),\pm(\delta,-1,-\delta),\pm(\delta,\delta,-1)\bigr\} && \hbox{\small cases with three unit squares},
\end{align*}
where the relevant parameter(s) $\delta,\delta_1,\epsilon_2,\delta_{\epsilon_2} >0$ should also be chosen based on $\B$ (in explicit ways, as we will see).
\end{theorem}

The theorem will follow from the proofs of Propositions \ref{prop:R^3-no-pairs} and \ref{prop:R^3-exactly-one-pair}, of Proposition \ref{prop:R^3-exactly-two-pairs}, which is treated as a special case of Proposition \ref{prop:exactly-two-n-1-tuples}, and of Proposition \ref{prop:all-n-1-tuples}. We will use the following terminology for the last two cases mentioned here: they concern bodies $\B\in \U^3$ which contain exactly two \emph{maximal} unit subcubes, or all possible \emph{maximal} unit subcubes, respectively (where `subcube' implies proper inclusion here, and where `maximality' is in terms of dimension). We treat these cases in the next section, in Propositions \ref{prop:exactly-two-n-1-tuples} and \ref{prop:all-n-1-tuples}, proving the analogous results in an arbitrary dimension $n\geq 3$.


\medskip

\begin{proposition}
\label{prop:R^3-no-pairs}
Let $\B \in \mathcal{U}^3$ and suppose 
\begin{equation*}
\hbox{$\norm{e_i+e_j}_\B > 1$ for every $i,j \in [3]$.} \tag{$\ast$} 
\end{equation*}
Then there exists $\delta > 0$ so that $\B$ can be illuminated by some coordinate permutation of the set
\begin{equation*}
{\cal F}_{\ref{prop:R^3-no-pairs},\ref{prop:R^3-exactly-one-pair},\delta}:= \bigl\{\pm(1,\delta,0),\ \pm(-\delta,1,0),\ \pm(0,0,1)\bigr\}.
\end{equation*}
\end{proposition}
\begin{proof}
For every $i\in [3]$ and for every $j\in [3]\backslash \{i\}$, set $a_{i,j}$ to be the supremum of non-negative numbers $x_j$ such that
\begin{equation*}
e_i + x_j e_j \in \B
\end{equation*}
(note that the set of such numbers is nonempty as $x_j=0$ belongs to it). Observe that by compactness the vector $e_i + a_{i,j}e_j$ is also in $\B$, and so are all its coordinate reflections. By assumption $(\ast)$, we have $a_{i,j} < 1$ for all $i\in [3]$ and for every $j\in [3]\backslash \{i\}$.

\bigskip

Let $a_{i_0,j_0}$ be the maximum of these numbers (not necessarily unique). For the rest of the proof, WLOG, we assume that $\{i_0,j_0\} = \{1,2\}$. We fix some $\delta < \frac{1-a_{i_0,j_0}}{2}$ and consider the corresponding set ${\cal F}_{\ref{prop:R^3-no-pairs},\ref{prop:R^3-exactly-one-pair},\delta}$.

Consider a boundary point $y=(y_1,y_2,y_3)$ of $\B$. We will show how to illuminate $y$ based on the number $\abs{\cZ_y}$ of zero coordinates of $y$.
\begin{itemize}
\item If $\abs{\cZ_y}=2$, then necessarily $y = ±e_s$ for some $s\in [3]$. If $y=\pm e_3$, then $\mp e_3$ from ${\cal F}_{\ref{prop:R^3-no-pairs},\ref{prop:R^3-exactly-one-pair},\delta}$ illuminates $y$. If $y=\pm e_1$, say $y=e_1$, then $y+(-1,-\delta,0) = (0,-\delta,0)\in \intr\B$ given that $\delta \in (0,1)$. Similarly we deal with $y=-e_1$ or $y=\pm e_2$.
\item Assume now that $\abs{\cZ_y}=1$. 
\begin{itemize}
\item Suppose first that $y = (y_1,y_2,0)$. Then we choose the direction $d_y$ in ${\cal F}_{\ref{prop:R^3-no-pairs},\ref{prop:R^3-exactly-one-pair},\delta}$ which has 1st and 2nd entries non-zero and with opposite signs to the corresponding entries of $y$. Clearly $y+\epsilon d_y \in \intr\B$ for $\epsilon \in (0, \min(|y_1|, |y_2|))$ (recall Corollary \ref{cor:uncond-illum}). 
\item If $y=(y_1,0,\pm 1)$ or $(0,y_2,\pm 1)$, we use the directions $\mp e_3$ (since by our assumptions we must have $|y_1| < 1$ or $|y_2|<1$ in these cases). We can use the same directions if $y=(y_1,0,y_3)$ or $(0,y_2,y_3)$ with $y_i$ representing numbers in $(-1,1)\setminus\{0\}$ here. 
\item Finally, suppose that $y$ is a boundary point of the form $(1,0,y_3)$. Then
\begin{equation*}
(1,0,y_3) + (-1,-\delta,0) = (0,-\delta,y_3).
\end{equation*}
Based on our notation, and because of the 1-unconditionality, we have $|y_3| \leq a_{1,3}\leq a_{i_0,j_0} < 1$. By convexity the vector
\begin{equation*}
\tfrac{1-a_{i_0,j_0}}{2}\,e_2 \ +\  \tfrac{1+a_{i_0,j_0}}{2}\, e_3
\end{equation*} 
is in $\B$. Moreover, it has strictly larger 2nd and 3rd entries in absolute value compared to $(1,0,y_3) + (-1,-\delta,0)$, thus, by Lemma \ref{lem:smaller-coordinates}, the latter vector is in $\intr\B$. 
\item Analogously to the last subcase, we deal with points of the form $(-1,0,y_3)$ and $(0,\pm 1,y_3)$.
\end{itemize}
\item Finally we consider a boundary point $y$ with $\abs{\cZ_y} =0$. 
\begin{itemize}
\item If $|y_3| <1$, we can choose the direction $d_y$ in ${\cal F}_{\ref{prop:R^3-no-pairs},\ref{prop:R^3-exactly-one-pair},\delta}$ which has 1st and 2nd entries non-zero and with opposite signs to the corresponding entries of $y$, and then employ Corollaries \ref{cor:uncond-illum} and \ref{cor:affine-set}. 
\item In the opposite case we have $y=(y_1, y_2,\pm 1)$; say, for illustration purposes, $y_3 = 1$. Note that, by the 1-unconditionality, the points $(|y_1|,0,1)$ and $(0,|y_2|,1)$ are in $\B$, thus $|y_1| \leq a_{3,1}\leq a_{i_0,j_0}$ and $|y_2| \leq a_{3,2}\leq a_{i_0,j_0}$. Hence, $(y_1,y_2,1) + (0,0,-1) = (y_1,y_2,0)$ and at least one of the following holds: $|y_1| < 1$ and $|y_2|\leq a_{1,2} < 1$, or $|y_2| < 1$ and $|y_1| \leq a_{2,1} < 1$. Thus, regardless of whether $(i_0,j_0) = (1,2)$ or $(2,1)$, we can use Lemmas \ref{lem:smaller-coordinates} and \ref{lem:affine-set} (with $\B\cap \{\xi\in \R^3: \xi_2=y_2\}$ or $\B\cap \{\xi\in \R^3:\xi_1=y_1\}$, respectively) to conclude that $(y_1,y_2,0)=y + (0,0,-1)\in \intr\B$.
\end{itemize}
\end{itemize}

In cases where $\{i_0,j_0\}\neq \{1,2\}$ (that is, if neither $a_{1,2}$ nor $a_{2,1}$ is the maximum of the numbers we defined above), we set $t_0$ for the remaining element of $[3]$ and consider the linear map/coordinate permutation
\begin{equation*}
\iota:\R^3\to \R^3,\quad x_{i_0}e_{i_0} + x_{j_0}e_{j_0} + x_{t_0}e_{t_0} \mapsto (x_{i_0}, x_{j_0}, x_{t_0}).
\end{equation*}
Then $\iota(\B)$ is illuminated by ${\cal F}_{\ref{prop:R^3-no-pairs},\ref{prop:R^3-exactly-one-pair},\delta}$ as above, and thus $\B$ is illuminated by $\iota^{-1}({\cal F}_{\ref{prop:R^3-no-pairs},\ref{prop:R^3-exactly-one-pair},\delta})$.
\end{proof}

\smallskip

Note that in the previous proposition we did not have to add the assumption that $\B$ is not an affine image of the cube: this is in fact implied from the other assumptions (namely that $\B$ is in $\U^3$ and does not contain any unit squares), which were enough to verify that $\B$ can be illuminated by (at most) 6 directions. 

The same happens in the next proposition.


\begin{proposition}
\label{prop:R^3-exactly-one-pair}
Let $\B \in \mathcal{U}^3$ and suppose that there is exactly one pair of indices $i_1,i_2\in [3]$ such that $\norm{e_{i_1}+e_{i_2}}_\B = 1$. Then there exists $\delta > 0$ so that $\B$ can be illuminated by some coordinate permutation of the set
\begin{equation*}
{\cal F}_{\ref{prop:R^3-no-pairs},\ref{prop:R^3-exactly-one-pair},\delta}:= \bigl\{\pm(1,\delta,0),\ \pm(-\delta,1,0),\ \pm(0,0,1)\bigr\}.
\end{equation*}
\end{proposition}
\begin{proof}
Let us assume WLOG that $(1,1,0)\in \B$. Moreover, let $a_1, a_2\geq 0$ be maximum possible such that $(1,0,a_1)$ and $(0,1,a_2)\in \B$ (by our main assumption we have that $a_1,a_2\in [0,1)$). WLOG we can assume that $a_1\geq a_2$.

Here we distinguish three main cases for boundary points $y=(y_1,y_2,y_3)$ of $\B$: (i) $\abs{\cZ_y} \geq 1$ and $|y_3| < 1$, (ii) $\abs{\cZ_y} = 0$ and $|y_3| < 1$, and (iii) $|y_3|=1$. We start with the following key observation.
\begin{itemize}
\item[$\blacklozenge$] We can choose $\delta_1 > 0$ small enough so that the directions $\pm(1,\delta,0)$ will illuminate both the point $(1,0,a_1)$ and all its coordinate reflections whenever $\delta \leq \delta_1$. Indeed, we have e.g. that $(1,0,a_1) + (-1,-\delta,0) = (0,-\delta, a_1) = -\delta e_2 + a_1 e_3 \in \intr\B$, as long as $0< \delta < 1-a_1$, so we can set $\delta_1 = (1-a_1)/2$.
\item[$\blacklozenge$] At the same time, as long as $\delta \in (0,1)$, the pair $\pm(-\delta,1,0)$ of directions illuminates $(0,1,a_2)$ and its coordinate reflections: this is because e.g. $(0,1,a_2) + (\delta,-1,0) = (\delta,0,a_2)$, which has strictly smaller 1st coordinate compared to $(1,0,a_1)$ and at most as large 3rd coordinate, which is also $<1$ (so we can either use Lemma \ref{lem:smaller-coordinates} alone, or combine it with Lemma \ref{lem:affine-set}). 
\end{itemize}

Now, suppose that $y=(y_1,0,y_3)$ with $0<|y_3|<1$.
\begin{itemize}
\item If $|y_1|=1$, then necessarily (by the maximality of $a_1$) we have that $|y_3|\leq a_1$. Thus:
\begin{itemize}
\item if $|y_3|=a_1$, we have already seen that the directions $\pm(1,\delta,0)$ illuminate $y$ (as long as $\delta\leq \delta_1$); 
\item if $|y_3| < a_1$, then $y$ is in the convex hull of the point $(1,0,a_1)$ and its coordinate reflections, and hence it is also illuminated by the directions $\pm(1,\delta,0)$. 
\end{itemize}
\item If $|y_1|<1$, then the direction $(0,0,-\sign(y_3))$ illuminates $y$ (since $(y_1,0,0)\in\intr\B$).
\end{itemize}

\smallskip

In an analogous way, we see that the directions $\pm(-\delta,1,0),\pm(0,0,1)$ illuminate boundary points of the form $y=(0,y_2,y_3)$ with $0<|y_3|<1$.

\medskip

\noindent Furthermore, the directions $\pm(1,\delta,0),\ \pm(-\delta,1,0)$, $\delta\leq \delta_1 < 1$, also illuminate: 
\begin{itemize}
\item the point $(1,1,0)$ and its coordinate reflections (and thus also any boundary point $y$ which satisfies $y_3=0$, since this will be in the convex hull of the former points);
\item any boundary point $y$ of $\B$ which satisfies $\abs{\cZ_y}=0$ and $|y_3| < 1$.
\end{itemize}

This takes care of the first two cases of boundary points in our breakdown.
It remains to figure out how to illuminate boundary points $y$ satisfying $|y_3|=1$. In such a case, by our main assumption we have $|y_1| <1$ and $|y_2|<1$, and thus the point $(y_1,y_2,0)\in \intr\B$ (we see this if we compare with the point $(1,1,0)\in \B$). Thus the direction $(0,0,-\sign(y_3))$ illuminates $y$.
\end{proof}


\smallskip

In contrast with the previous two propositions, the main assumptions in the next one can also be satisfied by affine images of the cube in $\R^3$. In fact, all bodies in $\U^3$ which satisfy the main assumption will contain such an affine image of the cube (which will also be in $\U^3$). Therefore, we have to explicitly rule out the cases where we don't have strict inclusion, and in the remaining cases we have to make crucial use of `special' boundary points which verify the strict inclusion. 

\begin{proposition}\label{prop:R^3-exactly-two-pairs}
Suppose that $\B\in {\cal U}^3$ is not an affine image of the cube but has the property that, for some permutation $(i_1,i_2,i_3)$ of $[3]$, $\norm{e_{i_1}+e_{i_2}}_\B = \norm{e_{i_1}+e_{i_3}}_\B = 1$, while $\norm{e_{i_2}+e_{i_3}}_\B>1$. Then there exist $\delta_1, \eta_{\delta_1}>0$ or $\epsilon_2, \delta_{\epsilon_2}>0$ such that $\B$ can be illuminated by a coordinate permutation of one of the following sets:
\begin{equation*}
{\cal F}^1_{\ref{prop:R^3-exactly-two-pairs},\delta_1} = \I_{ex2,1}^3(\delta_1) := \bigl\{\pm (\delta_1,\delta_1,\pm1), \pm(-\delta_1,1,0)\bigr\}
\end{equation*}
or 
\begin{equation*}
{\cal F}^2_{\ref{prop:R^3-exactly-two-pairs},\epsilon_2,\delta_{\epsilon_2}} = \I_{ex2,2}^3(\epsilon_2,\delta_{\epsilon_2}) := \big\{\pm(\epsilon_2,1,1), \pm(-\delta_{\epsilon_2},1,\delta_{\epsilon_2}), \pm(-\delta_{\epsilon_2},-\delta_{\epsilon_2}, 1)\bigr\}.
\end{equation*}
\end{proposition}

As before, we can check that, after applying some coordinate permutation on $\B$ (something which wouldn't ruin our main assumptions), we would be able to use one of the above sets exactly. We can thus assume WLOG that $\B$ contains the points $(1,1,0)$ and $(1,0,1)$ but not the point $(0,1,1)$. We first need the following

\begin{lemma}\label{lem:R^3-exactly-two-pairs-interior}
Suppose that $\B\in {\cal U}^3$ satisfies:
\begin{itemize}
\item $\|(1,1,0)\|_{\B}=\|(1,0,1)\|_{\B} = 1$, while $\|(0,1,1)\|_{\B}>1$;
\item $\B$ is NOT an affine image of the cube.
\end{itemize}
Then, for any $\epsilon\in (0,1]$, the point
\begin{equation*}
\bigl(1-\epsilon,\,\tfrac{1}{2},\,\tfrac{1}{2}\bigr)
\end{equation*} 
is an interior point of $\B$.
\end{lemma}
\begin{proof}
Since $\B$ contains the points $(1,1,0)$ and $(1,0,1)$ and all their coordinate reflections, it will contain their convex hull too, which is the set
\begin{equation*}
\bigl\{(x_1,x_2,x_3)\in {\mathbb R}^3: |x_1|\leq 1,\,|x_2|+|x_3|\leq 1\bigr\} = [-1,1]\times CP_1^2.
\end{equation*}
This is an affine image of the cube, therefore, by our last assumption for $\B$, we must have that $\B\setminus\bigl([-1,1]\times CP_1^2\bigr)\neq \emptyset$. 

Combined with the assumption that $\B\in {\cal U}^3$, this implies that $\B$ contains a point of the form $(0,z_2,z_3)$ where $z_2, z_3\in (0,1]$ and $z_2 + z_3 > 1$. From this we can obtain that
\begin{equation}\label{eqp1:lem:R^3-exactly-two-pairs-interior}
\bigl(0,\tfrac{1}{2},\tfrac{1}{2}\bigr)\in \intr\B
\end{equation}
as follows. WLOG we can assume that $z_2\geq z_3$. 
\begin{itemize}
\item[-] If $z_3 > \frac{1}{2}$, then \eqref{eqp1:lem:R^3-exactly-two-pairs-interior} follows immediately by Lemma \ref{lem:smaller-coordinates}.
\item[-] If $z_3\leq \frac{1}{2}$, then we can observe the following: $z_2 + z_3 > 1$ implies that $z_2 > 1-z_3 \geq \frac{1}{2}$, and hence $z_2+(1-z_3) > 2(1-z_3)\geq 1$. Therefore, we can find $\lambda\in (0,1)$ such that
\begin{equation*}
\lambda\cdot \bigl(z_2+(1-z_3)\bigr) = 1 \quad \Rightarrow\quad \lambda z_2 = \lambda z_3 + (1-\lambda).
\end{equation*}
From this we obtain that the point $(0,y_2,y_3)=\lambda (0,z_2,z_3)+ (1-\lambda)(0,0,1)\in \B$ has equal 2nd and 3rd coordinates, and moreover that
\begin{equation*}
y_2+y_3= \lambda z_2 + \lambda z_3 + (1-\lambda) = \lambda(z_2+z_3) + (1-\lambda) > 1.
\end{equation*}
In other words, $y_2=y_3 > \frac{1}{2}$, whence \eqref{eqp1:lem:R^3-exactly-two-pairs-interior} follows as in the previous case.
\end{itemize}

Finally, we note that $\bigl(1,\tfrac{1}{2},\tfrac{1}{2}\bigr) = \frac{1}{2}(1,1,0) + \frac{1}{2}(1,0,1)\in \B$, and thus, for any $\epsilon\in (0,1)$, the point
\begin{equation*}
\bigl(1-\varepsilon,\,\tfrac{1}{2},\,\tfrac{1}{2}\bigr) = (1-\epsilon)\bigl(1,\tfrac{1}{2},\tfrac{1}{2}\bigr) + \epsilon\bigl(0,\tfrac{1}{2},\tfrac{1}{2}\bigr)
\end{equation*}
is an interior point of $\B$ as a (non-trivial) convex combination of two points of $\B$ with one of them being interior.
\end{proof}

\noindent\textit{Comment on the proof of Proposition \ref{prop:R^3-exactly-two-pairs}.} As mentioned above, we can assume that $\B$ contains the points $(1,1,0)$ and $(1,0,1)$ but not the point $(0,1,1)$.

\smallskip

We will distinguish two main cases: 
\begin{itemize}
\item[$\blacklozenge$] either both $(1,1,0)$ and $(1,0,1)$ are extreme points of $\B$,
\item[$\blacklozenge$] or at least one of them is not an extreme point of $\B$.
\end{itemize}
In the former case, we will see that the set ${\cal F}^2_{\ref{prop:R^3-exactly-two-pairs},\epsilon_2,\delta_{\epsilon_2}}$ (which coincides with the set $\I_{ex2,2}^n(\epsilon_2,\delta_{\epsilon_2})$ of Proposition \ref{prop:exactly-two-n-1-tuples} when $n=3$) illuminates $\B$ for some explicit $\epsilon_2,\delta_{\epsilon_2}$ depending on $\B$.

\medskip

On the other hand, if e.g. $(1,1,0)$ is NOT an extreme point of $\B$, then necessarily we can find $a\in (0,1)$ such that $(1,1,a)\in \B$. In such a case we can show that the set ${\cal F}^1_{\ref{prop:R^3-exactly-two-pairs},\delta_1}$ illuminates $\B$ for some explicit $\delta_1$.

Similarly, if $(1,1,0)$ is an extreme point of $\B$, but $(1,0,1)$ is not, then (assuming what is claimed in the previous paragraph) it is not hard to deduce that a coordinate permutation of ${\cal F}^1_{\ref{prop:R^3-exactly-two-pairs},\delta_1}$ illuminates $\B$: indeed, it suffices to consider the transformation that swaps the 2nd and the 3rd coordinate.

\medskip

Full details can be found in the proof of Proposition \ref{prop:exactly-two-n-1-tuples}, which is the generalisation of Proposition \ref{prop:R^3-exactly-two-pairs} to all dimensions $n\geq 3$. \qed

\bigskip

The proof of Theorem \ref{thm:R^3-summary} will be completed with the proofs of Propositions \ref{prop:exactly-two-n-1-tuples} and \ref{prop:all-n-1-tuples} in 
the next section.

\section{Bodies with maximal unit subcubes}\label{sec:maximal-unit-subcubes}

In this section we deal with 1-unconditional convex bodies $\B$ in arbitrary dimensions $n$ which:
\begin{itemize}
\item have been normalised to be in $\U^n$, thus are contained in the unit cube $[-1,1]^n$;
\item are not affine images of the cube, and thus certainly satisfy $\B\subsetneq [-1,1]^n$;
\item and have at least one coordinate hyperplane projection (equivalently, coordinate hyperplane section) equal to $[-1,1]^{n-1}$ (this is equivalent to saying that $\B$ contains the point ${\bm 1}-e_i$ for at least one $i\in [n]$). As mentioned before, we will also say in this case that $\B$ contains a \textit{maximal unit subcube}.
\end{itemize}
Classifying 1-unconditional convex bodies in this way, by whether points of the form $e_{i_1}+e_{i_2}+\cdots + e_{i_m}$ are contained in $\B$ or not, for different $m\in \{0,1,\ldots,n-1\}$, is inspired by recent approaches to settle the Illumination Conjecture for bodies with many symmetries (starting with Tikhomirov's work \cite{Tikhomirov-2017}, and also adopted in the precursor \cite{Sun-Vritsiou-sym}  to this paper). 

Of course, we could also describe the instances that we are focusing on here without/before employing the `special' normalisation that we use. We are considering 1-unconditional convex bodies $\B$ which satisfy the following: if ${\cal R}_\B$ is the circumscribing rectangular box given by 
\[{\cal R}_\B= \bigl\{x\in {\mathbb R}^n: |x_i|\leq \|e_i\|_\B^{-1}\ \hbox{for all}\ i\in [n]\bigr\},\] 
then at least one coordinate hyperplane projection of $\B$ and the corresponding one of ${\cal R}_\B$ coincide.  
We will show that, for such bodies $\B$, $\II(\B) \leq 2^n-2$ unless $\B$ is an affine image of the cube.

The results of this section will also help us obtain the main result of the next section. We consider cases based on how many of the coordinate hyperplane projections of $\B$ coincide with the corresponding ones of $[-1,1]^n$ (or, before normalisation, with those of the circumscribing rectangular box); equivalently, based on how many maximal unit subcubes $\B$ contains.

\begin{proposition}\label{prop:all-n-1-tuples}
Let $n\geq 3$ and let $\B \in \mathcal{U}^n$ with the property that, for all $i\in [n]$, $\norm{\sum_{j\neq i} e_j}_\B = \norm{{\bm 1}-e_i}_\B = 1$ but $\B\neq [-1,1]^n$. Then there is $\delta=\delta_\B > 0$ such that $\B$ can be illuminated by the set
\begin{equation*}
\I^n_{-2}(\delta)=\I^n(\delta) \backslash\bigl\{\pm(+\delta, +\delta,\ldots,+\delta,+\delta, -\!\!\!\!\!\!-\delta, +1)\bigr\}.
\end{equation*}
In other words, $\II(\B) \leq 2^n -2$.
\end{proposition}
\begin{proof} Set $\gamma = \norm{{\bm 1}}_\B^{-1}$; then $\gamma{\bm 1}\in \partial\B$ (clearly $\gamma < 1$ since $\B\subsetneq [-1,1]^n$). We will see that $\B$ is illuminated by the set in the statement as long as $\delta < 1-\gamma$.

\smallskip

Observe that, for every boundary point $y = (y_1,y_2,\ldots, y_{n-1}, y_n)$ of $\B$, we can pick two different indices $i_y, j_y\in [n]$ such that $|y_{i_y}|\leq \gamma\leq |y_{j_y}|$ (this is because, if we had $|y_i| > \gamma$ for all $i\in [n]$, then by Lemma \ref{lem:smaller-coordinates} we would get that $\gamma{\bm 1}\in \intr(\B)$, which would contradict our choice of $\gamma$; moreover, if the first of the desired inequalities were satisfied by every coordinate of $y$, that is, if we had $|y_i|\leq \gamma$ for all $i\in [n]$, then for at least one index $j$ we should have $|y_j|=\gamma$, otherwise $y$ would not be a boundary point of $\B$).

\medskip

Fix now some $y\in\partial\B$, and pick the smallest, say, index $i_0$ such that $|y_{i_0}| \leq \gamma$. As recalled in Subsection \ref{subsec:sym-results}, Fact B, we can find a direction $d= d_{y,i_0}\in \I^n_{-2}(\delta)$ which deep illuminates ${\rm Proj}_{e_{i_0}^\perp}(y)$ (since the latter vector is \underline{non-zero} and has at most $n-1$ non-zero coordinates). 

\smallskip

Note that the direction $d=d_{y,i_0}$ that we just considered illuminates the point $y$. Indeed, if $t_0$ is the index at which $\norm{d}_\infty$ is attained (by the definition of deep illumination, this also ensures that $y_{t_0}\neq 0$), then 
\begin{itemize}
\item $(y+|y_{t_0}|d)_{t_0} = 0$, 
\item and at the same time $|(y+|y_{t_0}|d)_i| \leq \max\bigl(|y_i|-|y_{t_0}|\delta,\, |y_{t_0}|\delta\bigr) < 1$ for all $i\in [n]\backslash\{t_0,i_0\}$ (here this is satisfied as long as $\delta < 1$). 
\item Finally, $|(y+|y_{t_0}|d)_{i_0}| \leq |y_{i_0}| + |y_{t_0}|\delta \leq |y_{i_0}| + \delta < 1$, as long as we choose $\delta < 1-\gamma$. 
\end{itemize}
We can now invoke Lemma \ref{lem:affine-set} for the affine subspace $\{\xi\in \R^n: \xi_{t_0} = 0\}$: we compare $y+|y_{t_0}|d$ with ${\bm 1}-e_{t_0}$ (both points of $\B\cap \{\xi\in \R^n: \xi_{t_0} = 0\}$) to conclude that the former point is in the interior of $\B$.
\end{proof}

\begin{proposition}\label{prop:all-but-one-n-1-tuples}
Let $\bm{n\geq 4}$ and let $\B \in \mathcal{U}^n$ with the property that there is $s_0\in [n]$ such that $\norm{{\bm 1} -e_i}_\B = 1$ for all $i\neq s_0$, while at the same time $\norm{{\bm 1}-e_{s_0}}_\B > 1$. Then there is $\delta=\delta_\B > 0$ such that $\B$ can be illuminated by coordinate permutations of the set
\begin{equation*}
\I^{n-1}_{-2}(\delta)\times \{-\delta,+\delta\}.
\end{equation*}
In other words, $\II(\B) \leq 2\cdot (2^{n-1}-2) = 2^n - 4$.
\end{proposition}
\noindent{\bf Note.} Given our notation,
\begin{equation*}
\I^{n-1}_{-2}(\delta) = \I^{n-1}(\delta) \backslash\left\{\pm\left(e_{n-1} -\delta e_{n-2} + \delta\sum_{j\in [n-3]}e_j\right)\right\}
\end{equation*}
(where the standard basis vectors here are considered in $\R^{n-1}$).
\begin{proof}
WLOG assume that $s_0 =n$. Set $\theta_n = \norm{{\bm 1}-e_n}_\B^{-1}$. Then $\theta_n({\bm 1}-e_n)\in \partial \B$, and we have $0 < \theta_n < 1$ by our assumptions. We will see that $\B$ is illuminated by the set in the statement as long as $\delta < 1-\theta_n$ (and by coordinate permutations of this set if $s_0\neq n$).

Consider a boundary point $y$ of $\B$, and moreover suppose that $y$ is an extreme boundary point (recall Fact A from Section \ref{sec:prelims}, namely that it suffices to show how to illuminate these boundary points). Because $n\geq 4$, because $y$ is an extreme boundary point, and given that, by our assumptions, ${\bm 1}-e_i\in \partial \B$ for all $i\in [n-1]$, we can find two distinct indices $i_y, j_y\in [n-1]$ such that $y_{i_y}\cdot y_{j_y}\neq 0$. Moreover, since $\theta_n({\bm 1}-e_n)\in \partial \B$, and by the 1-unconditionality of $\B$, we can find an index $i_0\in [n-1]$ such that $y_{i_0}\leq \theta_n$. 

Note that at least one of the two indices $i_y, j_y$ is different from $i_0$, and thus the vector ${\rm Proj}_{[e_{i_0},e_n]^\perp}(y)$ is non-zero. At the same time, it has at most $n-2$ non-zero coordinates among its first $n-1$. Therefore, looking initially at these coordinates, we can find a direction $d^\prime=d^\prime_{y,i_0}\in \I^{n-1}_{-2}(\delta)$ which deep illuminates the subvector of the first $n-1$ coordinates of ${\rm Proj}_{e_{i_0}^\perp}(y)$ (recall Fact B from Subsection \ref{subsec:sym-results}), and then we can naturally rely on this to pick a direction $d=d_{y,i_0}\in \I^{n-1}_{-2}(\delta)\times \{-\delta,+\delta\}$ which deep illuminates ${\rm Proj}_{e_{i_0}^\perp}(y)$.

Given the way we selected $d=d_{y,i_0}$, if $t_0$ is the index at which $\|d\|_\infty$ is attained, then $t_0\in [n-1]\setminus\{i_0\}$. Moreover, $y_{t_0}\neq 0$. We can now check, just as in the previous proposition, that $y+|y_{t_0}|d\in \intr\B$ (by comparing coordinate-wise the displaced vector $y+|y_{t_0}|d$ with the boundary point ${\bm 1}-e_{t_0}$).
\end{proof}

\begin{remark}
\textup{The above proposition cannot be restated as simply as above when $n=3$ because in $\R^3$ there are $\widetilde{Q}\in \U^3$ which are affine images of the cube and satisfy the main assumption (namely that there exists $s_0\in [n]$ such that $\norm{{\bm 1}-e_i}_{\widetilde{Q}} = 1$ for all $i\neq s_0$, while at the same time $\norm{{\bm 1}-e_{s_0}}_{\widetilde{Q}} > 1$; this is the case when e.g. $\widetilde{Q} = CP_1^2\times [-1,1]$). Therefore, when $n=3$, we need to further assume that $\B$ is not an affine image of the cube, and we need to have extra steps in our proof which essentially encode and capitalise on this additional, necessary assumption. This leads to Proposition \ref{prop:R^3-exactly-two-pairs}, which finds a better high-dimensional analogue in Proposition \ref{prop:exactly-two-n-1-tuples}.}
\end{remark}

For the remaining cases, we need to introduce some further, combinatorially constructed, sets of directions in ${\mathbb R}^n$ that will serve as building blocks for the illuminating sets we will use.

\subsection{Constructing other illuminating sets}\label{subsec:illum-sets}

\begin{notation} Let us fix $n\geq 2$ and $\delta \in (0,1)$, and consider the set $\I^n(\delta)$ from Reminder \ref{reminder:In(delta)} (exactly as is described there).

Define a function $m.c.$ which maps each $d\in \I^n(\delta)$ to the index of its maximum (in absolute value) coordinate; in other words $m.c. : \I^n(\delta)\to [n]$ and e.g. $m.c.((-\delta,\delta,-1,-\delta,...,-\delta,-\delta)) = 3$. Note that the index of the maximum coordinate, as well as its sign, will not change no matter what value of $\delta\in (0,1)$ we pick. Thus, we can also identify directions $d_1\in \I^n(\delta_1)$ and $d_2\in \I^n(\delta_2)$ if their respective coordinates have the same signs and if their maximum coordinate is the same, and then we can view $m.c.$ as a function from the set of these equivalence classes to $[n]$.
\medskip\\
By abusing our notation, we also consider the set $\I^n(1)$, which is simply the set $\{-1,1\}^n$. For each $d^\prime\in \I^n(1)$, there is a unique direction $d_{d^\prime}\in \I^n(1/2)$ (say) which agrees in sign in each coordinate. Then we can also define a function $m.c. : \I^n(1) \equiv \{-1,1\}^n \to [n]$ by setting $m.c.(d^\prime) = m.c.(d_{d^\prime})$.
\end{notation}

\begin{notation} 
Starting from the set $\I^n(\delta)$, we construct a new, similar set $\widehat{\I}^n_{n-1,n}(\delta)$ in the following way: 
\begin{itemize}
\item if $d\in \I^n(\delta)$ satisfies $m.c.(d) \in \{n-1,n\}$, then we keep it in $\widehat{\I}^n_{n-1,n}(\delta)$ as well (there are $2^{n-1} + 2^{n-2}$ such directions). 
\item If $d\in \I^n(\delta)$ satisfies $m.c.(d) < n-1$, then replace $d$ by a direction $d^\prime$ which agrees in sign with $d$ in every coordinate, but has maximum (in absolute value) coordinate the $(n-1)$-th one, and place $d^\prime$ in $\widehat{\I}^n_{n-1,n}(\delta)$. Thus e.g. the direction $(-1,-\delta,-\delta,\ldots,-\delta,-\delta,-\delta)$ from $\I^n(\delta)$ will give the direction $(-\delta,-\delta,-\delta,\ldots,-\delta,-1,-\delta)$ in $\widehat{\I}^n_{n-1,n}(\delta)$.
\end{itemize}

\noindent Similarly, we construct a new set $\widehat{\I}^n_{n-2,n-1,n}(\delta)$ in the following way: 
\begin{itemize}
\item if $d\in \I^n(\delta)$ satisfies $m.c.(d) \in \{n-2,n-1,n\}$, then we keep it in $\widehat{\I}^n_{n-2,n-1,n}(\delta)$ as well (these are $2^{n-1} + 2^{n-2}+2^{n-3}$ directions). 
\item If $d\in \I^n(\delta)$ satisfies $m.c.(d) < n-2$, then replace $d$ by a direction $d^\prime$ which agrees in sign with $d$ in every coordinate, but has maximum (in absolute value) coordinate the $(n-2)$-th one, and place $d^\prime$ in $\widehat{\I}^n_{n-2,n-1,n}(\delta)$. Thus e.g. the direction $(-1,-\delta,-\delta,\ldots,-\delta,-\delta,-\delta)$ from $\I^n(\delta)$ will give the direction $(-\delta,-\delta,-\delta,\ldots,-1,-\delta,-\delta)$ in $\widehat{\I}^n_{n-2,n-1,n}(\delta)$.
\end{itemize}
\end{notation}

\begin{definition}
Let $n\geq 3$ be an integer. If $n$ is even, we set
\begin{equation*}
\J_n:=
\bigl\{d\in \I^n(1):\,\#\{i\in[n]: d_i=+1\}\in \{2,4,6\ldots,n-2\}\bigr\}.
\end{equation*}

If $n$ is odd, we distinguish cases.
\begin{itemize}
\item[--] If $n=3$, we set $\J_3=\{(1,1,1),\,(-1,-1,-1)\}$. 
\item[--] If $n>3$ and $(n-1)/2$ is odd, we first define
\begin{equation*}
\J_n^1:= 
\bigl\{d\in \I^n(1):\,\#\{i\in[n]: d_i=+1\}\in \{2,4,\ldots, \tfrac{n-1}{2}-1\}\bigr\}
\end{equation*}
and then
\begin{equation*}
\J_n : = \J_n^1 \cup (-\J_n^1).
\end{equation*}
\item[--] Finally, if $n\geq 3$ is odd and $(n-1)/2$ is even, then we set
\begin{equation*}
\J_n^1:= 
\bigl\{d\in \I^n(1):\,\#\{i\in[n]: d_i=+1\}\in \{1,3,\ldots, \tfrac{n-1}{2}-1\}\bigr\}
\end{equation*}
and afterwards we set
\begin{equation*}
\J_n : = \J_n^1 \cup (-\J_n^1).
\end{equation*}
\end{itemize}
\end{definition}

\begin{remarks}\label{observ:J_n-sets}
(a) With the above definition, we have ensured that the sets $\J_n$ are symmetric (that is, $\J_n = -\J_n$).

\medskip

(b) We also need some `efficient' bounds on the cardinalities of $\J_n$. Note that $\abs{\J_3}=2 = 2^{n-1}-2$. We will check that, for all $n\geq 3$, $\abs{\J_n}\leq 2^{n-1} -2$. 

\smallskip

When $n$ is even, observe that, to determine a direction in $\J_n$, we simply need to know which coordinates are equal to $+1$, and the subset of the corresponding indices will range over all subsets of $[n]$ of even cardinality $\geq 2$ and $\leq n-2$. Thus
\begin{equation*}
\abs{\J_n} = \sum_{s=1}^{(n-2)/2}\binom{n}{2s} = \frac{2^n}{2} - \binom{n}{0}-\binom{n}{n} = 2^{n-1}-2
\end{equation*}
(the fact that $\displaystyle \sum_{s=0}^{\lfloor m/2\rfloor}\binom{m}{2s} = \sum_{u\,\hbox{\footnotesize even}\,\leq m}\binom{m}{u}=\frac{2^m}{2}$, regardless of whether $m$ is even or odd, can be checked using induction in $m$).

\smallskip

To estimate the cardinality of $\J_n$ when both $n>3$ and $(n-1)/2$ are odd, we first observe that
\begin{equation*}
\sum_{u=0}^{(n-1)/2}\binom{n}{u} = 2^{n-1}
\end{equation*}
and also that, for $u < (n-1)/2$, we have that $\binom{n}{u} < \binom{n}{u+1}$. Therefore,
\begin{align*}
2^{n-1}= \binom{n}{0} +\binom{n}{1} + \sum_{u=2}^{(n-1)/2}\binom{n}{u}
\geq \,n+1\ +\  \frac{n-3}{4} + 2\sum_{s=1}^{(n-3)/4}\binom{n}{2s},
\end{align*}
which shows that
\begin{equation*}
\abs{\J_n^1} = \sum_{s=1}^{(n-3)/4}\binom{n}{2s} \leq 2^{n-2} - \frac{5n+1}{8}.
\end{equation*}
This gives that
\begin{equation*}
\abs{\J_n}= \sum_{s=1}^{(n-3)/4}\left[\binom{n}{2s}+\binom{n}{n-2s}\right] = 2\sum_{s=1}^{(n-3)/4}\binom{n}{2s} \leq 2^{n-1} - \frac{5n+1}{4}.
\end{equation*}

Finally, when $n>3$ is odd and $(n-1)/2$ is even, we similarly observe that
\begin{align*}
2^{n-1}= \binom{n}{0} + \sum_{u=1}^{(n-1)/2}\binom{n}{u}
\geq \,1\  +\ \frac{n-1}{4} + 2\sum_{s=1}^{(n-1)/4}\binom{n}{2s-1} 
\end{align*}
which shows that
\begin{equation*}
\abs{\J_n^1} = \sum_{s=1}^{(n-1)/4}\binom{n}{2s-1} \leq 2^{n-2} - \frac{n+3}{8}.
\end{equation*}
Thus, in this case
\begin{equation*}
\abs{\J_n} \leq 2^{n-1} - \frac{n+3}{4}.
\end{equation*}
As claimed, in all cases we have that $\abs{\J_n}\leq 2^{n-1} -2$.

\medskip

(c) The last key property of the sets $\J_n$ which we use in the sequel is the following: consider any direction $d^\prime
\in \{0,1,-1\}^n$ such that exactly two coordinates of $d^\prime$ are equal to $0$. Then we can find a direction $d\in \J_n$ such that $d_i=d^\prime_i$ whenever $d^\prime_i\neq 0$.

\smallskip

Indeed, let $i_1,i_2\in [n]$ be the two indices for which $d^\prime_{i_1}=d^\prime_{i_2}=0$ (WLOG suppose that $i_1<i_2$). Assume first that $n$ is even. The subset of indices $i$ for which $d^\prime_i=+1$ is a subset ${\cal P}_{d^\prime}$ of $[n]\backslash\{i_1,i_2\}$. If $\abs{{\cal P}_{d^\prime}}=0$, set $d_{i_1}=d_{i_2}=+1$, and we will have that the corresponding subset ${\cal P}_d$ for $d$ has cardinality 2, so $d\in \J_n$. If $\abs{{\cal P}_{d^\prime}}$ is even and $>0$, then we can set $d_{i_1}=d_{i_2}=-1$, and we will have that $0<\abs{{\cal P}_d}\leq n-2$, thus $d$ will be in $\J_n$ again. Finally, if $\abs{{\cal P}_{d^\prime}}$ is odd, then it must be an odd number between 1 and $n-3$, thus we can set $d_{i_1}=+1$, $d_{i_2}=-1$.

\medskip

Assume now that $n$ is odd. We first deal with the case $n=3$. In this case, $\abs{{\cal P}_{d^\prime}}=0$ or $=1$. If the former holds, set $d=(-1,-1,-1)$, while, if the latter holds, set $d=(+1,+1,+1)$. In both cases, $d$ will agree with $d^\prime$ in the unique entry of  $d^\prime$ that is non-zero.

\smallskip

If $n>3$ and $(n-1)/2$ is odd, again observe that ${\cal P}_{d^\prime}$ is a subset of $[n]\backslash\{i_1,i_2\}$. If $\abs{{\cal P}_{d^\prime}}=0$, set $d_{i_1}=d_{i_2}=+1$, and we will have that $\abs{{\cal P}_d}=2$, thus $d\in \J_n$. If $\abs{{\cal P}_{d^\prime}}$ is a positive even number $< (n-1)/2$, then we set $d_{i_1}=d_{i_2}=-1$. If $\abs{{\cal P}_{d^\prime}}$ is a positive even number $> (n-1)/2$, then it is also $<n-2$ (since $n-2$ is odd in this case), and we can set $d_{i_1}= +1$, $d_{i_2}=-1$, in which case $\abs{{\cal P}_d}$ will be an odd number between $\frac{n+3}{2}=n-(\frac{n-1}{2}-1)$ and $n-2$, thus $d\in \J_n$. If instead $\abs{{\cal P}_{d^\prime}}$ is an odd number $< (n-1)/2$, then set $d_{i_1}= +1$, $d_{i_2}=-1$, in which case $\abs{{\cal P}_d}$ will be an even number between $2$ and $\frac{n-3}{2}=\frac{n-1}{2}-1$. If $\abs{{\cal P}_{d^\prime}}$ is an odd number $\geq (n-1)/2$ and $< n-2$, set $d_{i_1}=d_{i_2}=+1$, in which case $\abs{{\cal P}_d}$ will be an odd number between $\frac{n+3}{2}$ and $n-2$. Finally, if $\abs{{\cal P}_{d^\prime}}=n-2$, set $d_{i_1}=d_{i_2}=-1$. In all these cases, we end up with a direction $d\in \J_n$.

\smallskip

Similarly we deal with the last case, where $n > 3$ is odd and $(n-1)/2$ is even. If $\abs{{\cal P}_{d^\prime}}$ is an even number $< (n-1)/2$, then set $d_{i_1}= +1$, $d_{i_2}=-1$. If $\abs{{\cal P}_{d^\prime}}$ is an even number $\geq (n-1)/2$, then it will also be $<n-2$ (since $n-2$ is odd), and thus we will be able to set $d_{i_1}=d_{i_2}=+1$ to get that $\abs{{\cal P}_d}$ is an even number $\geq \frac{n+3}{2}=n-(\frac{n-1}{2}-1)$ and $\leq n-1$. If $\abs{{\cal P}_{d^\prime}}$ is an odd number $< (n-1)/2$, then set $d_{i_1}=d_{i_2}=-1$. If $\abs{{\cal P}_{d^\prime}}$ is an odd number $> (n-1)/2$ (and obviously $\leq n-2$ since ${\cal P}_{d^\prime}\subseteq [n]\backslash\{i_1,i_2\}$), then set $d_{i_1}=+1$ and $d_{i_2}=-1$. In all cases, we pick a direction $d\in \J_n$ as needed.
\end{remarks}

\smallskip

With these properties at hand, we are now ready to illuminate the remaining cases of bodies in $\U^n$ which contain a maximal unit subcube.

\subsection{Remaining cases with maximal unit subcubes}

\begin{proposition}\label{prop:exactly-n-k-n-1-tuples}
Let $k\geq 2$ and $n\geq k+3$, and consider $\B\in {\cal U}^n$ such that we can find $k$ indices $j_1 < j_2 < \cdots < j_k$ in $[n]$ with the property that 
\begin{equation*}
{\bm 1} - e_j \in \B
\end{equation*}
for all $j\in [n]\backslash\{j_1,j_2,\ldots,j_k\}$, while, if $j_s\in \{j_1,j_2,\ldots,j_k\}$,
\begin{equation*}
{\bm 1} - e_{j_s} \notin \B.
\end{equation*}
Then there is $\delta=\delta_\B$ such that: 
\begin{itemize}
\item if $k=2$, then $\B$ can be illuminated by a coordinate permutation of the set
\begin{gather}\label{eq:prop:exactly-n-k-n-1-tuples:illumset1}
\bigl[[\I^{n-2}(\delta)\backslash \{\pm\underbrace{(\delta,\delta,\ldots,\delta,-\delta,1)}_{n-2}\}]\times\{\delta,-\delta\}^2\bigr]\ \bigcup\  \bigl\{\pm(\delta,\delta,\ldots,\delta,-\delta,1,0,0)\bigr\},
\end{gather}
\item
and if $k\geq 3$, then $\B$ can be illuminated by a coordinate permutation of the set
\begin{gather}\label{eq:prop:exactly-n-k-n-1-tuples:illumset2}
\bigl[[\I^{n-k}(\delta)\backslash \{\pm\underbrace{(\delta,\delta,\ldots,\delta,-\delta,1)}_{n-k}\}]\times\{\delta,-\delta\}^k\bigr]
\ \bigcup\  \bigl[\{\pm\underbrace{(\delta,\delta,\ldots,\delta,-\delta,1)}_{n-k}\}\times \delta\cdot\J_k\bigr].
\end{gather}
\end{itemize}
Thus $\II(\B) \leq (2^{n-2}-2)\cdot 4 + 2 = 2^n-6$ if $k=2$, while $\II(B)\leq (2^{n-k}-2)\cdot 2^k + 2\cdot (2^{k-1}-2) = 2^n - 2^k - 4$ if $k\geq 3$.
\end{proposition}
\begin{proof} \textit{Case $k=2$:} WLOG suppose that $\{j_1,j_2\}=\{n-1,n\}$. In other words, ${\bm 1}-e_{n-1}$ and ${\bm 1}-e_n$ are not in $\B$, while ${\bm 1}-e_j\in \partial \B$ for all $j\in [n-2]$.

Set $\theta_{n-1}=\|{\bm 1}-e_{n-1}\|_\B^{-1}$ and $\theta_n=\|{\bm 1}-e_n\|_\B^{-1}$; by our assumptions $\theta_{n-1},\theta_n\in (0,1)$. Set also $\Theta_0=\max\{\theta_{n-1},\theta_n\}$. We will show that $\B$ is illuminated by the set in \eqref{eq:prop:exactly-n-k-n-1-tuples:illumset1} as long as $\delta < 1-\Theta_0$.

Consider a boundary point $y$ of $\B$, and moreover suppose that $y$ is an extreme boundary point. Because $n\geq k+3\geq 5$, because $y$ is an extreme boundary point, and given that, by our assumptions, ${\bm 1}-e_i\in \partial \B$ for all $i\in [n-2]$, we can find two distinct indices $i_y, j_y\in [n-2]$ such that $y_{i_y}\cdot y_{j_y}\neq 0$. Moreover, since $\theta_n({\bm 1}-e_n)\in \partial \B$, and by the 1-unconditionality of $\B$, we can find an index $i_n\in [n-1]$ such that $|y_{i_n}|\leq \theta_n\leq \Theta_0$. Similarly, because $\theta_{n-1}({\bm 1}-e_{n-1})\in \partial \B$, we can find an index $i_{n-1}\in [n]\setminus\{n-1\}$ such that $|y_{i_{n-1}}|\leq \theta_{n-1}\leq \Theta_0$ (note that it's possible that $i_{n-1}=i_n$).

\medskip

Assume first that either $i_n$ or $i_{n-1}$ can be picked from $[n-2]$, and denote the smallest such index by $i_0$. 
Note that at least one of the two indices $i_y, j_y$ is different from $i_0$, and thus the vector ${\rm Proj}_{[e_{i_0},e_{n-1}, e_n]^\perp}(y)$ is non-zero. At the same time, it has at most $n-3$ non-zero coordinates among its first $n-2$. Therefore, looking initially at these coordinates, we can find a direction 
\begin{equation*}
d^\prime=d^\prime_{y,i_0}\in \I^{n-2}_{-2}(\delta)=\I^{n-2}(\delta)\backslash \{\pm\underbrace{(\delta,\delta,\ldots,\delta,-\delta,1)}_{n-2}\}
\end{equation*} which deep illuminates the subvector of the first $n-2$ coordinates of ${\rm Proj}_{e_{i_0}^\perp}(y)$ (recall Fact B), and then we can rely on this to pick a direction $d=d_{y,i_0}\in \I^{n-2}_{-2}(\delta)\times \{-\delta,+\delta\}^2$ which deep illuminates ${\rm Proj}_{e_{i_0}^\perp}(y)$.

Given the way we selected $d=d_{y,i_0}$, if $t_0$ is the index at which $\|d\|_\infty$ is attained, in other words, if $t_0=m.c.(d)$, then $t_0\in [n-2]\setminus\{i_0\}$. Moreover, $y_{t_0}\neq 0$. We can now check that $y+|y_{t_0}|d\in \intr\B$ (by comparing coordinate-wise the displaced vector $y+|y_{t_0}|d$ with the boundary point ${\bm 1}-e_{t_0}$).

\medskip

Next assume that $|y_i| > \Theta_0$ for all $i\in [n-2]$. Then necessarily $i_n=n-1$ and $i_{n-1}=n$. In other words, $\max\{|y_{n-1}|,\,|y_n|\}\leq \Theta_0 < 1$. We can then pick a direction $d$ from the set in \eqref{eq:prop:exactly-n-k-n-1-tuples:illumset1} such that $d_i\cdot y_i< 0$ for all $i\in [n-2]$, and, as before, we can check that $y+|y_{t_0}|d\in \intr\B$ (where $t_0=m.c.(d)$).

\medskip

\textit{Cases where $k\geq 3$:} WLOG we suppose that $\{j_1,j_2,\ldots,j_k\}= \{n-k+1,n-k+2,\ldots,n\}$, and for each $j\in \{n-k+1,n-k+2,\ldots,n\}$ we set $\theta_j=\norm{{\bm 1}-e_j}_\B^{-1}$. We also set $\Theta_0=\max\{\theta_j: n-k+1\leq j\leq n\}$; by our assumptions $\Theta_0\in (0,1)$. We will show that $\B$ is illuminated by the set in \eqref{eq:prop:exactly-n-k-n-1-tuples:illumset2} as long as $\delta < 1-\Theta_0$.

Again consider an extreme boundary point $y$ of $\B$. We can find two distinct indices $i_y, j_y\in [n-k]$ such that $y_{i_y}\cdot y_{j_y}\neq 0$. If we can also find an index $i_0\in [n-k]$ such that $|y_{i_0}|\leq \Theta_0$, then as before we can pick a direction
\begin{equation*}
d\in [\I^{n-k}(\delta)\backslash \{\pm\underbrace{(\delta,\delta,\ldots,\delta,-\delta,1)}_{n-k}\}]\times\{\delta,-\delta\}^k
\end{equation*}
which will deep illuminate ${\rm Proj}_{e_{i_0}^\perp}(y)$, and then we can check that $y+|y_{t_0}|d\in \intr\B$ (where $t_0=m.c.(d)$).

\smallskip

Assume now that $|y_i| > \Theta_0$ for all $i\in [n-k]$. Since $\theta_n({\bm 1}-e_n)\in \partial\B$, there is an index $i_n\in [n-1]$ such that $|y_{i_n}|\leq\theta_n\leq \Theta_0$. Given our prior assumption, we also have that $i_n\in \{n-k+1,\ldots, n-1\}$. Next we also use the assumption that $\theta_{i_n}({\bm 1}-e_{i_n})\in \partial \B$, which implies that there is an index
\begin{equation*}
s_n\in \{n-k+1,n-k+2,\ldots,n\}\setminus\{i_n\}
\end{equation*}
such that $|y_{s_n}|\leq \theta_{i_n}\leq \Theta_0$. 

We can find a (unique) direction $d^\prime\in \I^{n-k}(\delta)$ such that $d^\prime_i\cdot y_i < 0$ for all $i\in [n-k]$. 
\begin{itemize}
\item[$\blacklozenge$] If 
\begin{equation*}
d^\prime\notin \{\pm(\delta,\delta,\ldots,\delta,-\delta,1)\},
\end{equation*}
then we can extend $d^\prime$ to a direction
\begin{equation*}
d\in [\I^{n-k}(\delta)\backslash \{\pm\underbrace{(\delta,\delta,\ldots,\delta,-\delta,1)}_{n-k}\}]\times\{\delta,-\delta\}^k
\end{equation*}
which deep illuminates $y$.
\item[$\blacklozenge$] On the other hand, if $d^\prime\in \{\pm(\delta,\delta,\ldots,\delta,-\delta,1)\}$, then, given the third main property of the set $\J_k$, we can find a direction 
\begin{equation*}
d\in \{\pm(\delta,\delta,\ldots,\delta,-\delta,1)\}\times \delta\cdot\J_k
\end{equation*}
which deep illuminates the vector ${\rm Proj}_{[e_{i_n},e_{s_n}]^\perp}(y)$. 
\end{itemize}
In both subcases, we can check that, for the direction $d$ that we ended up picking, $y+|y_{t_0}|d\in \intr\B$ (where $t_0=m.c.(d)\in [n-k]$).
The proof is complete.
\end{proof}

It is clear from the above proof that the assumption that $n-k\geq 3$ (that is, the number of maximal unit subcubes contained in $\B$ is $\geq 3$) was crucially used (so that, given an extreme boundary point $y$, we could find two distinct indices $i_y,j_y\in [n-k]$ such that $y_{i_y}\cdot y_{j_y}\neq 0$). This is not an artefact of the proof though. As we will see, if the number of maximal unit subcubes contained in $\B$ is exactly 2, then $\B$ could also be an affine image of the cube (e.g. $\B=CP_1^2\times[-1,1]^{n-2}$). Thus we need to argue more carefully about how to illuminate such $\B$ which are not parallelepipeds. 

Not for the same reason, but we also have to argue more carefully when $\B$ contains only one maximal unit subcube; this is the case that we deal with now.


\begin{proposition}\label{prop:exactly-one-n-1-tuple}
Let $n\geq 4$, and let $\B\in {\cal U}^n$ with the property that there exists 
exactly one index $i_0\in [n]$ such that ${\bm 1}-e_{i_0}\in \B$.
Then there are $\delta, \eta_\delta$ and $\widetilde{\delta}>0$ such that $\B$ can be illuminated by a coordinate permutation of the set
\begin{equation*}
\I_{\delta,\eta,\widetilde{\delta}} = \I_{\delta, \eta} \cup \bigl[\bigl(\widetilde{\delta}\cdot\J_{n-1}\bigr)\times\{1,-1\}\bigr]
\end{equation*}
where $\J_{n-1}$ is defined as in the previous subsection (note that $n-1\geq 3$ here), and where
\begin{equation*}
\I_{\delta, \eta}:=\bigl\{\bigl(\pm(1,\eta_\delta),\pm\delta,\pm\delta,\ldots,\pm\delta,0\bigr),\,\bigl(\pm(-\eta_\delta,1),\pm\delta,\pm\delta,\ldots,\pm\delta,0\bigr)\bigl\}.
\end{equation*}
{\bf Note.} Recall that $\abs{\J_{n-1}} \leq 2^{n-2}-2$ and thus $\abs{\I_{\delta,\eta,\widetilde{\delta}}} \leq 2^{n-1} + 2(2^{n-2}-2) = 2^n-4$. 

Recall also that, since $\J_{n-1}$ is symmetric, the set $\bigl(\widetilde{\delta}\cdot\J_{n-1}\bigr)\times\{1,-1\}$ is formed from pairs of opposite directions, and so is $\I_{\delta,\eta,\widetilde{\delta}}$.
\end{proposition}
\begin{proof}
WLOG we can assume that $i_0=n$, and thus ${\bm 1}-e_n\in \B$, while ${\bm 1}-e_j\notin\B$ for $j\in [n-1]$. For each such $j$, we set $\alpha_j$ to be the supremum of non-negative numbers $x_n$ such that
\begin{equation*}
x_ne_n + \sum_{i\in[n-1]\backslash\{j\}}e_i\in \B.
\end{equation*}
By compactness the point $w_j:=\alpha_je_n + \sum_{i\in[n-1]\backslash\{j\}}e_i\in \B$ (in fact, $\in\partial\B$), and, by our assumptions, $\alpha_j\in [0,1)$. WLOG we assume that $\alpha_1\geq \alpha_2\geq \alpha_j$ for each $j\in [n-1]\backslash\{1,2\}$. 

\smallskip

Given any $\delta \in (0,1)$, the directions $(\pm 1,0,\pm \delta,\pm \delta,\ldots,\pm\delta,0)$ illuminate the point $w_2=\alpha_2e_n+\sum_{i\in [n-1]\backslash\{2\}}e_i$ and all its coordinate reflections (we can use Corollary \ref{cor:affine-set} here). Furthermore, by Lemma \ref{lem:perturb}, we can find $\eta_{0,\delta}>0$ sufficiently small so that the directions \[\bigl(\pm(1,\eta),\pm \delta,\pm \delta,\ldots,\pm\delta,0\bigr)\] will illuminate the same points if $0<\eta \leq \eta_{0,\delta}$ (note that, in the latter subset of (perturbed) directions, we have $\sign(d_1)=\sign(d_2)$ for each direction $d$). 

Similarly, the directions $(0,\pm 1,\pm \delta,\pm \delta,\ldots,\pm\delta,0)$ illuminate $w_1=\alpha_1e_n+\sum_{i\in [n-1]\backslash\{1\}}e_i$ and all its coordinate reflections, and if we pick $\eta_{0,\delta}$ even smaller if needed, so will the directions \[\bigl(\pm(-\eta_{0,\delta},1),\pm\delta,\pm\delta,\ldots,\pm\delta,0\bigr)\] (note that, in the latter subset of directions, $\sign(d_1^\prime)=-\sign(d_2^\prime)$ for each direction $d^\prime$; 
this shows that the set $\I_{\delta,\,\eta_{0,\delta}}$ we just formed contains any combination of signs for the first $n-1$ coordinates).

\medskip

Consider now $j\in [n-1]\backslash\{1,2\}$, and suppose $x$ is one of the coordinate reflections of the point
\begin{equation*}
w_j=\alpha_je_n + \sum_{i\in [n-1]\backslash\{j\}}e_i.
\end{equation*}
Based on whether $\sign(x_1)=\sign(x_2)$ or not, choose 
a direction $d_x$ 
from $\I_{\delta,\,\eta_{0,\delta}}$
such that $d_{x,i}\cdot x_i\leq 0$ for all $i$ \emph{(to avoid any ambiguity later, we can agree to set $d_{x,j} = +\delta$)}. Then the coordinates of $x+d_x$ satisfy the following: 
\begin{itemize}
\item[--] one of $(x+d_x)_1$, $\,(x+d_x)_2$ is equal to 0, 
\item[--] while the other one is $<1$ in absolute value (in fact equal to $1-\eta_{0,\delta}$ in absolute value). 
\item[--] In addition, $\abs{(x+d_x)_j}=\delta$,
\item[--] and for all $i\in [n-1]\backslash\{1,2,j\}$, $\abs{(x+d_x)_i}=1-\delta$. 
\item[--] Finally $\abs{(x+d_x)_n}=\abs{x_n} = \alpha_j\leq \alpha_2\leq \alpha_1< 1$. 
\end{itemize}
Hence we can use Lemma \ref{lem:affine-set} (combined with Lemma \ref{lem:smaller-coordinates}), and compare with one of the points $\alpha_1e_n+\sum_{i\in [n-1]\backslash\{1\}}e_i$ and $\alpha_2e_n+\sum_{i\in [n-1]\backslash\{2\}}e_i$ to conclude that $x+d_x\in \intr\B$.

\bigskip

For the rest of the proof we fix $\delta_0\in (0,1)$ and $\eta_0\equiv \eta_{0,\delta_0}$ such that the set $\I_{\delta_0,\eta_0}$ illuminates all the points $w_j$ and their coordinate reflections. Moreover, for each of these points $x$ we fix a direction $d_x$ from $\I_{\delta_0,\eta_0}$ which illuminates $x$ (we pick it as before, when it is not unique), and then, recalling Lemma \ref{lem:perturb}(b), we also find $\tau_x>0$ with the property that, if $z\in \partial \B$ and $\norm{x-z}_\infty \leq \tau_x$, then $z$ is also illuminated by $d_x$.
Given that the points $\alpha_je_n + \sum_{i\in [n-1]\backslash\{j\}}e_i$ and their coordinate reflections are finitely many, we can set
\begin{equation*}
\tau_0 = \min\left(\frac{1-\alpha_1}{2},\ \min\left\{\tau_x: x\ \hbox{is a coordinate reflection of one of the points $w_j$}\right\}\right).
\end{equation*}

Next, for each $j\in [n-1]$, set $\beta_j=\beta_j(\alpha_j,\tau_0)$ to be the supremum of positive numbers $u$ such that
\begin{equation*}
(\alpha_j+\tau_0)e_n + u\!\!\!\! \sum_{i\in[n-1]\backslash\{j\}}e_i\in \B
\end{equation*}
(such positive numbers exist since we can consider convex combinations of the point $e_n$ and the point $\sum_{i=1}^{n-1}e_i$, and focus on those convex combinations which are closer to $e_n$). By compactness we have that the point $(\alpha_j+\tau_0)e_n + \beta_j\cdot \sum_{i\in[n-1]\backslash\{j\}}e_i\in \B$ (in fact $\in\partial\B$); moreover $\beta_j<1$ since $(\alpha_j+\tau_0)e_n +\sum_{i\in[n-1]\backslash\{j\}}e_i\notin \B$ (by how we chose $\alpha_j$ previously).

Finally, as before, for each $j\in [n-1]$ set $\theta_j:=\norm{{\bm 1}-e_j}_\B^{-1}$; by our main assumption, we have that $\Theta_0:=\max\{\theta_j:j\in [n-1]\} \in (0,1)$. We can finally choose
\begin{equation*}
\widetilde{\delta}_0 < \min\bigl\{\tau_0,\, \min\{1-\beta_j:j\in [n-1]\},\,1-\Theta_0\bigr\}.
\end{equation*}

\medskip

We are now ready to illuminate all boundary points of $\B$ using the set $\I_{\delta_0,\,\eta_0,\,\widetilde{\delta}_0}$. Let $y\in \partial \B$.
\begin{itemize}
\item[$\blacklozenge$] If $y={\bm 1}-e_n$ or one of its coordinate reflections, then $y$ is illuminated by some direction in $\I_{\delta_0,\,\eta_0}$ (since we can find any combination of signs for the first $n-1$ coordinates).
\item[$\blacklozenge$] If $n\in \cZ_y$, then $y$ is contained in the convex hull of ${\bm 1}-e_n$ and its coordinate reflections, and thus it is also illuminated by some direction in $\I_{\delta_0,\,\eta_0}$.
\item[$\blacklozenge$] Assume now that $n\notin \cZ_y$. 
\begin{itemize}
\item[$\bullet$] Suppose also that $y$ has the property:
\begin{equation*}
\hbox{for two distinct indices $s_1,s_2\in [n-1]$, $\,\max\{|y_{s_1}|,\,|y_{s_2}|\} < 1 - \widetilde{\delta}_0$.} \tag{$\ast$}
\end{equation*}
By the third main property of the set $\J_{n-1}$, we can find a direction $d\in \bigl(\widetilde{\delta}_0\cdot\J_{n-1}\bigr)\times\{1,-1\}$ which deep illuminates the vector ${\rm Proj}_{[e_{s_1},e_{s_2}]^\perp}(y)$. We can then check that the displaced vector $y+|y_n|d\in \intr\B$ (by comparing coordinate-wise with the vector ${\bm 1}-e_n$).

Note that property ($\ast$) is satisfied in several instances, including the following:
\begin{itemize}
\item[--] when $\abs{\cZ_y}\geq 2$, given also our previous assumption that $n\notin \cZ_y$.
\item[--] When $|y_n| > \Theta_0$. Indeed, in this case, since $\theta_{n-1}({\bm 1}-e_{n-1})\in \partial\B$, we must be able to find an index $s_1\in [n-2]$ such that $|y_{s_1}| \leq \theta_{n-1}\leq \Theta_0 < 1- \widetilde{\delta}_0$.

Similarly, because $\theta_{s_1}({\bm 1}-e_{s_1})\in \partial\B$, we should be able to find an index $s_2\in [n-1]\setminus\{s_1\}$ such that $|y_{s_2}|\leq \theta_{s_1}\leq \Theta_0$.
\item[--] When $\cZ_y = \{t_0\}\subset [n-1]$, and $|y_n| > \alpha_{t_0} +\tau_0$. Indeed, by our choice of $\beta_{t_0}$, we have that $v_{t_0}:=(\alpha_{t_0}+\tau_0)e_n + \beta_{t_0}\sum_{i\in[n-1]\backslash\{t_0\}}e_i\in \partial\B$. Given that the $n$-th coordinate of ${\rm Proj}_{e_{t_0}^\perp}(y)$ is strictly bigger in absolute value than the $n$-th coordinate of $v_{t_0}$, while their $t_0$-th coordinates are both zero, we must have that some other coordinate of $v_{t_0}$ exceeds the respective one of ${\rm Proj}_{e_{t_0}^\perp}(y)$ in absolute value. Thus, we can find some $t_1\in [n-1]\setminus\{t_0\}$ such that $|y_{t_1}|\leq \beta_{t_0} < 1- \widetilde{\delta}_0$.
\end{itemize}
\smallskip
\item[$\bullet$] Suppose now that, for at least $n-2$ indices $i\in [n-1]$, $|y_i|\geq 1-\widetilde{\delta}_0 > 1-\tau_0$.

First of all, this implies that $|y_n|\leq \Theta_0 < 1$. Therefore,
\begin{itemize}
\item[--] if $\cZ_y=\emptyset$, then, since $|y_n|\in (0,1)$, we can use Corollary \ref{cor:affine-set} to conclude that a direction $d$ from $\I_{\delta_0,\eta_0}$ illuminates $y$ (it suffices to pick the unique direction from $\I_{\delta_0,\eta_0}$ which satisfies $d_i\cdot y_i < 0$ for all $i\in [n-1]$).
\item[--] If instead $\cZ_y\neq \emptyset$, then necessarily, given our main assumptions here, we will have that $\cZ_y=\{t_1\}\subset [n-1]$. From the previous remarks, we know that this implies that $|y_n|\leq \alpha_{t_1}+\tau_0$. 

If $|y_n| \leq \alpha_{t_1}$, then $y$ is in the convex hull of the point 
\begin{equation*}
w_{t_1} = \alpha_{t_1}e_n + \sum_{i\in [n-1]\backslash\{t_1\}}e_i
\end{equation*}
and its coordinate reflections, and thus it is illuminated by some direction in $\I_{\delta_0,\eta_0}$.

If instead $|y_n| \in (\alpha_{t_1},\,\alpha_{t_1}+\tau_0]$, then we will have that $\norm{x-y}_\infty \leq \tau_0 \leq \tau_x$ with $x$ some coordinate reflection of $w_{t_1}$. Thus, $y$ will again be illuminated by some direction $d\in \I_{\delta_0,\eta_0}$ (in fact, the selected direction $d_x$ which illuminates $x$). 
\end{itemize}
\end{itemize}
\end{itemize}
In this way, we have illuminated all boundary points of $\B$.
\end{proof}


We now turn to the case where a convex body $\B\in \U^n$ contains exactly two maximal unit subcubes. The proof of the following proposition will also encompass the full proof of Proposition \ref{prop:R^3-exactly-two-pairs}, that is, the relevant result in $\R^3$.

\begin{proposition}\label{prop:exactly-two-n-1-tuples}
Let $n\geq 3$, and suppose $\B\in {\cal U}^n$ is not an affine image of the cube but has the property that there are exactly two distinct indices $i_1, i_2\in [n]$ such that ${\bm 1}-e_{i_s}\in \B$ for $s=1,2$, while ${\bm 1}-e_j\notin \B$ for any $j\in [n]\backslash\{i_1,i_2\}$. Then, up to coordinate permutations, one of the following sets illuminates $\B$.
\begin{itemize}
\item[1.] The set $\I_{ex2,1}^n(\delta)$, for some $\delta=\delta_\B\in (0,1)$, which consists of the following directions:
\begin{itemize}
\item[-] the directions $\bigl\{\delta\cdot\bigl(d,\,\sign(d_{m.c.(d)})\bigr) : d\in \I^{n-2}(1)\bigr\}\times\{\pm 1\}$ for some $\delta>0$ which depends only on $\B$ (these are $2^{n-2}\cdot 2 = 2^{n-1}$ directions);
\item[-] the directions $\bigl\{\bigl(\delta\cdot d,\,-\sign(d_{m.c.(d)}),0\bigr) : d\in \I^{n-2}(1)\bigr\}$ for the same $\delta>0$ as above (these are $2^{n-2}$ directions);
\item[-] if $n>4$, the directions $\bigl\{\delta\cdot\bigl(d,\,-\sign(d_{m.c.(d)})\bigr): d\in \J_{n-2}\bigr\}\times\{\pm 1\}$ for the same $\delta>0$ as above (these are $\leq (2^{n-3}-2)\cdot 2 = 2^{n-2}-4$ directions).
\end{itemize}
{\bf Important Clarification.}
Note that, in the cases of $n=3$ and $n=4$, we have not defined a set $\J_{n-2}$. In these cases, as we will see, we don't need a third subfamily of directions, and we can quickly write down what the sets $\I_{ex2,1}^3(\delta)$ and $\I_{ex2,1}^4(\delta)$ look like: $\I_{ex2,1}^3(\delta) :=\bigl\{\pm(\delta,\delta,+1),\pm(\delta,\delta,-1),\,\pm(\delta,-1,0)\bigr\}$ and
\begin{align*}
\I_{ex2,1}^4(\delta) := \bigl\{&\pm(\delta,\delta,\delta,1),\pm(\delta,\delta,\delta,-1),\pm(-\delta,\delta,\delta,1),\pm(-\delta,\delta,\delta,-1),
\\
&\pm(\delta,\delta,-1,0),\pm(-\delta,\delta,-1,0)\bigr\}.
\end{align*}
\item[2.] The set $\I_{ex2,2}^n(\epsilon,\delta_\epsilon)$, for some $\epsilon>0$ and $\delta_\epsilon>0$ which depend only on $\B$, which is defined as below:
\begin{align*}
\I_{ex2,2}^n(\epsilon,\delta_\epsilon) &:=\{\pm(\epsilon,\epsilon,\ldots,\epsilon, 1, 1)\} 
\\
&\qquad\ \  \bigcup \bigl[\widehat{\I}^n_{n-1,n}(\delta_\epsilon)\backslash\{\pm(\delta_\epsilon,\delta_\epsilon,\ldots,\delta_\epsilon,1,\delta_\epsilon),\,\pm(\delta_\epsilon,\delta_\epsilon,\ldots,\delta_\epsilon,-\delta_\epsilon,1)\}\bigr].
\end{align*}
\end{itemize}
We then see that in both cases $\II(\B) \leq 2^n-2$.
\end{proposition}
\begin{proof}
WLOG assume that $i_1=n-1$ and $i_2=n$. In other words, $\B$ contains the points $(1,1,1,\ldots,1,0,1)$ and $(1,1,1,\ldots,1,1,0)$, but does not contain other similar points, that is, points of the form ${\bm 1} - e_j$ for $j\in [n-2]$. We distinguish two cases.
\begin{itemize}
\item[Case 1.] At least one of the points ${\bm 1}-e_{n-1}$ and ${\bm 1}-e_n$ is NOT an extreme point of $\B$. In other words, there is $a>0$ such that $ae_{n-1}+\sum_{i\in [n]\backslash\{n-1\}}e_i \in\B$ or $ae_n+\sum_{i\in [n-1]}e_i\in \B$. WLOG assume that $\B$ certainly contains a point of the form $ae_n+\sum_{i\in [n-1]}e_i$ for some positive $a$. 

Let $\alpha_{n-1}$ be the supremum of all $y_{n-1}\geq 0$ such that \[y_{n-1}e_{n-1}+\sum_{i\in [n]\backslash\{n-1\}}e_i \in\B.\] By compactness we have that $\alpha_{n-1}e_{n-1}+\sum_{i\in [n]\backslash\{n-1\}}e_i \in\B$, and thus $0\leq \alpha_{n-1} <1$. Similarly, set $\alpha_n$ to be the supremum of all $y_n\geq 0$ such that \[y_ne_n+\sum_{i\in [n-1]}e_i \in\B.\] By compactness we have that $\alpha_ne_n+\sum_{i\in [n-1]}e_i \in\B$, and by our assumptions it follows that $1>\alpha_n > 0$.

\medskip

As in the previous propositions, for every $j\in [n-2]$, set $\theta_j=\norm{{\bm 1}-e_j}_\B^{-1}$, and $\Theta_0=\max_{j\in [n-2]}\theta_j$. Also, set $\gamma=\norm{{\bm 1}}_\B^{-1}$. Clearly, for all $j\in [n-2]$, $\,0 < \gamma\leq \theta_j\leq \Theta_0 < 1$. On the other hand, $\max(\alpha_{n-1},\alpha_n) \leq \gamma$ (because otherwise $\gamma{\bm 1}$ would not be a boundary point of $\B$). We will show that, in the setting here, $\I_{ex2,1}^n(\delta)$ illuminates $\B$ as long as $\delta < 1-\Theta_0$.

Consider a boundary point $x$ of $\B$.
\begin{itemize}
\item[$\blacklozenge$] Assume first that $x=\alpha_{n-1}e_{n-1}+\sum_{i\in [n]\backslash\{n-1\}}e_i $ or one of its coordinate reflections. Then we find the unique direction $d_x$ in $\I^{n-2}(1)$ which has opposite signs to $x$ in each of the first $n-2$ coordinates. Observe that $d^\prime_x=\bigl(\delta\cdot d_x,\delta\sign(d_{x,m.c.(d_x)}),-\sign(x_n)\bigr)$ illuminates $x$ (we can compare the displaced vector $x+d^\prime_x$ with the point ${\bm 1}-e_n$ to confirm this; this is because $\alpha_{n-1}+\delta\leq \gamma+\delta\leq \Theta_0 + \delta < 1$).
\item[$\blacklozenge$] Similarly, the first two types of directions in $\I_{ex2,1}^n(\delta)$ illuminate all the coordinate reflections of $\alpha_ne_n+\sum_{i\in [n-1]}e_i$ \textit{(here is where we use the assumption that ${\bm 1}-e_n$ is NOT an extreme point of $\B$, and thus that $\alpha_n > 0$; indeed, this allows us to have access to all of the $2^{n-1}$ combinations of signs for the first $n-1$ coordinates, that are all needed, without having to include $2^{n-1}$ all new directions, which would otherwise lead to an illuminating set of larger than desired size).} 

In more detail now, if $x$ is one of these coordinate reflections, again we pick the unique direction $d_x$ in $\I^{n-2}(1)$ which has opposite signs to $x$ in each of the first $n-2$ coordinates. In the case that $\sign(d_{x,m.c.(d_x)})=\sign(x_{n-1})$, then the direction $d^\prime=\bigl(\delta\cdot d_y,-\sign(d_{y,m.c.(d_y)}),0\bigr)$ (of the second type that we included) illuminates $x$ (this is because $|x_n| = \alpha_n \in (0,1)$, and thus we can use Corollary \ref{cor:affine-set}). 

Otherwise, the direction $d^{\prime\prime}=\bigl(\delta\cdot d_x,\delta\sign(d_{x,m.c.(d_x)}),-\sign(x_n)\bigr)$ (of the first type) will work instead.
\medskip
\item[$\blacklozenge$] Next, assume that either $|x_{n-1}|\leq \alpha_{n-1}$ or $|x_n|\leq \alpha_n$ (or both). Then $x$ is in the convex hull of the points $\alpha_{n-1}e_{n-1}+\sum_{i\in [n]\backslash\{n-1\}}e_i$ and $\alpha_ne_n+\sum_{i\in [n-1]}e_i$ and their coordinate reflections, and thus it is illuminated by some of the directions we have already used.
\item[$\blacklozenge$] We now suppose that $|x_{n-1}| > \alpha_{n-1}\geq 0$ AND $|x_n| > \alpha_n > 0$. 
\begin{itemize}
\item[$\bullet$] Assume first that $\min\{|x_{n-1}|,\,|x_n|\}\leq \Theta_0$.
\begin{itemize}
\item[--] In this case, if $|x_n| >\Theta_0$, then necessarily $|x_{n-1}|\leq \Theta_0$. Thus, as for the point $\alpha_{n-1}e_{n-1}+\sum_{i\in [n]\backslash\{n-1\}}e_i$ and its coordinate reflections, we can use a direction of the first type in $\I_{ex2,1}^n(\delta)$ to illuminate the boundary point $x$ that we are considering now (it suffices to pick $d_x$ such that $d_{x,i}\cdot x_i \leq 0$ for all $i\in [n-2]\cup\{n\}$, and then compare the displaced vector $x+|x_n|d_x$ with the point ${\bm 1}-e_n$).
\item[--] If $|x_n|\leq \Theta_0 < 1$, then, since we also have $|x_{n-1}| > 0$, we can pick a direction $d_x$ from the first two types in $\I_{ex2,1}^n(\delta)$ so that $d_{x,i}\cdot x_i\leq 0$ for all $i\in [n]$. Depending on whether $m.c.(d_x)=n-1$ or not, we consider the displaced vector $x+|x_{n-1}|d_x$ or the displaced vector $x+|x_n|d_x$, and we compare with the points ${\bm 1}-e_{n-1}$ or ${\bm 1}-e_n$ respectively. 
\end{itemize}
\item[$\bullet$] We finally consider the cases where $\min\{|x_{n-1}|,\,|x_n|\} > \Theta_0\geq \gamma$. Then we can find $s_1\in [n-2]$ such that $|x_{s_1}|\leq\gamma < 1-\delta$.
\begin{itemize}
\item[--] In the case where $n=3$, we can quickly confirm that $x$ is illuminated by one of the directions $\pm(\delta,\delta,\pm 1)$.
\item[--] If instead $n\geq 4$, then we can find one more index $s_2\in [n-2]\setminus\{s_1\}$ such that $|x_{s_2}|\leq \Theta_0$. Indeed, since $\theta_{s_1}({\bm 1}-e_{s_1})\in \partial \B$, we should be able to find an index $s_2\in [n]\setminus\{s_1\}$ such that $|x_{s_2}|\leq \theta_{s_1}\leq \Theta_0$. Since both $|x_{n-1}|$ and $|x_n|$ are $>\Theta_0$, we must have $s_2\in [n-2]\setminus\{s_1\}$. 

\medskip

We can now conclude the following: if $n=4$, then $\{s_1,s_2\}=\{1,2\}$, and we have that $x$ is illuminated by one of the directions $\pm(\delta,\delta,\delta,1),\,\pm(\delta,\delta,\delta,-1)$ (indeed, set $d_x$ to be the unique direction among these that satisfies $d_{x,i}\cdot x_i < 0$ for $i\in \{n-1,n\}$, and then compare $x+|x_n|d_x$ with the point ${\bm 1}-e_n$ to conclude that the former is in $\intr\B$).

\bigskip

If instead $n>4$, then we have access to the set $\J_{n-2}$ that we previously defined, and hence we can find a direction $d^\prime\in \J_{n-2}$ such that $d^\prime_j\cdot x_j \leq 0$ for all $j\in [n-2]\setminus\{s_1,s_2\}$. Clearly $d^\prime\in \I^{n-2}(1)$ as well, and thus both vectors 
\begin{equation*}
\bigl(\delta\cdot d^\prime,\,\delta\sign(d^\prime_{m.c.(d^\prime)}),\,-\sign(x_n)\bigr) \ \  \hbox{and}\ \  
\bigl(\delta\cdot d^\prime,\,-\delta\sign(d^\prime_{m.c.(d^\prime)}),\,-\sign(x_n)\bigr)
\end{equation*}
are in $\I_{ex2,1}^n(\delta)$. The one which has $(n-1)$-th coordinate of opposite sign to $x_{n-1}$ will illuminate $x$ (and we can compare the displaced vector to the point ${\bm 1}-e_n$ to confirm this).
\end{itemize}
\end{itemize}
\end{itemize}
We have thus illuminated all boundary points of $\B$ (regardless of what the dimension $n\geq 3$ is), when Case 1 holds.
\item[Case 2.] Both of the points ${\bm 1}-e_{n-1}$ and ${\bm 1}-e_n$ are extreme points of $\B$. Recall that we have assumed that $\B$ is not an affine image of the cube, however the convex hull of these two points and of their coordinate reflections is an affine image of the cube, the convex body $[-1,1]^{n-2}\times CP_1^2$. Thus there must exist a point $x\in \B$ outside of this convex hull, which in particular implies that $|x_{n-1}|+|x_n| > 1$. 
Just as in Lemma \ref{lem:R^3-exactly-two-pairs-interior}, we can show that this entails that
\begin{equation*}
\tfrac{1}{2}(e_{n-1}+e_n)\in \intr\B, \quad \hbox{or equivalently} \ 
\norm{e_{n-1}+e_n}_\B^{-1} > \frac{1}{2},
\end{equation*}
which further implies that
\begin{equation*}
\hbox{for every $\epsilon>0$, the point $(1-\epsilon,1-\epsilon,1-\epsilon,\ldots,1-\epsilon,\,\tfrac{1}{2},\,\tfrac{1}{2})\in\intr\B$.}
\end{equation*}
Fix some $\epsilon_0\in (0,1)$, and note that, because the point
\begin{equation*}
(1-\tfrac{\epsilon_0}{2},\,1-\tfrac{\epsilon_0}{2},\,1-\tfrac{\epsilon_0}{2},\ldots,1-\tfrac{\epsilon_0}{2},\ \tfrac{1}{2},\,\tfrac{1}{2})\in\intr\B,
\end{equation*}
we can find $\zeta_0\in (0,\tfrac{1}{2})$ such that
\begin{equation*}
(1-\tfrac{\epsilon_0}{2},1-\tfrac{\epsilon_0}{2},1-\tfrac{\epsilon_0}{2},\ldots,1-\tfrac{\epsilon_0}{2},\ \tfrac{1}{2}+\zeta_0,\,\tfrac{1}{2}+\zeta_0)
\end{equation*}
is also an interior point of $\B$. 

\medskip

Next, we set $\beta_{n-1}$ to be the supremum of positive constants $c$ such that
\begin{equation*}
\zeta_0 e_{n-1} + c\!\!\!\sum_{i\in [n]\setminus\{n-1\}}\!\!\!e_i\in \B.
\end{equation*}
By compactness $\zeta_0 e_{n-1} + \beta_{n-1}\sum_{i\in [n]\setminus\{n-1\}}e_i\in \B$, and by the main assumption of Case 2, $\beta_{n-1}\in (0,1)$. Similarly, we set $\beta_n$ to be the supremum of positive constants $c^\prime$ such that
\begin{equation*}
\zeta_0e_n + c^\prime\sum_{i\in [n-1]}e_i\in \B.
\end{equation*}
Again, we have $\zeta_0e_n + \beta_n\sum_{i\in [n-1]}e_i\in \B$ and $\beta_n\in (0,1)$.

\medskip

Now we pick 
\begin{equation*}
\delta_0 < \min(1-\zeta_0, \,1-\beta_{n-1},\,1-\beta_n)
\end{equation*}
and consider the corresponding set $\I_{ex2,2}^n(\epsilon_0,\delta_0)$. We will show that this illuminates $\B$.

\smallskip

Consider a boundary point $x$ of $\B$.
\begin{itemize}
\item[$\blacklozenge$] If $x$ is the point ${\bm 1}-e_{n-1}$ or one of its coordinate reflections, then we find a direction $d_x\in \I_{ex2,2}^n(\epsilon_0,\delta_0)$ which has opposite signs to $x$ in all the first $n-2$ entries, as well as the $n$-th entry, and has maximum entry the $n$-th entry. Unless $\sign(x_1)=\sign(x_2)=\sign(x_3)=\cdots=\sign(x_{n-2})=\sign(x_n)$, we take $d_x$ from the set 
\begin{equation*}
\widehat{\I}^n_{n-1,n}(\delta_0)\backslash\{\pm(\delta_0,\delta_0,\ldots,\delta_0,1,\delta_0),\pm(\delta_0,\delta_0,\ldots,\delta_0,-\delta_0,1)\}
\end{equation*}
and it is clear that $x+d_x\in \intr\B$ if we compare to the point ${\bm 1}-e_n$. In the remaining case, $d_x$ will be one of the directions $\pm(\epsilon_0,\epsilon_0,\ldots,\epsilon_0, 1, 1)$, and again we can see that
\begin{equation*}
x+\tfrac{1}{2}d_x\in \intr\B
\end{equation*}
since $x+\tfrac{1}{2}d_x$ will be a coordinate reflection of the point $(1-\tfrac{\epsilon_0}{2},1-\tfrac{\epsilon_0}{2},\ldots,1-\tfrac{\epsilon_0}{2},\,\tfrac{1}{2},\,\tfrac{1}{2})$.
\item[$\blacklozenge$] Analogously we illuminate the point ${\bm 1}-e_n$ and its coordinate reflections. 
\item[$\blacklozenge$] Given the above, we have now also illuminated all boundary points $x$ of $\B$ which fall in the convex hull of the points ${\bm 1}-e_{n-1}$, $\,{\bm 1}-e_n$ and their coordinate reflections; thus all boundary points $x$ which satisfy $|x_{n-1}|+|x_n| \leq 1$.
\medskip
\item[$\blacklozenge$] Suppose now that $|x_{n-1}|+|x_n| > 1$.
\begin{itemize}
\item[$\bullet$] Assume first that $|x_{n-1}|\leq \zeta_0$. Again, we find a direction $d_x$ which has opposite signs to $x$ in all the first $n-2$ entries, as well as the $n$-th entry, and has maximum entry the $n$-th entry. Unless $\sign(x_1)=\sign(x_2)=\sign(x_3)=\cdots=\sign(x_{n-2})=\sign(x_n)$, we take $d_x$ from the set 
\begin{equation*}
\widehat{\I}^n_{n-1,n}(\delta_0)\backslash\{\pm(\delta_0,\delta_0,\ldots,\delta_0,1,\delta_0),\pm(\delta_0,\delta_0,\ldots,\delta_0,-\delta_0,1)\}
\end{equation*}
and then compare $x+|x_n|d_x$ to the point ${\bm 1}-e_n$ (since in particular $|(x+|x_n|d_x)_{n-1}|\leq |x_{n-1}|+|x_n|\delta_0 \leq \zeta_0 + \delta_0 < 1$). 

In the remaining case, we choose $d_x$ from $\pm(\epsilon_0,\epsilon_0,\ldots,\epsilon_0, 1, 1)$. We can see that
\begin{equation*}
x+\tfrac{1}{2}d_x\in \intr\B
\end{equation*}
because we will have 
\begin{itemize}
\item $\big|(x+\frac{1}{2}d_x)_i\big| \leq \max(1-\tfrac{\epsilon_0}{2},\tfrac{\epsilon_0}{2}) = 1-\tfrac{\epsilon_0}{2}$ for all $i\in [n-2]$, 
\item $\big|(x+\frac{1}{2}d_x)_n\big|\leq \max(1-\tfrac{1}{2},\tfrac{1}{2}) = \tfrac{1}{2}$, 
\item and finally $\big|(x+\frac{1}{2}d_x)_{n-1}\big|\leq |x_{n-1}|+\tfrac{1}{2} \leq \tfrac{1}{2}+\zeta_0$, 
\end{itemize}
while the point
\begin{equation*}
(1-\tfrac{\epsilon_0}{2},1-\tfrac{\epsilon_0}{2},1-\tfrac{\epsilon_0}{2},\ldots,1-\tfrac{\epsilon_0}{2},\ \tfrac{1}{2}+\zeta_0,\,\tfrac{1}{2}+\zeta_0)
\end{equation*}
is an interior point of $\B$, given our choice of $\zeta_0$.
\smallskip
\item[$\bullet$] We argue analogously if $|x_n|\leq \zeta_0$.
\medskip
\item[$\bullet$] Finally, let us assume that $\min(|x_{n-1}|,|x_n|)>\zeta_0$. 
\begin{itemize}
\item[--] Suppose also that $\max(|x_{n-1}|,|x_n|)\leq \max(\beta_{n-1},\beta_n) < 1-\delta_0$.

If $x_{n-1}\cdot x_n < 0$, then we can pick a direction $d_x$ from the set
\begin{equation*}
\widehat{\I}^n_{n-1,n}(\delta_0)\backslash\{\pm(\delta_0,\delta_0,\ldots,\delta_0,1,\delta_0),\pm(\delta_0,\delta_0,\ldots,\delta_0,-\delta_0,1)\}
\end{equation*}
to illuminate $x$. In fact, in most cases we can take $d_x$ from the second half of this set, which contains directions with maximum (in absolute value) coordinate the $n$-th one (and consider the vector $x+|x_n|d_x$, comparing it to ${\bm 1}-e_n$). This will work in all cases except when \underline{all entries of $x$ are non-zero} and
\begin{equation*}
\sign(x_1)=\sign(x_2)=\sign(x_3)=\cdots=\sign(x_{n-2})=\sign(x_n) = -\sign(x_{n-1}).
\end{equation*}
In this last subcase, we instead choose $d_x$ from $\pm(\delta_0,\delta_0,\delta_0,\ldots,\delta_0,-1,-\delta_0)$ (and consider the vector $x+|x_{n-1}|d_x$, comparing it to ${\bm 1}-e_{n-1}$, given that we have $\big|(x+|x_{n-1}|d_x)_n\big| = |x_n|+|x_{n-1}|\delta_0\leq |x_n|+\delta_0 < 1$ by our last assumption above).

\smallskip

Similarly, if $x_{n-1}\cdot x > 0$, we choose $d_x$ from the first half of the set $\widehat{\I}^n_{n-1,n}(\delta_0)$, except in the case where \underline{all entries of $x$ are non-zero} and
\begin{equation*}
\sign(x_1)=\sign(x_2)=\sign(x_3)=\cdots=\sign(x_{n-2})=\sign(x_{n-1}) = \sign(x_n).
\end{equation*}
In this last subcase we can instead choose $d_x$ to be one of the directions $\pm(\epsilon_0,\epsilon_0,\ldots,\epsilon_0, 1, 1)$ and consider the vector $x+\varepsilon_xd_x$ where $\varepsilon_x = \min_{i\in [n]}|x_i|$.
\item[--] Now, suppose that $\max(|x_{n-1}|,|x_n|) > \max(\beta_{n-1},\beta_n)$. Given also our `parent' assumption that $\min(|x_{n-1}|,|x_n|) > \zeta_0$, by the choice of the constants $\beta_{n-1}$ and $\beta_n$ we can find $j_0\in [n-2]$ such that $|x_{j_0}| \leq \max(\beta_{n-1},\beta_n) < 1-\delta_0$. But then it is possible to find a direction $d_x$ from 
\begin{equation*}
\widehat{\I}^n_{n-1,n}(\delta_0)\backslash\{\pm(\delta_0,\delta_0,\ldots,\delta_0,1,\delta_0),\pm(\delta_0,\delta_0,\ldots,\delta_0,-\delta_0,1)\}
\end{equation*}
such that $d_{x,i}\cdot x_i \leq 0$ for all $i\in [n]\backslash\{j_0\}$, and we can check that this $d_x$ illuminates $x$.
\end{itemize}
\end{itemize}
\end{itemize}
\end{itemize}
We have completed the proof in both Case 1 and Case 2.
\end{proof}

\smallskip

\begin{remark}\label{rem:inductive-hyp}
(I) Combining Propositions \ref{prop:all-n-1-tuples}, \ref{prop:all-but-one-n-1-tuples} and \ref{prop:exactly-n-k-n-1-tuples}, Proposition \ref{prop:exactly-one-n-1-tuple} (and its 3-dimensional version, Proposition \ref{prop:R^3-exactly-one-pair}) and Proposition \ref{prop:exactly-two-n-1-tuples}, we reach the following conclusion: for any dimension $n\geq 3$, if $\B\in \U^n$ is NOT an affine image of the cube, but has the property that there is \underline{at least one} index $i_0\in [n]$ such that ${\bm 1}-e_{i_0}\in\B$, then $\II(\B)\leq 2^n-2$. In other words, if $\B$ contains at least one maximal unit subcube (according to our terminology) but is not a parallelepiped, then $\II(\B)\leq 2^n-2$.

\medskip

(II) Moreover, Proposition \ref{prop:all-n-1-tuples} and the proof of Proposition \ref{prop:exactly-two-n-1-tuples} allow us to complete the discussion of Section \ref{sec:3-dim-bodies}: we can conclude that, for ALL 1-unconditional bodies $\B$ in $\R^3$, $\II(B)\leq 6$, unless $\B$ is a parallelepiped.

\medskip

(III) As a `bonus', we have also confirmed all the above results while using illuminating sets which consist of pairs of opposite directions.
\end{remark}

\subsection{Bodies with all unit subcubes of dimension $n-2$}

We finish this section by proving a similar result to the above, namely Proposition \ref{prop:main-all-unit-subcubes-dim-n-2} of the Introduction, since a similar argument can work here as well. 
To keep the proof simpler, we assume that the convex bodies that we will consider are not already covered by any of the previous propositions, or in other words by Theorem \ref{thm:main-unit-subcubes}.

\begin{proposition}\label{prop:all-n-2-tuples}
\textup{Let $n\geq 4$ and let $\B\in {\cal U}^n$ with the property that $\B$ contains
\begin{equation*}
e_1 + e_2 + \cdots + e_{n-2}
\end{equation*}
and all its coordinate permutations, but it does not contain any coordinate permutation of 
\begin{equation*}
e_1 + e_2 + \cdots + e_{n-2} + e_{n-1}.
\end{equation*}
For $\delta\in (0,1)$ and $\zeta=\zeta_\delta\in (0,\delta)$, let $\widetilde{\I}_{\delta,\zeta}$ be the following set of directions: 
\begin{itemize}
\item[---] we will include most of the directions in $\widehat{\I}^n_{n-2,n-1,n}(\delta)$ {\bf except}
\begin{itemize}
\item[(i$^\prime$)] the directions $d_+=(\delta,\delta,\ldots,\delta,-\delta,1)$ and $d_-=-(\delta,\delta,\ldots,\delta,-\delta,1)$, 
\item[(ii$^\prime$)] and all the other directions $d\in \widehat{\I}^n_{n-2,n-1,n}(\delta)$ whose sequence of coordinate signs differs from that of either $d_+$ or $d_-$ in exactly one place (e.g. $\pm(\delta,\delta,\ldots,\delta,1,\delta,\delta)$ or $\pm(\delta,\delta,\ldots,\delta,-\delta,-\delta,1)$); note that, for each of these directions $d$, the place where its sequence of signs differs from that of $d_+$ or $d_-$ is not the place where $d$ has its maximum in absolute value coordinate.
\end{itemize}
\item[---] For each of the directions $d$ in (ii$^\prime$) we introduce a `replacement' direction $d^\prime$ as follows: we first distinguish whether $d$ has an almost identical sequence of signs to that of $d_+$ or to that of $d_-$ (suppose for illustration purposes that it is $d_+$ here). We then set $d^\prime$ to be the direction which has the same respective entries as $d$, except for the one entry $d_{i_d}$ of $d$ which differs in sign from the respective entry of $d_+$, in which case we set $d^\prime_{i_d}=\sign(d_{i_d})\zeta$ (e.g. if $d= (\delta,\delta,\ldots,\delta,-\delta,-\delta,1)$, then $d^\prime = (\delta,\delta,\ldots,\delta,-\zeta,-\delta,1)$).
\end{itemize}
Then there are $\delta=\delta_\B>0$ and $\zeta=\zeta_\B>0$ 
such that $\B$ can be illuminated by the corresponding set $\widetilde{\I}_{\delta,\zeta}$.}
\end{proposition}
\begin{remark}
For illustration purposes, let's write down what the sets $\widetilde{\I}_{\delta,\zeta}$ would look like in $\R^4$: we get the set of directions
\begin{align*}
\bigl\{\!&\pm(\delta, 1,\zeta,\delta),\ \pm(\delta,-1,-\delta,-\delta),\ \,\pm(\delta,\delta,-1,-\zeta),\ \pm(\delta,-\delta, 1,\delta),
\\
&\pm(\delta,\delta,\delta,-1),\ \pm(\delta,-\zeta,-\delta,1),\ \ \pm(-\zeta,\delta,-\delta,1)\bigr\}.
\end{align*}
\end{remark}
\begin{proof}
We first verify the following 
\\
{\bf Claim.} Let $y\in {\mathbb R}^n\setminus\{\vec{0}\}$ such that $\abs{\cZ_y}=2$. Then we can find $d\in\widetilde{\I}_{\delta,\zeta}$ which deep illuminates $y$, and which in addition satisfies the following: if $d_i = \zeta$ for some $i\in [n]$ (unique in our setting), then $y_i=0$ (in other words, if $d$ is one of the modified directions in $\widetilde{\I}_{\delta,\zeta}$, and $d_i$ is the modified coordinate, then this corresponds to one of the two zero coordinates of $y$).
\smallskip\\
\emph{Proof of Claim.} Note that, due to the way we construct $\widetilde{\I}_{\delta,\zeta}$ from $\widehat{\I}^n_{n-2,n-1,n}(\delta)$ (and ultimately from $\I^n(\delta)$), all sequences of signs whose last 3 terms take one of the forms $\pm (1,1,1)$ or $\pm(-1,1,1)$ or $\pm(1,1,-1)$ are still there (there are $2^{n-2}$ sequences of signs of each such form). Moreover, there is exactly one pair of opposite directions in each of these subgroups which comes from the modified directions: in fact, 
\begin{itemize}
\item[$\cdot$] in the first subgroup the coordinate which may be equal to $\zeta$ is the $(n-1)$-th one, 
\item[$\cdot$] in the second subgroup the coordinate which may be equal to $\zeta$ is the $n$-th one, 
\item[$\cdot$] and in the third subgroup the coordinate which may be equal to $\zeta$ is the $(n-2)$-th one. 
\end{itemize} 
At the same time, we observe that \textit{(because of the specific, combinatorial construction of $\I^n(\delta)$ that we rely on in this paper, and then the construction of $\widehat{\I}^n_{n-2,n-1,n}(\delta)$ from that)} 
\begin{itemize}
\item[$\cdot$] the maximum (in absolute value) coordinate of all directions in the first subgroup is the $(n-2)$-th one, 
\item[$\cdot$] the maximum coordinate of all directions in the second subgroup is the $(n-1)$-th one, 
\item[$\cdot$] and the maximum coordinate of all directions in the third subgroup is the $n$-th one.
\end{itemize}

We can now analyse what $d$ should be, based on where the zero coordinates of $y$ are found.
\begin{itemize}
\item[Case 1:] The two zero coordinates of $y$ are among the last three ones. Then the remaining coordinate from these, say coordinate $i_0\in \{n-2,n-1,n\}$, is non-zero, and so are all the coordinates with index $<n-2$. Thus we can focus on one of the first three subgroups of directions in $\widetilde{\I}_{\delta,\zeta}$, which contains directions with maximum (in absolute value) coordinate the $i_0$-th one, and pick the unique direction $d$ whose signs on the non-zero coordinates of $y$ are opposite to the corresponding signs of $y$ (in such a case, even if $d$ has a coordinate equal to $\pm \zeta$, this will have index in $\{n-2,n-1,n\}\setminus \{i_0\}$). 
\item[Case 2:] Only one of the zero coordinates of $y$ is among the last three ones, say the coordinate with index $i_1\in \{n-2,n-1,n\}$. Let us also write $i_2$ for the index of the other zero coordinate of $y$: $i_2 < n-2$. For illustration purposes, let's assume that $i_1=n$ (the other cases can be treated completely analogously). Then we can use directions either from the first subgroup (if the $(n-2)$-th and $(n-1)$-th coordinates of $y$ have the same sign), or from the second subgroup (if these coordinates of $y$ have opposite signs). To avoid the one pair of directions in these subgroups which has an $(n-1)$-th coordinate equal to $\zeta$ (in absolute value), we can focus on directions $d$ which satisfy $\sign(d_{i_2})\neq \sign(d_{n-1})$ (since $d_{i_2}$ will correspond to a zero coordinate of $y$, so it can have either positive or negative sign without issue).
\item[Case 3:] Both of the zero coordinates of $y$ have indices $< n-2$, say indices $i_3$ and $i_4$ (where $1\leq i_3 < i_4\leq n-3$). Clearly this case can occur only when $n\geq 5$. We first focus on the subgroup of directions in $\widetilde{\I}_{\delta,\zeta}$ whose sequences of signs in the last three coordinates match or are exactly opposite to the respective sequence for $y$. From within this subgroup, it suffices to consider those directions $d$ which satisfy $\sign(d_{i_3})\neq \sign(d_{i_4})$ (because in this way we both avoid the one pair of opposite directions/sign-sequences missing from $\widetilde{\I}_{\delta,\zeta}$ compared to $\widehat{\I}^n_{n-2,n-1,n}(\delta)$, which we wouldn't have been able to pick anyway, and also we make sure that, even if a suitable direction $d$ has a coordinate equal to $\pm \zeta$, this will be its $i_3$-th or its $i_4$-th one, as desired).
\end{itemize}
The proof of the claim is complete.

\medskip

Since $\B$ does not contain any of the coordinate permutations of ${\bm 1}-e_n$, as previously we can set, for each $j\in [n]$, $\,\theta_j:=\norm{{\bm 1}-e_j}_\B^{-1}$; we will have $\theta_j\in (0,1)$. We also set $\Theta_0:=\max_{j\in[n]}\theta_j$, and pick 
\begin{equation*}
0 < \delta < \min\left\{\frac{1}{6},\ \frac{1-\Theta_0}{2}\right\}.
\end{equation*}
Furthermore, set $\gamma:=\norm{{\bm 1}}_\B^{-1}$. 
A suitable value for $\zeta$ will become clear towards the end of the proof, but for now we just make sure that $\zeta < \delta$.

\smallskip

Let $x$ be an \underline{extreme} point of $\B$. We distinguish two main cases.
\begin{itemize}
\item[$\abs{\cZ_x}\geq 1$.] Suppose $x_{i_0}=0$ for some $i_0\in [n]$. Because $\theta_{i_0}\bigl({\bm 1}-e_{i_0}\bigr)\in \partial\B$, for at least one index $i_1\in [n]\setminus\{i_0\}$ we must have $\abs{x_{i_1}}\leq \theta_{i_0}\leq \Theta_0 < 1$. Moreover, for every $j\in [n]\setminus\{i_0,i_1\}$, we must have $x_j\neq 0$, because otherwise $x$ would not be an extreme point of $\B$ (it would be in the convex hull of a point of the form ${\bm 1}-e_{i_0}-e_j$, $\,j\in [n]\setminus\{i_0,i_1\}$, and of its coordinate reflections, without being any of those points). Thus the point
\begin{equation*}
y_x:= x \,-\,x_{i_0}e_{i_0} - x_{i_1}e_{i_1}  = x\,-\,x_{i_1}e_{i_1}
\end{equation*}
has exactly two zero coordinates, and hence, by the above claim, we can find a direction $d\in\widetilde{\I}_{\delta,\zeta}$ which deep illuminates $y_x$. Moreover, we can make sure that, if $d_s = \pm\zeta$ for some $s\in [n]$, then $y_{x,s}=0$ (or, in other words, $|d_s|=\zeta \Rightarrow s\in \{i_0,i_1\}$). 

Let $j_0 = m.c.(d)$. Then $x+\abs{x_{j_0}}d$
satisfies the following: 
\begin{itemize}
\item[$\cdot$] its $j_0$-th coordinate is zero, 
\item[$\cdot$] while $\abs{(x+\abs{x_{j_0}}d)_{i_0}}\leq \abs{x_{j_0}}\delta$
\item[$\cdot$] and $\abs{(x+\abs{x_{j_0}}d)_{i_1}}\leq \Theta_0 + \abs{x_{j_0}}\delta$. Moreover, by our assumptions on $\delta$, we have $\Theta_0+\abs{x_{j_0}}\delta < 1-\abs{x_{j_0}}\delta$. 
\item[$\cdot$] Finally, for any $j\in [n]\setminus\{j_0,i_0,i_1\}$, we will have that
\begin{equation*}
\abs{(x+\abs{x_{j_0}}d)_j} \leq \max\bigl\{\abs{x_j}-\abs{x_{j_0}}\delta,\ \abs{x_{j_0}}\delta\bigr\}\leq 1- \abs{x_{j_0}}\delta.
\end{equation*}
\end{itemize}
We now compare this point to the convex combinations
\begin{multline*}
c_1(x):=(1-\abs{x_{j_0}}\delta)\bigl({\bm 1}-e_{i_0}-e_{j_0}\bigr) \ +\ \abs{x_{j_0}}\delta\cdot e_{i_0}
\\
\hbox{and}\ \ c_2(x):= (1-\abs{x_{j_0}}\delta)\bigl({\bm 1}-e_{i_0}-e_{j_0}-e_{i_1}\bigr) \ +\ \abs{x_{j_0}}\delta\cdot e_{i_0}\ +\ \abs{(x+\abs{x_{j_0}}d)_{i_1}}\cdot e_{i_1}
\end{multline*}
which are contained in $\B$. All the respective coordinates of $c_1(x)$ and $c_2(x)$ are equal, except for the $i_1$-th coordinate: in that case, $c_2(x)$ has a strictly smaller coordinate than $c_1(x)$. Moreover, the coordinates of $x+\abs{x_{j_0}}d$ do not exceed the corresponding ones of $c_2(x)$ in absolute value. Thus it suffices to show that $c_2(x)\in \intr\B$ to also obtain that $x+\abs{x_{j_0}}d\in \intr\B$. To do this, we will use the fact that
\begin{equation*}
c_3(x):=(1-\abs{x_{j_0}}\delta)\bigl({\bm 1}-e_{i_0}-e_{j_0}-e_{i_1}\bigr) \ +\ \abs{x_{j_0}}\delta\cdot e_{i_0}
\end{equation*}
is an interior point of $\B$ (which is verified if we compare with the point ${\bm 1}-e_{j_0}-e_{i_1}\in \B$). Thus the desired conclusions follow by applying Lemma \ref{lem:affine-set} on the section 
\[\bigl\{z\in \B: z_s= (1-\abs{x_{j_0}}\delta)\ \hbox{for all}\ s\in [n]\setminus\{i_0,i_1,j_0\},\ z_{j_0}=0,\ z_{i_0}=\abs{x_{j_0}}\delta\bigr\},\]
which contains all three points $c_1(x), c_2(x), c_3(x)$.
\item[$\abs{\cZ_x}=0$.] Here, it suffices to pick a direction $d$ from $\widetilde{\I}_{\delta,\zeta}$ which deep illuminates $x$ (which, in this case, simply means that $\sign(d_i)=-\sign(x_i)$ for all $i\in [n]$). This will not be possible only in the case that $\sign(x_i)=\sign(x_n)$ for all $i\in [n-2]$ and $\sign(x_{n-1})=-\sign(x_n)$ (since we removed without any replacement the only two directions in $\widehat{\I}^n_{n-2,n-1,n}(\delta)$ which had exactly this property). To deal with this remaining case, we distinguish two subcases.
\begin{itemize}
\item At least two coordinates of $x$ are $\leq \frac{1}{3}$ in absolute value. Say $\max(\abs{x_{i_1}}, \abs{x_{i_2}})\leq \frac{1}{3}$ with $1\leq i_1 < i_2 \leq n$. 
Recall that we have also assumed that $x$ has no zero coordinates.
For the point 
\begin{equation*}
y_x=x-x_{i_1}e_{i_1}-x_{i_2}e_{i_2}
\end{equation*}
we again use our initial claim, and find a direction $d\in\widetilde{\I}_{\delta,\zeta}$ which deep illuminates $y_x$, and is such that, if $d_s = \pm\zeta$ for some $s\in [n]$, then $y_{x,s}=0$ (in other words, $|d_s|=\zeta \Rightarrow s\in \{i_1,i_2\}$). 

Let $j_0=m.c.(d)$. Then $x+\abs{x_{j_0}}d$
satisfies the following: 
\begin{itemize}
\item[$\cdot$] its $j_0$-th coordinate is zero, 
\smallskip
\item[$\cdot$] while
\begin{equation*}
\abs{(x+\abs{x_{j_0}}d)_{i_s}} \leq \abs{x_{i_s}} + \abs{x_{j_0}}\delta \leq \frac{1}{3} + \delta < \frac{1}{2}
\end{equation*}
for both $s=1$ or $2$ (the last inequality holds because of our assumptions on $\delta$). 
\item[$\cdot$] Finally, for any $j\in [n]\setminus\{j_0,i_1,i_2\}$, we will have that
\begin{equation*}
\abs{(x+\abs{x_{j_0}}d)_j} \leq \max\bigl\{\abs{x_j}-\abs{x_{j_0}}\delta,\ \abs{x_{j_0}}\delta\bigr\}\leq 1- \abs{x_{j_0}}\delta.
\end{equation*}
\end{itemize}
We can now compare to the convex combination
\begin{equation*}
\frac{1}{2}\bigl({\bm 1}-e_{j_0}-e_{i_1}\bigr) \ +\  \frac{1}{2}\bigl({\bm 1}-e_{j_0}-e_{i_2}\bigr)
\end{equation*}
in $\B$, to conclude that $x+\abs{x_{j_0}}d\in \intr\B$.
\item At most one of the coordinates of $x$ is $\leq \frac{1}{3}$ in absolute value. Let $i_0$ be the (smallest) index of $\min_{i\in [n]}\abs{x_i}$. 
We must have $\abs{x_{i_0}}\leq \gamma$ (otherwise the point $\gamma {\bm 1}$ would not be a boundary point of $\B$). 
Moreover, recall that the only subcase that we have to still consider here is the one satisfying the following: 
\begin{itemize}
\item[(i)] for all $j\in [n]\setminus \{i_0\}$ we have $\abs{x_j} > \frac{1}{3} > \delta$; 
\item[(ii)] $\sign(x_i)=\sign(x_n)$ for all $i\in [n-2]$ and $\sign(x_{n-1})=-\sign(x_n)$. 
\end{itemize}

\smallskip

We now pick the {\bf unique} direction $d\in\widetilde{\I}_{\delta,\zeta}$ which satisfies $\abs{d_{i_0}}=\zeta$ and $\sign(d_{i_0})=\sign(x_{i_0})$, while $\sign(d_j)=-\sign(x_j)$ for all $j\in [n]\setminus\{i_0\}$. Let $j_0=m.c.(d)$. Then $x+\abs{x_{j_0}}d$
satisfies the following: for all $j\in [n]\setminus\{i_0\}$,
\begin{equation*}
\abs{(x+\abs{x_{j_0}}d)_j} \leq \abs{x_j}-\abs{x_{j_0}}\delta,
\end{equation*}
while
\begin{equation*}
\abs{(x+\abs{x_{j_0}}d)_{i_0}} = \abs{x_{i_0}}+\abs{x_{j_0}}\zeta \leq \abs{x_{i_0}} + \zeta.
\end{equation*}
Fix some $\lambda_0\in (0,\frac{\delta}{3})$; then, for all $j\in [n]\setminus\{i_0\}$, we can write
\[\lambda_0\abs{x_j} \leq \lambda_0 < \frac{\delta}{3} <  \abs{x_{j_0}}\delta,\]
and thus $\abs{x_j}-\abs{x_{j_0}}\delta < (1-\lambda_0)\abs{x_j}$. Based on this, we consider the convex combination
\begin{equation*}
c_0(x):=(1-\lambda_0)\sum_{i\in [n]}\abs{x_i}e_i \ +\ \lambda_0\cdot e_{i_0}  = (1-\lambda_0)\big|\vec{x}\big| \ +\ \lambda_0\cdot e_{i_0}
\end{equation*}
which is contained in $\B$, and we observe that:
\begin{itemize}
\item[$\cdot$] each of its coordinates with index $j\in [n]\setminus\{i_0\}$ is equal to $(1-\lambda_0)\abs{x_j} > \abs{(x+\abs{x_{j_0}}d)_j}$, 
\smallskip
\item[$\cdot$] while the coordinate with index $i_0$ is equal to
\begin{equation*}
(1-\lambda_0)\abs{x_{i_0}} + \lambda_0 = \abs{x_{i_0}} + \lambda_0(1-\abs{x_{i_0}}) \geq \abs{x_{i_0}} + \lambda_0(1-\gamma).
\end{equation*}
\end{itemize}
Hence, if we fix some $\zeta <  \lambda_0(1-\gamma)$ \textit{(e.g. fix $\lambda_0 = \frac{\delta}{4}$ and $\zeta = \frac{\delta}{5}(1-\gamma)$; note that none of these quantities has to depend on the point $x$ that we are considering)}, we can also ensure that $\abs{(x+\abs{x_{j_0}}d)_{i_0}}$ is strictly smaller than the $i_0$-th coordinate of the convex combination $c_0(x)$. This will imply that $x+\abs{x_{j_0}}d\in \intr\B$. 
\end{itemize}
\end{itemize}
We have addressed all possibilities for the arbitrary extreme point of $\B$, and we have shown how to illuminate it in each case using some direction from $\widetilde{\I}_{\delta,\zeta}$ (with $\delta$ and $\zeta$ suitably chosen with respect to $\B$). Therefore $\widetilde{\I}_{\delta,\zeta}$ is an illuminating set for $\B$, and $\II(\B) \leq 2^n - 2$.
\end{proof}

\section{Cubelike 1-unconditional convex bodies}\label{sec:cubelike}

Here we prove Theorem \ref{thm:main-cubelike-bodies}, namely we settle the Illumination Conjecture for 1-unconditional convex bodies $K$ which have the following property:
\begin{equation*}
\hbox{if $x$ is an extreme point of $K$, then $x_i\neq 0$ for all $1\leq i\leq n$} \tag{$\dagger$}
\end{equation*}
(we called these \emph{cubelike bodies}).

Recall that, because of standard results such as Fact A and Corollary \ref{cor:uncond-illum}, it is well-known that such convex bodies can be illuminated by $2^{{\rm dim}(K)}$ directions. Therefore, the novelty in Theorem \ref{thm:main-cubelike-bodies} is that we also verify the conjectured equality cases of the Illumination Conjecture.

\smallskip

The proof of Theorem \ref{thm:main-cubelike-bodies} could be summarised as follows: we will use induction in the dimension, and in the inductive step we will rely on combining two key lemmas, which we present first. 

The first of these lemmas could be of independent interest as well, since it applies in a broader setting than that of `cubelike' 1-unconditional convex bodies.

\begin{lemma}\label{lem:proj-sec}
Let $n\geq 3$, and let $K$ be a convex body in $\R^n$ and $H$ a(n) (affine) hyperplane of $\R^n$. Suppose that: (i) ${\rm Proj}_H(K) = K\cap H$ (where projection of any given vector here means translating the vector parallel to a normal vector to $H$ until we hit $H$), and (ii) $K$ has NO extreme points in $H$, that is, ${\rm ext}(K) \cap H=\emptyset$.

\smallskip

Then $\II(K) \leq 2\cdot \II(K\cap H) = 2\cdot\II({\rm Proj}_H(K))$ (note that $\II(K\cap H)$ is the illumination number of an $(n-1)$-dimensional convex body, found by illuminating $K\cap H= {\rm Proj}_H(K)$ within $H$).
\end{lemma}
\begin{proof}
WLOG we can assume that $H= e_n^\perp + a e_n$ for some $a\in \R$, and then, by translating both $H$ and $K$ by $-ae_n$, we can assume (for simplicity) that $H= e_n^\perp$. From now on, we will write $K_{e_n}$ instead of $K\cap H= K\cap e_n^\perp$.

Note also that, because $K\cap H= {\rm Proj}_H(K)$, we have that ${\rm aff}(K\cap H) = H$.

Next we observe that $K$ contains points $x$ with $x_n = \langle x, e_n\rangle > 0$, as well as points $y$ with $y_n < 0$. This is because, if this were not true, we would have that $K_{e_n}$ is a support set of $K$, and thus it would have to contain some extreme points of $K$, contrary to our second main assumption. As a consequence of this, we also get that $(\intr K) \cap e_n^\perp \neq \emptyset$.

\smallskip

Set now $N_0=\II(K_{e_n})$ (where we initially view $K_{e_n}$ as a subset of $\R^{n-1}$ instead of $e_n^\perp$). We can find a set ${\cal D} = \{d_1,d_2,\ldots, d_{N_0}\}$ of directions in $\R^{n-1}$ which illuminates $K_{e_n}$. Let us also restate this as a statement about subsets and directions of $\R^n$: ${\cal D}$ can be viewed as a (minimum-size) subset of $e_n^\perp$ which has the property that, for every $p\in {\rm relbd}(K_{e_n})$, we can find $d_i\in {\cal D}$ and $\varepsilon> 0$ such that $p+\varepsilon d_i \in {\rm relint}(K_{e_n})$.

\smallskip

{\bf Claim.} We can find a common $\varepsilon_0 > 0$ with the property that, for every $p\in K_{e_n}$ (and not just $\in {\rm relbd}(K_{e_n})$), there will be some $d_i\in {\cal D}$ such that $p+\varepsilon_0d_i\in {\rm relint}(K_{e_n})$. By convexity, this will also imply that, for $p$ and $d_i$ as before, $p+\varepsilon^\prime d_i\in {\rm relint}(K_{e_n})$ for every $\varepsilon^\prime \in (0, \varepsilon_0)$.

\smallskip

\emph{Proof of the claim.} We will use compactness (working with relatively open sets in the subspace topology on $e_n^\perp$, so as not to further complicate our notation). The following is essentially the (core of the) standard argument that shows that the illumination number and the covering number of a convex body coincide.

For every $p\in K_{e_n}$ we can find $d_{i_p}\in {\cal D}$ and some $\varepsilon_p>0$ such that $p+\varepsilon_pd_{i_p}\in {\rm relint}(K_{e_n})$. For points on the relative boundary of $K_{e_n}$, this is already guaranteed by our choice of the set ${\cal D}$. On the other hand, if $p\in {\rm relint}(K_{e_n})$, then, no matter which direction $d$ we choose from ${\rm aff}(K_{e_n}) = e_n^\perp$ \emph{(more accurately, from ${\rm aff}(K_{e_n}) -p$, which just happens to coincide with ${\rm aff}(K_{e_n})$ here)}, we can get the desired conclusion as long as we pick $\varepsilon_p=\varepsilon_{p,d}$ small enough.

We can rewrite this as $p\in -\varepsilon_pd_{i_p} + {\rm relint}(K_{e_n})$, and thus
\begin{equation*}
K_{e_n} \subset \bigcup_{p\in K_{e_n}} \bigl(-\varepsilon_pd_{i_p} + {\rm relint}(K_{e_n})\bigr).
\end{equation*}
Therefore, by compactness, we can find finitely many positive numbers $\varepsilon_1,\varepsilon_2,\ldots, \varepsilon_M$, $M\geq 1$, such that, for every $p\in K_{e_n}$, it will be possible to write
\begin{equation*}
p \in -\varepsilon_jd_i + {\rm relint}(K_{e_n}) \quad \Leftrightarrow \quad p+\varepsilon_jd_i\in {\rm relint}(K_{e_n})
\end{equation*}
for some $j\in \{1,2,\ldots, M\}$ and $i\in\{1,2,\ldots,N_0\}$. Finally, if we set $\varepsilon_0=\min\{\varepsilon_j: 1\leq j\leq M\}$, by convexity we will have that $p+\varepsilon_0d_i\in {\rm relint}(K_{e_n})$ as well, while $\varepsilon_0$ will not depend on the point $p$ anymore. The proof of the claim is complete.

\bigskip

We can finally define an illuminating set for the convex body $K$. Set
\begin{equation*}
a_0 = \max\{|x_n|:x\in K\}.
\end{equation*}
By our remarks at the beginning, we know that $a_0 > 0$. Set $\eta_0 = \frac{\varepsilon_0}{a_0}$. We claim that the set
\begin{equation*}
\eta_0{\cal D}\times \{\pm 1\} = \bigl\{(\eta_0 d_i, 1), (\eta_0d_i, -1): 1\leq i\leq N_0\bigr\}
\end{equation*}
illuminates $K$ (where we abuse our notation a bit again, and view ${\cal D}$ as a subset of $(n-1)$-dimensional vectors now).

Indeed, let $x$ be an \underline{extreme} point of $K$. Then, by our assumptions $x_n \neq 0$, and also
\begin{equation*}
{\rm Proj}_{e_n^\perp}(x) \in {\rm Proj}_{e_n^\perp}(K) = K\cap e_n^\perp = K_{e_n}.
\end{equation*}
Hence, we can find $d_i\in {\cal D}$ such that ${\rm Proj}_{e_n^\perp}(x) + \varepsilon^\prime d_i\in {\rm relint}(K_{e_n})$ for any $\varepsilon^\prime \in (0,\varepsilon_0]$.

\medskip

We will show that $x$ is illuminated by the direction $\bigl(\eta_0d_i, -\sign(x_n)\bigr)$:
\begin{equation*}
x+ |x_n|\bigl(\eta_0 d_i, -\sign(x_n)\bigr) = \bigl({\rm Proj}_{e_n^\perp}(x) + |x_n|\tfrac{\varepsilon_0}{a_0}\,d_i,\ 0\bigr)\in {\rm relint}(K_{e_n})
\end{equation*}
given that $|x_n|\frac{\varepsilon_0}{a_0} \leq \varepsilon_0$ (again, we abuse our notation and view ${\rm Proj}_{e_n^\perp}(x)$ as a vector with $n-1$ coordinates). It remains to recall that, since $(\intr K) \cap e_n^\perp\neq \emptyset$ due to our assumptions on $e_n^\perp = H$, Lemma \ref{lem:affine-set} gives us that ${\rm relint}(K_{e_n}) = {\rm relint}(K\cap e_n^\perp) \subset \intr K$.
\end{proof}

\smallskip

\begin{remark}\label{rem:cubelike-balanced}
Let $K$ be a 1-unconditional convex body in $\R^n$ which has Property $(\dagger)$. It is not hard to see that we can find at least one affine image $\widetilde{K}$ of $K$ which belongs to the subclass $\U^n$ and still has Property $(\dagger)$.

Indeed, as mentioned in Section \ref{sec:prelims}, an obvious choice for an affine image of $K$ from $\U^n$ is the convex body we get if we multiply $K$ by the diagonal matrix ${\rm diag}\bigl(\|e_1\|_K^{-1},\|e_2\|_K^{-1},\ldots,\|e_n\|_K^{-1}\bigr)$. Let us write $T_0$ for the linear transformation of $\R^n$ which corresponds to this matrix. It is clear that
\begin{equation*}
x\in {\rm ext}(K) \quad \Leftrightarrow \quad T_0(x)\in {\rm ext}(T_0(K)),
\end{equation*}
while this transformation $T_0$ is such that, for all $y\in \R^n$ and $i\in [n]$, $y_i\neq 0\,\Leftrightarrow\,(T_0(y))_i\neq 0$.
\end{remark}

\begin{lemma}\label{lem:cubelike-maximal-proj}
Let $\B$ be a cubelike $1$-unconditional convex body in $\U^n$ (in other words, assume that $\B$ has Property $(\dagger)$, and that $\|e_i\|_\B=1$ for all $i\in [n]$). 

Suppose that for some $i_0\in [n]$ we have that ${\rm Proj}_{e_{i_0}^\perp}(\B)$ is a(n) ($(n-1)$-dimensional) parallelepiped. Then we will have that ${\rm Proj}_{e_{i_0}^\perp}(\B)$ coincides with ${\rm Proj}_{e_{i_0}^\perp}([-1,1]^n)$ (and not only that they are affinely equivalent). Equivalently, $\B$ must contain the point ${\bm 1} - e_{i_0}$.
\end{lemma}
\begin{proof}
WLOG we can assume that $i_0=n$. By the linearity of projections, we can observe that
\begin{equation*}
{\rm ext}\bigl({\rm Proj}_{e_n^\perp}(\B)\bigr) \subseteq \bigl\{{\rm Proj}_{e_n^\perp}(x): x\in {\rm ext}(\B)\bigr\}.
\end{equation*}
Moreover, since we have assumed that ${\rm Proj}_{e_n^\perp}(\B)$ is an $(n-1)$-dimensional parallelepiped, we know that it has exactly $2^{n-1}$ extreme points. Let $v_1$ be one of them; as already observed, we can find an extreme point $x_1$ of $\B$ such that $v_1 = {\rm Proj}_{e_n^\perp}(x_1)$. 

Now note that, because of Property $(\dagger)$, all coordinates of $x_1$ are non-zero, and hence the first $n-1$ coordinates of $v_1$ will be non-zero. Moreover, by the 1-unconditionality, we know that all coordinate reflections of $v_1$ must also be extreme points of ${\rm Proj}_{e_n^\perp}(\B)$, and given what we just remarked, we have that there are $2^{n-1}$ different such coordinate reflections (including $v_1$ itself).

We conclude that the extreme points of ${\rm Proj}_{e_n^\perp}(\B)$ are precisely $v_1 = {\rm Proj}_{e_n^\perp}(x_1)$ and its coordinate reflections. WLOG we can assume that $v_1$ has only positive coordinates (except for its last one).

We finally observe that, for all $j\in [n-1]$, $e_j\in {\rm Proj}_{e_n^\perp}(\B)$, and thus it must be possible to write it as a convex combination of $v_1$ and its coordinate reflections. This implies that $|v_{1,j}| = v_{1,j} \geq 1$ (and since $\B\in \U^n$, we also have $v_{1,j} = x_{1,j} \leq 1$). We conclude that $v_1 = {\rm Proj}_{e_n^\perp}(x_1) = {\bm 1}-e_n$, and by the 1-unconditionality we know that this is contained in $\B$.
\end{proof}

We are now ready to give the

\noindent \emph{Proof of Theorem \ref{thm:main-cubelike-bodies}.} Let $n\geq 3$, and consider a cubelike 1-unconditional convex body $\B$ in $\R^n$ which is NOT a parallelepiped. Because of Remark \ref{rem:cubelike-balanced}, we can also assume that $\B\in \U^n$ without ruining any of the other assumptions. We will show that $\II(\B) \leq 2^n -2$ by using induction in the dimension $n$.

{\bf Base of induction:} If $n=3$, then Theorem \ref{thm:R^3-summary} gives us that $\II(\B) \leq 6$ (without even having to use the assumption that $\B$ is cubelike), and it also guarantees that we can illuminate (any affine image of) $\B$ using 3 pairs of opposite directions.

\smallskip

{\bf Induction Step:} Let us now assume that the theorem holds in dimension $n - 1$ for some $n > 3$, and consider $\B\in \U^n$ as above. 

Because of the 1-unconditionality, we have that ${\rm Proj}_{e_n^\perp}(\B) = \B\cap e_n^\perp$. Moreover, as we also observed in the proof of Lemma \ref{lem:cubelike-maximal-proj}, it holds that
\begin{equation*}
{\rm ext}\bigl({\rm Proj}_{e_n^\perp}(\B)\bigr) \subseteq \bigl\{{\rm Proj}_{e_n^\perp}(x): x\in {\rm ext}(\B)\bigr\},
\end{equation*}
and thus ${\rm Proj}_{e_n^\perp}(\B)$ (viewed as a convex body in $\R^{n-1}$) is 1-unconditional and cubelike.
\begin{itemize}
\item[$\blacklozenge$] If ${\rm Proj}_{e_n^\perp}(\B)$ is a parallelepiped of $\R^{n-1}$, then Lemma \ref{lem:cubelike-maximal-proj} gives us that $\B$ contains the point ${\bm 1}-e_n$. Note that $\B$ is assumed to NOT be a parallelepiped of $\R^n$, and thus Theorem \ref{thm:main-unit-subcubes} applies in this case, allowing us to illuminate $\B$ using $\frac{1}{2}(2^n-2)$ pairs of opposite directions.
\item[$\blacklozenge$] If instead ${\rm Proj}_{e_n^\perp}(\B)$ is NOT a parallelepiped of $\R^{n-1}$, then we can invoke the inductive hypothesis and conclude that $\II({\rm Proj}_{e_n^\perp}(\B)) \leq 2^{n-1} -2$. In fact, we obtain that ${\rm Proj}_{e_n^\perp}(\B)$ can be illuminated using $2^{n-2}-1$ pairs of opposite directions.

We can then combine this with Lemma \ref{lem:proj-sec} (given our assumptions, which imply the conditions of that lemma), and this gives that $\II(\B) \leq 2\cdot \II({\rm Proj}_{e_n^\perp}(\B))\leq 2^n - 4$. 

In addition, by looking at the proof of Lemma \ref{lem:proj-sec} as well, we can check that, if we start with an illuminating set of ${\rm Proj}_{e_n^\perp}(\B)$ which consists of pairs of opposite directions, then we pass to an illuminating set for $\B$ which also consists of pairs of opposite directions (and has double cardinality).
\end{itemize}
This completes the proof. \qed

\medskip

It is worth remarking that we only looked at the particular hyperplane projection ${\rm Proj}_{e_n^\perp}(\B)$ of $\B$ for simplicity. By the 1-unconditionality assumption (which is a rather strong symmetry assumption from certain points of view), we have that ${\rm Proj}_{e_j^\perp}(\B) = \B\cap e_j^\perp$ for all $j\in [n]$. Therefore, we could have applied the inductive hypothesis and Lemma \ref{lem:proj-sec}, exactly as we did above, for any index $j_0\in [n]$ for which we would know that ${\rm Proj}_{e_{j_0}^\perp}(\B)$ is NOT a parallelepiped. And if it turned out that no such index exists, then (by also recalling Lemma \ref{lem:cubelike-maximal-proj}) we would deduce that we are in the setting of Proposition \ref{prop:all-n-1-tuples}, which is just a special case of Theorem \ref{thm:main-unit-subcubes} (in fact, the case with the easiest, most direct proof). Thus, the proof of Theorem \ref{thm:main-cubelike-bodies} only truly requires Proposition \ref{prop:all-n-1-tuples}.

Still, in this paper we sought to give a full proof of Theorem \ref{thm:main-unit-subcubes}, because this allows us to settle more high-dimensional cases of the Illumination Conjecture in the class of 1-unconditional bodies.

\section{1-unconditional convex bodies in $\R^4$}\label{sec:4-dim-bodies}

It remains to complete the proof of Theorem \ref{thm:main-dim3-dim4}: note that we have already fully settled the part of the theorem that concerns $\R^3$. Recall also that Propositions \ref{prop:all-n-1-tuples}, \ref{prop:all-but-one-n-1-tuples}, \ref{prop:exactly-one-n-1-tuple} and \ref{prop:exactly-two-n-1-tuples} apply with $n=4$ too, and cover all cases where we have at least one coordinate permutation of $e_1+e_2+e_3$ contained in a 1-unconditional convex body $\B\in \mathcal{U}^4$. Analogously, Proposition \ref{prop:all-n-2-tuples} corresponds to, and settles, the case of 1-unconditional convex bodies in $\mathcal{U}^4$ which contain \underline{\bf all} coordinate permutations of $e_1+e_2$ and no coordinate permutations of $e_1+e_2+e_3$. 

Therefore, to also fully confirm the Illumination Conjecture for 1-unconditional convex bodies in $\R^4$, it remains to address the cases where $\B\in \U^4$ contains only some of the coordinate permutations of $e_1+e_2$ or none of them. We summarise the conclusions that we reach in this section in the following

\smallskip

\begin{theorem}\label{thm:R^4-summary}
Let $\B\in \U^4$ which is not a parallelepiped. 
\begin{itemize}
\item If $\B$ contains at least one coordinate permutation of $e_1+e_2+e_3$, then, as we have already seen, there exist $\delta\in (0,1)$, or $\delta_1\in (0,1)$, or $\epsilon_2$ and $\delta_{\epsilon_2}\in (0,1)$, or $\delta_3,\widetilde{\delta}_3\in (0,1)$ and $\eta_3\in (0,\delta_3)$, which depend only on $\B$, so that $\B$ is illuminated by one of the following sets:
\begin{gather*}
{\cal F}_{\ref{prop:all-n-1-tuples},\delta}^4= \I^4_{-2}(\delta) =\bigl\{\pm(1,\delta,\delta,\delta),\,\pm(\delta,-1,-\delta,-\delta),\,\pm(\delta,\delta,-1,-\delta),\,\pm(\delta,-\delta,1,\delta),
\\
\phantom{\I^4_{-2}(\delta)=(1)} \pm(\delta,\delta,\delta,-1),\,\pm(\delta,-\delta,-\delta,1),\,\pm(\delta,-\delta,\delta,-1)\bigr\},
\\[0.4em]
\hbox{or}\ \, {\cal F}_{\ref{prop:all-but-one-n-1-tuples},\delta}^4 = \I^3_{-2}(\delta)\times \{\pm\delta\} = \bigl\{\pm(1,\delta,\delta),\,\pm(\delta,-1,-\delta),\,\pm(\delta,\delta,-1)\bigr\}\times\{\pm \delta\},
\\[0.5em]
\hbox{or}\ \, {\cal F}_{\ref{prop:exactly-two-n-1-tuples},\delta_1}^{4,1} = \I_{ex2,1}^4(\delta_1) = \bigl\{\pm(\delta_1,\delta_1,\delta_1,1),\,\pm(\delta_1,\delta_1,\delta_1,-1),\,\pm(\delta_1,-\delta_1,-\delta_1,1),
\\
\phantom{{\cal F}_{\ref{prop:exactly-two-n-1-tuples},\delta_1}^{4,1}= \I_{ex2,1}^4(\delta_1) = (1,1,1)} \,\pm\!(\delta_1,-\delta_1,-\delta_1,-1), \,\pm(\delta_1,\delta_1,-1,0),\,\pm(\delta_1,-\delta_1,1,0)\bigr\},
\\[0.5em]
\hbox{or}\ \, {\cal F}_{\ref{prop:exactly-two-n-1-tuples},\epsilon_2,\delta_{\epsilon_{_{\tiny 2}}}}^{4,2} = \I_{ex2,2}^4(\epsilon_2,\delta_{\epsilon_{_{\tiny 2}}}) = \bigl\{\pm(\epsilon_2,\epsilon_2,1,1),\,\pm(\delta_{\epsilon_{_{\tiny 2}}},-\delta_{\epsilon_{_{\tiny 2}}},-1,-\delta_{\epsilon_{_{\tiny 2}}}),\,\pm(\delta_{\epsilon_{_{\tiny 2}}},\delta_{\epsilon_{_{\tiny 2}}},-1,-\delta_{\epsilon_{_{\tiny 2}}}),
\\
\pm(\delta_{\epsilon_{_{\tiny 2}}},-\delta_{\epsilon_{_{\tiny 2}}},1,\delta_{\epsilon_{_{\tiny 2}}}),\,\pm(\delta_{\epsilon_{_{\tiny 2}}},\delta_{\epsilon_{_{\tiny 2}}},\delta_{\epsilon_{_{\tiny 2}}},-1),\,\pm(\delta_{\epsilon_{_{\tiny 2}}},-\delta_{\epsilon_{_{\tiny 2}}},-\delta_{\epsilon_{_{\tiny 2}}},1),\,\pm(\delta_{\epsilon_{_{\tiny 2}}},-\delta_{\epsilon_{_{\tiny 2}}},\delta_{\epsilon_{_{\tiny 2}}},-1)\bigr\},
\\[0.6em]
\hbox{or}\ \, {\cal F}_{\ref{prop:exactly-one-n-1-tuple},\delta_3,\eta_3,\widetilde{\delta}_3}^4 = \bigl\{\bigl(\pm(1,\eta_3),\pm\delta_3,0\bigr),\, \bigl(\pm(-\eta_3,1),\pm \delta_3,0\bigr)\bigr\}\,\cup\,\bigl[\bigl\{\pm(\widetilde{\delta}_3,\widetilde{\delta}_3,\widetilde{\delta}_3)\bigr\}\times\{\pm 1\}\bigr],
\end{gather*}
or by some coordinate permutation of those.
\item Again, as we have seen, if $\B$ contains all coordinate permutations of $e_1+e_2$ (and none of the coordinate permutations of $e_1+e_2+e_3$), then $\B$ can be illuminated by the set
\begin{align*}
{\cal F}_{\ref{prop:all-n-2-tuples},\delta,\zeta}^4 = \widetilde{\I}_{\delta,\zeta}^4 = \bigl\{\!&\pm(\delta, 1,\zeta,\delta),\ \pm(\delta,-1,-\delta,-\delta),\ \,\pm(\delta,\delta,-1,-\zeta),\ \pm(\delta,-\delta, 1,\delta),
\\
&\pm(\delta,\delta,\delta,-1),\ \pm(\delta,-\zeta,-\delta,1),\ \ \pm(-\zeta,\delta,-\delta,1)\bigr\}
\end{align*}
for some $\delta\in (0,1)$ and $\zeta\in (0,\delta)$ which depend only on $\B$. Note that an equivalent description for $\B$ here is that it contains all possible unit squares.
\item In the remaining 4-dimensional cases, where $\B$ may contain only some unit squares, or none at all, we can use a coordinate permutation of one of the following illuminating sets:
\begin{align*}
{\cal F}_{\ref{prop:R^4-no-pairs},\ref{prop:R^4-exactly-five-pairs},\delta, \eta,\zeta}:= \bigl\{\!&\pm(1,\delta,\eta,0),\ \pm(\delta,-1,-\eta,0),\ \pm(\delta,-\eta,1,\zeta),\ \pm(-\delta,-\eta,1,\zeta),
\\
&\quad \pm(0,\pm(\eta,\delta),1),\ \pm(0,1,-\delta,\eta)\bigr\},
\end{align*}
\begin{align*}
{\cal F}_{\ref{prop:R^4-exactly-one-pair},\ref{prop:R^4-exactly-two-pairs},\delta, \eta}:= \bigl\{\!&\pm(1,\delta,\eta,0),\ \pm(\delta,-1,\eta,0),\ \pm(-\eta,\delta,1,0),\ \pm(\eta,\delta,-1,0),
\\
&\pm(-\eta,0,\delta,1), \ \pm(-\eta,0,\delta,-1),\ \pm(1,0,1,0)\bigr\},
\end{align*}
\begin{equation*}
{\cal F}_{\ref{prop:R^4-exactly-two-pairs}, {\rm alt}, \delta_1} := \bigl\{\pm(1,\delta_1,0,0),\pm(-\delta_1,1,0,0),\pm(0,0,1,\delta_1),\pm(0,0,-\delta_1,1)\bigr\},
\end{equation*}
\begin{align*}
{\cal F}_{\ref{prop:R^4-exactly-three-pairs},\ref{prop:R^4-exactly-four-pairs},\delta, \eta}:= \bigl\{\!&\pm(1,\delta,\eta,0),\ \pm(\delta,-1,-\eta,0),\ \pm(\delta,-\eta,1,0),\ \pm(\delta,\eta,-1,0),
\\
&\quad \pm(0,\pm(\eta,\delta),1),\ \pm(0,1,-1,0)\bigr\},
\end{align*}
\begin{gather*}
{\cal F}_{\ref{prop:R^4-exactly-three-pairs}, {\rm alt}, \delta_2,\eta_2} := \bigl\{\pm(\eta_2,1,\delta_2,0),\ \pm(-\eta_2,1,\delta_2,0),\ \pm(-\eta_2,1,-\delta_2,0),\qquad \qquad \qquad
\\
\qquad\quad \pm(\eta_2,\eta_2,1,\delta_2),\ \pm(-\eta_2,-\eta_2,1,\delta_2),\ \pm(\eta_2,\eta_2,-\delta_2,1),\ \pm(-\eta_2,-\eta_2,-\delta_2,1)\bigr\},
\end{gather*}
\begin{align*}
\hbox{or}\ \ {\cal F}_{\ref{prop:R^4-exactly-four-pairs}, {\rm alt}, \delta_3,\eta_3} := \bigl\{\!&\pm(1,-\eta_3,-\delta_3,-\delta_3),\ \pm(-\eta_3,1,-\delta_3,-\delta_3),
\\&\pm(\delta_3,0,1,-\eta_3),\ \pm(\delta_3,0,-\eta_3,1),\  \pm(0,\delta_3,1,-\eta_3),\ \pm(0,\delta_3,-\eta_3,1)\bigr\}.
\end{align*}
\end{itemize}
\end{theorem}

Note that the order in which we include the parameters/`small' constants that we use as subscripts should indicate how they (potentially) depend on each other: parameters appearing earlier do not depend on later ones, but how small one may need to choose the later one(s) depends on the values of the earlier parameters.

In the same manner as in earlier sections, we divide the remaining cases into the following propositions (based on how many unit squares the given unconditional body contains).

\begin{proposition}
\label{prop:R^4-no-pairs}
Suppose that $\B \in \mathcal{U}^4$ satisfies $\norm{e_i+e_j}_\B > 1$ for every $i,j \in [4]$. Then there exist $\delta > 0$, $\eta=\eta_\delta>0$ and $\zeta=\zeta_{\delta,\eta}>0$ so that $\B$ can be illuminated by 
some coordinate permutation of 
the set
\begin{align*}
{\cal F}_{\ref{prop:R^4-no-pairs},\ref{prop:R^4-exactly-five-pairs},\delta, \eta,\zeta}:= \bigl\{\!&\pm(1,\delta,\eta,0),\ \pm(\delta,-1,-\eta,0),\ \pm(\delta,-\eta,1,\zeta),\ \pm(-\delta,-\eta,1,\zeta),
\\
&\quad \pm(0,\pm(\eta,\delta),1),\ \pm(0,1,-\delta,\eta)\bigr\}.
\end{align*}
Observe that $\abs{{\cal F}_{\ref{prop:R^4-no-pairs},\ref{prop:R^4-exactly-five-pairs},\delta, \eta,\zeta}} = 14$.
\end{proposition}
\begin{proof} For all $i,j\in [4], i\neq j$, set $\theta_{i,j}= \norm{e_i+e_j}_\B^{-1}$; this is equivalent to $\theta_{i,j}(e_i+e_j)\in \partial \B$, and thus $\theta_{i,j} < 1$ for all $i\neq j$ by our main assumption. Set $\Theta_0 = \max_{i<j}\theta_{i,j}$. 

Next, similarly to the proof of Proposition \ref{prop:R^3-no-pairs}, for each $i\in [4]$ and for each $j\in [4]\backslash \{i\}$, set $\alpha_{i,j}$ to be the supremum of non-negative numbers $x_j$ such that
\begin{equation*}
\frac{1+\Theta_0}{2}e_i + x_j e_j \in \B.
\end{equation*}
Then $\frac{1+\Theta_0}{2}e_i+\alpha_{i,j}e_j\in\B$, and we must have $\alpha_{i,j} <1$ \textit{(in fact, $\alpha_{i,j}$ must be $<\theta_{i,j}$, since otherwise the point $\theta_{i,j}(e_i+e_j)$ would not be a boundary point of $\B$; indeed, since $\theta_{i,j} < 1$, if we had that $\frac{1+\Theta_0}{2}e_i + \theta_{i,j} e_j\in \B$, we could use Lemma \ref{lem:affine-set} to conclude that $y_ie_i+\theta_{i,j}e_j\in \intr\B$ for any $y_i\in (0,\frac{1+\Theta_0}{2})$).} 

Set $\alpha_0=\max_{1\leq i\neq j\leq 4}\alpha_{i,j}$ and WLOG assume that (at least) one of $\alpha_{1,2},\, \alpha_{2,1}$ is equal to $\alpha_0$. We will now show that ${\cal F}_{\ref{prop:R^4-no-pairs},\ref{prop:R^4-exactly-five-pairs},\delta, \eta,\zeta}$ (with suitably chosen $\delta,\eta,\zeta$) illuminates $\B$. Let $x\in \partial\B$. Consider the following cases for the index set $\cZ_x$ of zero coordinates of $x$.
\smallskip
\begin{itemize}
\item[$\abs{\cZ_x}\geq 2$.] First of all, if $x=x_ie_i$ for some $i\in [4]$, where $x_i=\pm 1$, then we find a direction $d\in {\cal F}_{\ref{prop:R^4-no-pairs},\ref{prop:R^4-exactly-five-pairs},\delta, \eta,\zeta}$ satisfying $m.c.(d)=i$ and $d_i\cdot x_i < 0$. Observe that $\|d-d_ie_i\|_\infty = \delta$, and thus, as long as we choose $\zeta < \eta < \delta < \frac{1}{4}$, we will have that $x+d\in \intr\B$.

\smallskip

Similarly, if $x=x_ie_i + x_je_j$, and we assume WLOG that $|x_i|\geq |x_j| > 0$, we will have that $|x_j|\leq \theta_{i,j}\leq \Theta_0$. Again we choose $d\in {\cal F}_{\ref{prop:R^4-no-pairs},\ref{prop:R^4-exactly-five-pairs},\delta, \eta,\zeta}$ which satisfies $m.c.(d)=i$ and $d_i\cdot x_i < 0$. Then, for the displaced vector $x+|x_i|d$, we will have that
\begin{itemize}
\item[$\cdot$] $|(x+|x_i|d)_i|= 0$,
\item[$\cdot$] $|(x+|x_i|d)_j|\leq |x_j| + |x_i|\delta \leq |x_j|+\delta\leq \Theta_0+\delta$,
\item[$\cdot$] and for $s\in [4]\setminus\{i,j\}$, $\,|(x+|x_i|d)_s|\leq \delta$.
\end{itemize}
Thus, if we also choose $\zeta < \eta < \delta < \frac{1-\Theta_0}{3}$, we will have that $x+|x_i|d\in \intr\B$ \textit{(note that these restrictions on $\delta, \eta$ and $\zeta$ do not depend on what the coordinates of $x$ are).}
\medskip
\item[$\abs{\cZ_x} = 0$.] We consider two subcases here.
\begin{itemize}
\item[--] First, assume that $|x_1| \leq \frac{1+\Theta_0}{2}$. Then $|x_1| < 1$, and we can invoke Corollary \ref{cor:affine-set} to conclude that a direction $d$ from
\begin{equation*}
\pm(0,\pm(\eta,\delta),1),\ \ \pm(0,1,-\delta,\eta)
\end{equation*}
will illuminate $x$ if it holds that $d_s\cdot x_s < 0$ for all $s\in [4]\setminus\{1\}$. Otherwise, if none of these directions can satisfy this requirement, it will mean that $\sign(x_2)=-\sign(x_3)=-\sign(x_4)$, in which case one of the directions
\begin{equation*}
\pm(\delta,-\eta,1,\zeta),\ \  \pm(-\delta,-\eta,1,\zeta)
\end{equation*}
will illuminate $x$ (the one direction $d$ which also satisfies $d_1\cdot x_1 < 0$).
\smallskip
\item[--] Now we assume that $|x_1| > \frac{1+\Theta_0}{2}$. Then, by our choice of the numbers $\alpha_{i,j}$, we will have that $|x_4| \leq \alpha_{1,4}\leq \alpha_0 < 1$. 

If $\sign(x_1)=\sign(x_2)=\sign(x_3)$ or $\sign(x_1)=-\sign(x_2)=-\sign(x_3)$, then one of $\pm(1,\delta,\eta,0),\ \pm(\delta,-1,-\eta,0)$ will illuminate $x$ (by Corollary \ref{cor:affine-set}). Thus, assume for the remainder of this case that $x$ does not satisfy either of these sign assumptions. We consider further subcases here.
\begin{itemize}
\item[$\cdot$] If $\max(|x_2|,|x_3|) < \frac{1-\alpha_0}{8}\leq\frac{1-\alpha_{1,4}}{8}$, then, as long as we also make sure that $\zeta < \eta <\delta<\frac{1-\alpha_0}{8}$, we can choose a direction $d\in \bigl\{\pm(1,\delta,\eta,0)\bigr\}$ which satisfies $d_1\cdot x_1 < 0$ (and we can check that $x+|x_1|d\in \intr\B$ by comparing to the point $\frac{1+\alpha_0}{2} e_4 + \frac{1-\alpha_0}{4}(e_2+e_3)\in \B$).
\medskip
\item[$\cdot$] Assume now that $\max(|x_2|,|x_3|) \geq \frac{1-\alpha_0}{8} \geq \min(|x_2|,|x_3|)$. Also assume first that $\max(|x_2|,|x_3|)=|x_2|$. Then we choose the unique direction $d$ from $\bigl\{\pm(1,\delta,\eta,0),$ $\,\pm(\delta,-1,-\eta,0)\bigr\}$ which satisfies $d_s\cdot x_s < 0$ for $s\in [2]$ (recall that, by our last assumption on the signs of $x$, we will also have here that $d_3\cdot x_3 > 0$). For the displaced vector $x+\frac{1-\alpha_0}{8}d$ we observe that
\begin{center}
$\big|\big(x+\frac{1-\alpha_0}{8}d\big)_s\big| < (1-\lambda_0)|x_s|$ for $s\in \{1,2\}$, as long as 
\smallskip\\
$\lambda_0 < \frac{1-\alpha_0}{8}\delta$.
\end{center}
Based on this, we can compare with the vector
\begin{equation*}
u_{\alpha_0}:=(1-\lambda_0)\big|\vec{x}\big| + \lambda_0\left(\frac{3+\alpha_0}{4}e_4 + \frac{1-\alpha_0}{4} e_3\right)\in \B.
\end{equation*}
We have that $\big|\big(x+\frac{1-\alpha_0}{8}d\big)_s\big| < u_{\alpha_0,s}$ for $s\in \{1,2,4\}$, and also that
\begin{equation*}
u_{\alpha_0,3} = |x_3| + \lambda_0\left(\frac{1-\alpha_0}{4}-|x_3|\right) \geq |x_3|+\lambda_0\frac{1-\alpha_0}{8}.
\end{equation*}
Thus, if we also choose
\begin{equation*}
\eta < \lambda_0 < \frac{1-\alpha_0}{8}\delta,
\end{equation*}
we will have that $\big|\big(x+\frac{1-\alpha_0}{8}d\big)_3\big| \leq |x_3|+\frac{1-\alpha_0}{8}\eta < u_{\alpha_0,3}$, which finally shows that $x+\frac{1-\alpha_0}{8}d\in \intr\B$.

\bigskip

On the other hand, if $\max(|x_2|,|x_3|)=|x_3|$, then we pick the unique direction $d^\prime$ from $\bigl\{\pm(\delta,-\eta,1,\zeta),\,\pm(-\delta,-\eta,1,\zeta)\bigr\}$ which satisfies $d^\prime_s\cdot x_s< 0$ for $s\in \{1,3\}$. Similarly, we compare the displaced vector $x+\frac{1-\alpha_0}{8}d^\prime$ with the vector
\begin{equation*}
u^\prime_{\alpha_0} = (1-\lambda_1)\big|\vec{x}\big| + \lambda_1\left(\frac{7+\alpha_0}{8}e_4 + \frac{1-\alpha_0}{8} e_2\right)\in \B,
\end{equation*}
where $\lambda_1$ needs to be $< \frac{1-\alpha_0}{8}\delta$.
Observe that, as long as $\zeta < 7\lambda_1 < \frac{7(1-\alpha_0)}{8}\delta$, we will have that $x+\frac{1-\alpha_0}{8}d^\prime\in \intr\B$.
\medskip
\item[$\cdot$] Finally, assume that $\min(|x_2|,|x_3|)> \frac{1-\alpha_0}{8}$. Then we pick the unique direction $d\in \bigl\{\pm(\delta,-\eta,1,\zeta),\,\pm(-\delta,-\eta,1,\zeta)\bigr\}$ which satisfies $d_s\cdot x_s < 0$ for $s\in [3]$. Again, we can check that $x+\frac{1-\alpha_0}{8}d\in\intr\B$ as long as $\zeta < (1-\alpha_0)\eta$.
\end{itemize}
\end{itemize}
\end{itemize}

We now turn to the cases where:
\begin{itemize}
\item[$\abs{\cZ_x} = 1$.] We will argue similarly to the previous case, but will now rely on the existence of the points
\begin{equation*}
v_{i;\Theta_0}:= \frac{3+\Theta_0}{4}e_i + \frac{1-\Theta_0}{12}\sum_{j\neq i}e_j
\end{equation*}
in $\B$.
\begin{itemize}
\item[--] Assume first of all that $x_4=0$. 
\begin{itemize}
\item[$\cdot$] If $|x_3| \leq \frac{1+\Theta_0}{2}$, then we cannot have $\max(|x_1|,|x_2|)<\frac{1-\Theta_0}{24}$ \textit{(because otherwise we could compare the entries of $x$ with those of $v_{3;\Theta_0}$ and we would obtain that $x$ is not a boundary point).} If $|x_2|=\min(|x_1|,|x_2|)$ and it is $<\frac{1-\Theta_0}{24}$, then we simply choose a direction from $\pm(1,\delta,\eta,0)$ such that $d_1\cdot x_1 < 0$. Then $x+|x_1|d\in \intr\B$ as long as \[\zeta<\eta<\delta<\frac{1-\Theta_0}{24}\] (this can be seen by comparing with the corresponding coordinates of the point $v_{3;\Theta_0}$ again).

Similarly, if $|x_1|=\min(|x_1|,|x_2|)<\frac{1-\Theta_0}{24}$, then we illuminate $x$ with a direction from $\pm(\delta,-1,-\eta,0)$. 

\medskip

Finally, if $\min(|x_1|,|x_2|)\geq \frac{1-\Theta_0}{24}$, then we choose $d$ from the same directions as previously, $\pm(1,\delta,\eta,0),\ \pm(\delta,-1,-\eta,0)$, so that $d_s\cdot x_s < 0$ for $s\in [2]$. We compare the coordinates of the displaced vector $x+\frac{1-\Theta_0}{24}d$ with those of a convex combination of the form
\begin{equation}\label{eqp1:prop:R^4-no-pairs}
\lambda e_3+(1-\lambda)\big|\vec{x}\big| = \bigl((1-\lambda)|x_1|,\,(1-\lambda)|x_2|,\,|x_3|+\lambda(1-|x_3|),\,0\bigr)
\end{equation}
where $\lambda < \frac{1-\Theta_0}{24}\delta$. But then, as long as $\frac{1-\Theta_0}{24}\eta < \lambda\frac{1-\Theta_0}{2}\leq \lambda(1-|x_3|)$, which can be achieved if $\eta < \frac{1-\Theta_0}{2}\delta$, we will have that $x+\frac{1-\Theta_0}{24}d\in\intr\B$.
\medskip
\item[$\cdot$] If instead $|x_3| > \frac{1+|\Theta_0}{2}$, then we observe the following: by 1-unconditionality, we will have that $|x_1|e_1 + |x_3|e_3\in \B\,\Rightarrow\,|x_1|\leq \alpha_{3,1}\leq\alpha_0=\max(\alpha_{1,2},\alpha_{2,1})$. Similarly, we see that $|x_2|\leq \alpha_{3,2}\leq \max(\alpha_{1,2},\alpha_{2,1})$. Therefore, we can pick the unique direction $d\in \{\pm(\delta,-\eta,1,\zeta)\}$ which satisfies $d_3\cdot x_3 < 0$, and we will have that $x+|x_3|d\in \intr\B$ as long as \[\eta <\delta < \frac{\theta_{1,2}-\alpha_0}{2} = \frac{\theta_{1,2}-\max(\alpha_{1,2},\alpha_{2,1})}{2}\quad \hbox{and}\quad \zeta < \frac{1}{2}\left(1-\frac{\alpha_0}{\theta_{1,2}}\right)\] (we can confirm this by comparing to the point \[\frac{\theta_{1,2}+\alpha_0}{2}(e_1+e_2) + \left(\frac{1}{2}-\frac{\alpha_0}{2\theta_{1,2}}\right)e_4\]
which is a convex combination of $\theta_{1,2}(e_1+e_2)$ and of $e_4$).
\end{itemize}
\item[--] Now we assume that $x_3=0$. 
\begin{itemize}
\item[$\cdot$] If $|x_4|\leq \frac{1+\Theta_0}{2}$, we argue exactly as before, and we illuminate $x$ using one of the directions $\pm(1,\delta,\eta,0),\ \pm(\delta,-1,-\eta,0)$. The only change that we make is that, in the subcases where $\min(|x_1|,|x_2|)\geq \frac{1-\Theta_0}{24}$, instead of comparing the displaced vector $x+\frac{1-\Theta_0}{24}d$ with a vector analogous to the one in \eqref{eqp1:prop:R^4-no-pairs}, we compare with a vector of the form
\begin{equation*}
\lambda^\prime v_{4;\Theta_0} + (1-\lambda^\prime)\big|\vec{x}\big|.
\end{equation*}
\item[$\cdot$] If $|x_4| > \frac{1+\Theta_0}{2}$, again we argue similarly to the previous case, while illuminating $x$ with the unique direction $d\in \{\pm(0,\eta,\delta,1)\}$ which satisfies $d_4\cdot x_4 < 0$.
\end{itemize}
\item[--] Next, assume that $x_2=0$. In this case we illuminate $x$ with a direction $d$ from
\begin{equation*}
\pm(1,\delta,\eta,0),\ \ \pm(\delta,-\eta,1,\zeta),\ \ \pm(-\delta,-\eta,1,\zeta),\ \ \pm(0,\pm(\eta,\delta),1)
\end{equation*}
and we distinguish subcases based on whether $|x_4|\leq \frac{1+\Theta_0}{2}$ or not (in fact, in the latter subcase we illuminate $x$ with one of the directions from $\pm(0,\pm(\eta,\delta),1)$, relying on the fact that $|x_1|\leq \alpha_{4,1}\leq \alpha_0 < \frac{1+\Theta_0}{2}$).
\smallskip
\item[--] Finally assume that $x_1=0$. Then we illuminate $x$ with one of the directions from
\begin{equation*}
\pm(0,\pm(\eta,\delta),1),\ \ \pm(0,1,-\delta,\eta),\ \ \pm(\delta,-\eta,1,\zeta).
\end{equation*}
If $\sign(x_2)=\sign(x_3)$, then we use one of the first 4 directions. We also use one of these directions in the cases where $\sign(x_2)=-\sign(x_3)$, but at the same time $|x_2| \leq \frac{1+\Theta_0}{2}$ and $|x_4|\geq \frac{1-\Theta_0}{24}$, whereas we use the last pair of directions here if, instead of the last assumption, we have $\min(|x_3|,|x_4|) = |x_4| < \frac{1-\Theta_0}{24}$.

\smallskip

Finally, if $\sign(x_2)=-\sign(x_3)$ and $|x_2| > \frac{1+\Theta_0}{2}$, then we will have that $|x_4| \leq \alpha_{2,4}\leq \alpha_0$. Thus, regardless of whether $\min(|x_2|,|x_3|) = |x_3| < \frac{1-\Theta_0}{24}$ or not, we will illuminate $x$ by a direction $d\in \{\pm(0,1,-\delta,\eta)\}$ (it's just that in the former subcase we will compare the displaced vector $x+|x_2|d$ with the vector $v_{4;\Theta_0}$, whereas in the latter subcase we will compare the same displaced vector with a convex combination of the form $\lambda e_4 + (1-\lambda)|\vec{x}|$; note that in the latter subcase we again consider the displacement $x+|x_2|d$ because $|x_3|\geq \frac{1-\Theta_0}{24} > \delta$ by our restrictions so far, and thus $|(x+|x_2|d)_3| = |x_3| - |x_2|\delta$, which we can write as $< (1-\lambda)|x_3|$ for a suitably chosen $\lambda\in (0,1)$).
\end{itemize}
\end{itemize}

We have thus illuminated all boundary points of $\B$.
We finally remark that, if $\B$ doesn't satisfy the assumption that $\max_{1\leq i\neq j\leq 4}\alpha_{i,j} =\max(\alpha_{1,2},\alpha_{2,1})$, then clearly a coordinate permutation $\iota$ of $\R^4$ suffices to give an affine image $\iota(\B)$ of $\B$ which does. This completes the proof. 
\end{proof}


\begin{proposition}
\label{prop:R^4-exactly-one-pair}
Suppose that for a given $\B \in \mathcal{U}^4$ there is exactly one pair of indices $i_1,i_2\in [4]$ such that $\norm{e_{i_1}+e_{i_2}}_\B = 1$. Then there exist $\delta > 0$ and $\eta=\eta_\delta>0$ so that $\B$ can be illuminated by some coordinate permutation of the set
\begin{align*}
{\cal F}_{\ref{prop:R^4-exactly-one-pair},\ref{prop:R^4-exactly-two-pairs},\delta, \eta}:= \bigl\{\!&\pm(1,\delta,\eta,0),\ \pm(\delta,-1,\eta,0),\ \pm(-\eta,\delta,1,0),\ \pm(\eta,\delta,-1,0),
\\
&\pm(-\eta,0,\delta,1), \ \pm(-\eta,0,\delta,-1),\ \pm(1,0,1,0)\bigr\}.
\end{align*}
\end{proposition} 
\begin{proof}
WLOG we can assume that $e_1+e_2\in \B$. For any $1\leq i\neq j\leq 4$ with $\{i,j\}\neq \{1,2\}$, we set 
\[\theta_{i,j}:=\|e_i+e_j\|_\B^{-1},\]
which by our assumption will be $<1$. Set $\Theta_0 =\max\bigl\{\theta_{i,j}: 1\leq i\neq j\leq 4, \{i,j\}\neq \{1,2\}\bigr\}$ and fix some 
$\delta \in \bigl(0,\, \frac{1-\Theta_0}{4}\bigr)$ and some $\eta_\delta \in (0, \delta/2)$ (we will further restrict $\eta_\delta$ by the end of the proof). Let $x$ be a boundary point of $\B$, and consider the following cases.
\begin{itemize}
\item[$\abs{\cZ_x}=3$.] Since $\eta < \delta < \frac{1}{4}$, to illuminate a boundary point $x=x_i e_i$ with $x_i=\pm1$, it suffices to pick a direction $d_x$ from the first 12 in ${\cal F}_{\ref{prop:R^4-exactly-one-pair},\ref{prop:R^4-exactly-two-pairs},\delta, \eta}$ which satisfies $m.c.(d_x)=m.c.(x) = i$ and $d_{x,i}\cdot x_i < 0$.
\item[$\abs{\cZ_x}=2$.] We first deal with the subcase where
\begin{itemize}
\item[--] $\cZ_x = \{3,4\}$. Here we pick the unique direction $d\in \pm(1,\delta,\eta,0),\ \pm(\delta,-1,\eta,0)$ which satisfies $d_s\cdot x_s < 0$ for $s\in [2]$. Then the non-zero coordinates of $x+d$ do not exceed in absolute value the non-zero coordinates of either $(1-\delta,0,\eta,0)$ or $(0,1-\delta,\eta,0)$. The latter points are interior points of $\B$ since $1-\delta + \eta < 1$, and thus $x+d\in \intr\B$ too.
\item[--] $\cZ_x\neq\{3,4\}$. If we write $x= x_{i_1}e_{i_1} + x_{i_2}e_{i_2}$ with $\{i_1,i_2\}=[4]\setminus\cZ_x$, and we assume WLOG that $|x_{i_1}|\geq |x_{i_2}|$, then we will have that $|x_{i_2}|\leq \theta_{i_1,i_2}\leq \Theta_0$. Thus, we can pick a direction $d_x$ from the first 12 in ${\cal F}_{\ref{prop:R^4-exactly-one-pair},\ref{prop:R^4-exactly-two-pairs},\delta, \eta}$ so that $m.c.(d_x) = i_1$ and $d_{x,i_1}\cdot x_{i_1} < 0$. We can compare the displaced vector $x+|x_{i_1}|d_x$ to either $(\Theta_0+\delta)e_{i_2}+\delta e_{i_3}$ or $\Theta_0e_{i_2} + \delta(e_{i_3}+e_{i_4})$, where $\{i_3,i_4\}=\cZ_x$. Given that $\Theta_0 + 2\delta < 1$, the latter points are in $\intr\B$, and thus the same is true for $x+|x_{i_1}|d_x$.
\end{itemize}
\item[$\abs{\cZ_x}=0$.] Here we distinguish subcases based on the magnitude of $|x_4|$.
\begin{itemize}
\item[--] If $|x_4|\leq \Theta_0 < 1$, then, by employing Corollary \ref{cor:affine-set} (combined with Corollary \ref{cor:uncond-illum}), we can illuminate $x$ using one of the first 8 directions in ${\cal F}_{\ref{prop:R^4-exactly-one-pair},\ref{prop:R^4-exactly-two-pairs},\delta, \eta}$ (we choose the unique direction $d$ among these which satisfies $d_s\cdot x_s < 0$ for $s\in [3]$).
\item[--] Assume now that $|x_4| > \Theta_0$. Then $|x_1| \leq \theta_{1,4}\leq \Theta_0$ and $|x_2|\leq \theta_{2,4}\leq \Theta_0$.
\begin{itemize}
\item[$\cdot$] If $|x_3| < \frac{1-\Theta_0}{4}$, we use the unique direction $d\in \{\pm(-\eta,0,\delta,1)\}$ which satisfies $d_4\cdot x_4 < 0$. We compare the displaced vector $x+|x_4|d$ with the convex combination
\begin{equation*}
u_{3;\Theta_0}:=\frac{1+\Theta_0}{2}(e_1+e_2) + \frac{1-\Theta_0}{2}e_3\in \B.
\end{equation*}
\item[$\cdot$] If instead $|x_3|\geq \frac{1-\Theta_0}{4}$, then we use the unique direction $d\in \{\pm(-\eta,0,\delta,1),$ $\pm(-\eta,0,\delta,-1)\}$ which satisfies $d_s\cdot x_s < 0$ for $s\in \{3,4\}$. For the same indices $s$, we have that
\begin{equation*}
\big|(x+|x_4|d)_s\big| < (1-\lambda_0)|x_s|
\end{equation*}
as long as $\lambda_0\in (0,\delta)$ (since $(x+|x_4|d)_4=0$, and since $\lambda_0|x_3|\leq \lambda_0\theta_{3,4}\leq \lambda_0\Theta_0$ and $|x_4|\delta \geq \Theta_0\delta$). But then we compare $x+|x_4|d$ with the vector
\begin{equation*}
(1-\lambda_0)\big|\vec{x}\big| + \lambda_0(e_1+e_2)\in \B,
\end{equation*}
and we can conclude that $x+|x_4|d\in \intr\B$ as long as $\eta < \lambda_0(1-\Theta_0)< \delta(1-\Theta_0)$ (so that $\big|(x+|x_4|d)_1\big|\leq |x_1| + |x_4|\eta < |x_1| +\lambda_0(1-\Theta_0) \leq |x_1| + \lambda_0(1-|x_1|)$).
\end{itemize}
\end{itemize}
\item[$\abs{\cZ_x}=1$.] First of all, if $\cZ_x=\{4\}$, then we argue as in the cases where $\abs{\cZ_x}=0$ and $|x_4| \leq \Theta_0$: the first 8 directions of ${\cal F}_{\ref{prop:R^4-exactly-one-pair},\ref{prop:R^4-exactly-two-pairs},\delta, \eta}$ illuminate $x$.

The remaining subcases are the following.
\begin{itemize}
\item[--] Assume that $x_3=0$. If $|x_4|\leq \Theta_0$, while $\min(|x_1|,|x_2|)\leq \frac{1-\Theta_0}{4}\leq \max(|x_1|,|x_2|)$, then we use a direction $d$ from $\pm(1,\delta,\eta,0),\,\pm(\delta,-1,\eta,0)$ so that $m.c.(d)$ is the same index where $\max(|x_1|,|x_2|)$ is attained, and so that $d_{m.c.(d)}\cdot x_{m.c.(d)} < 0$. We will have that $x+|x_{m.c.(d)}|d\in \intr\B$ based on our restrictions on $\delta$ and $\eta$. 

\smallskip

If instead $|x_4|\leq \Theta_0$ and $\min(|x_1|,|x_2|)\geq \frac{1-\Theta_0}{4}$, then we choose from the same directions a direction $d$ such that $d_s\cdot x_s< 0$ for $s\in [2]$. We compare the displaced vector $x+\frac{1-\Theta_0}{4}d$ with a convex combination of the form
\begin{equation*}
(1-\lambda)\big|\vec{x}\big| + \lambda\left(\frac{1+\Theta_0}{2}e_4 + \frac{1-\Theta_0}{2}e_3\right)
\end{equation*}
where $\lambda < \frac{1-\Theta_0}{4}\delta$, to conclude that $x+\frac{1-\Theta_0}{4}d\in \intr\B$ as long as $\eta < \frac{1-\Theta_0}{2}\delta$.

\smallskip

Finally, if $|x_4| > \Theta_0$, then we use a direction $d\in\{\pm(-\eta,0,\delta,1)\}$ and compare $x+|x_4|d$ with $u_{3;\Theta_0} = \frac{1+\Theta_0}{2}(e_1+e_2) + \frac{1-\Theta_0}{2}e_3$ again.
\item[--] Assume now that $x_1=0$. Then we will illuminate $x$ using one of the directions
\begin{equation*}
\pm(\delta,-1,\eta,0),\ \ \pm(-\eta,\delta,1,0),\ \ \pm(\eta,\delta,-1,0),\ \ \pm(-\eta,0,\delta,1), \ \ \pm(-\eta,0,\delta,-1)
\end{equation*}
while distinguishing subcases based on whether $|x_4| \leq \Theta_0$ or not, and whether in the former case $\min(|x_2|,|x_3|) \leq \frac{1-\Theta_0}{4}$ or not, or whether in the latter case $|x_3|\leq \frac{1-\Theta_0}{4}$ or not.
\item[--] Finally, we assume that $x_2=0$. If $|x_4| > \Theta_0$, then we will have that $|x_1| \leq \Theta_0$, and thus we can use a direction $d$ from $\pm(-\eta,0,\delta,1), \,\pm(-\eta,0,\delta,-1)$ so that $d_s\cdot x_s < 0$ for $s\in \{3,4\}$: we can conclude that $x+|x_4|d\in \intr\B$ (while distinguishing subcases in our analysis based on whether $|x_3|\leq \frac{1-\Theta_0}{4}$ or not).

\smallskip

If instead $|x_4| \leq \Theta_0 < 1$, then we can also rely on Corollary \ref{cor:affine-set}. We distinguish cases based on whether $\sign(x_1)=-\sign(x_3)$ or not. In the former case, we use again one of the directions $\pm(-\eta,0,\delta,1), \,\pm(-\eta,0,\delta,-1)$ to illuminate $x$ (here we can find a direction $d$ such that $d_s\cdot x_s < 0$ for all $s\in \{1,3,4\}$). If instead $\sign(x_1) = \sign(x_3)$, then Corollary \ref{cor:affine-set} guarantees that one of the directions $\pm(1,0,1,0)$ illuminates $x$.
\end{itemize}
\end{itemize}
In the end, by examining our analysis more carefully, we can see that the restrictions $\delta\in  \bigl(0,\, \frac{1-\Theta_0}{4}\bigr)$ and $\eta< \frac{1-\Theta_0}{2}\delta$ are sufficient to complete the proof.
\end{proof}


\begin{proposition}\label{prop:R^4-exactly-two-pairs}
Suppose that for a given $\B \in \mathcal{U}^4$ there are exactly two pairs of indices $i_1,i_2\in [4]$ such that $\norm{e_{i_1}+e_{i_2}}_\B = 1$. Then at least one of the following two statements holds:
\smallskip\\
(i) there exist $\delta_1 > 0$ and $\eta_1=\eta_{\delta_1}>0$ so that $\B$ can be illuminated by some coordinate permutation of the set
\begin{align*}
{\cal F}_{\ref{prop:R^4-exactly-one-pair},\ref{prop:R^4-exactly-two-pairs},\delta_1, \eta_1}= \bigl\{\!&\pm(1,\delta_1,\eta_1,0),\ \pm(\delta_1,-1,\eta_1,0),\ \pm(-\eta_1,\delta_1,1,0),\ \pm(\eta_1,\delta_1,-1,0),
\\
&\pm(-\eta_1,0,\delta_1,1), \ \pm(-\eta_1,0,\delta_1,-1),\ \pm(1,0,1,0)\bigr\};
\end{align*}
(ii) there exists $\delta_2 > 0$ so that $\B$ can be illuminated by some coordinate permutation of the set
\begin{equation*}
{\cal F}_{\ref{prop:R^4-exactly-two-pairs}, {\rm alt}, \delta_2} := \bigl\{\pm(1,\delta_2,0,0),\pm(-\delta_2,1,0,0),\pm(0,0,1,\delta_2),\pm(0,0,-\delta_2,1)\bigr\}.
\end{equation*}
\end{proposition}
\begin{proof}
We first deal with the cases where statement (ii) definitely applies. These are the cases where the two pairs of indices $i_1\neq i_2$ and $j_1\neq j_2\in [4]$ for which we have 
$\norm{e_{i_1}+e_{i_2}}_\B = \norm{e_{j_1}+e_{j_2}}_\B = 1$ 
satisfy $\{i_1,i_2\}\cap \{j_1,j_2\}=\emptyset$.

\smallskip

WLOG we can assume $\{i_1,i_2\}=\{1,2\}$ and $\{j_1,j_2\}=\{3,4\}$. For each other pair $(i,j)$ of (distinct) indices, we set $\theta_{i,j}=\|e_i+e_j\|_\B^{-1}$, and we know from our assumptions that $\theta_{i,j} < 1$. Let us set $\Theta_0=\max\{\theta_{i,j}:1\leq i<j\leq 4, e_i+e_j\notin\B\}$, and let us fix $\delta_2 < 1-\Theta_0$.

\smallskip

Consider now a boundary point $x$ of $\B$, and suppose $|x_1|$ is its maximum coordinate (in absolute value); note that this doesn't have to be unique. Necessarily $\max\{|x_3|,|x_4|\}\leq\max\{\theta_{1,3},\theta_{1,4}\}\leq \Theta_0$, and thus $(0,0,x_3,x_4)$ is an interior point of $\B$ (since $e_3+e_4\in \B$). 

If $x_2=0$, then we illuminate $x$ using the direction $d_x=-\sign(x_1)(1,\delta_2,0,0)$: we will have $x+|x_1|d_x = (0,-x_1\delta_2,x_3,x_4)$, which by our assumptions is an interior point of $\B$ (this can be seen by comparing with the point $(1-\Theta_0)e_2+\Theta_0(e_3+e_4)\in \B$).

If $x_2\neq 0$, then we pick instead a direction $d_x$ from $\pm(1,\delta_2,0,0),\pm(-\delta_2,1,0,0)$ which satisfies $d_{x,s}\cdot x_s < 0$ for both $s=1$ and $=2$. Then, using Corollary \ref{cor:affine-set}, we see that $d_x$ illuminates $x$.

\smallskip

We argue analogously if the maximum (in absolute value) coordinate of $x$ is its 2nd or 3rd or 4th one.

\bigskip

Consider now the cases where $\{i_1,i_2\}\cap \{j_1,j_2\}\neq \emptyset$. WLOG let $e_1+e_2, \,e_2+e_3\in \B$. Again let us set $\Theta_0=\max\{\theta_{i,j}:1\leq i<j\leq 4, e_i+e_j\notin\B\}$, where $\theta_{i,j}=\|e_i+e_j\|_\B^{-1}$, and observe that $\Theta_0<1$. 
Fix now some 
\begin{equation*}
\delta_1 < \frac{1-\Theta_0}{4}
\end{equation*}
and suppose also that we have chosen some $\eta_1 < \delta_1/2$ (we will soon see that we need to restrict $\eta_1$ further, but this will be done in an unambiguous manner).

The argument that ${\cal F}_{\ref{prop:R^4-exactly-one-pair},\ref{prop:R^4-exactly-two-pairs},\delta_1, \eta_1}$ illuminates $\B$ is very analogous to that of the previous proposition: let us fix an \underline{extreme} point $x\in \partial \B$.
\begin{itemize}
\item[$\blacklozenge$] If $\abs{\cZ_x} = 3$, then $x$ is an extreme point only if $x=\pm e_4$. Then $x$ is illuminated by the directions $\mp (-\eta_1,0,\delta_1,1)$, as long as $\eta_1,\delta_1 < 1$.
\item[$\blacklozenge$] If $\abs{\cZ_x} = 2$, and $x=\pm e_1+ \pm e_2$ or $\pm e_2+ \pm e_3$, then one of the first 8 directions of ${\cal F}_{\ref{prop:R^4-exactly-one-pair},\ref{prop:R^4-exactly-two-pairs},\delta_1, \eta_1}$ illuminates $x$. The `trickiest' case here is if $x=\pm(e_1-e_2)$. Then we have to use the directions $\mp(\delta_1,-1,\eta_1,0)$: e.g. $(e_1-e_2)+  (-\delta_1,1,-\eta_1,0) = (1-\delta_1,0,-\eta_1,0)$, which is $\in {\rm int}(\B)$ since $1-\delta_1 + \eta_1 < 1$.

Note that all other points $x\in \partial \B$ with $\cZ_x=\{3,4\}$ or $\cZ_x=\{1,4\}$ are not extreme, but in the convex hull of $\pm e_1+ \pm e_2$ and $\pm e_2+ \pm e_3$ (so they are also illuminated by the same first 8 directions). Other subcases that we need to consider here are the following.
\begin{itemize}
\item[$\bullet$] If $x= (x_1,0,x_3,0)$, then by our assumptions $\min(|x_1|,|x_3|)\leq \Theta_0$. If $|x_1|\leq |x_3|$, then we illuminate $x$ using the directions $\pm(-\eta_1,\delta_1,1,0)$: 
\begin{equation*}
x+(-\sign(x_3)|x_3|)(-\eta_1,\delta_1,1,0) = (x_1+x_3\eta_1,-x_3\delta_1,0,0)\in {\rm int}(\B)
\end{equation*}
since $\abs{x_1+x_3\eta_1}\leq \Theta_0 + \frac{1-\Theta_0}{4} < 1- \frac{1-\Theta_0}{2}$, while $\abs{x_3\delta_1}\leq \delta_1 < \frac{1-\Theta_0}{2}$. 

Similarly, when $|x_3|\leq |x_1|$, we illuminate $x$ with the directions $\pm(1,\delta_1,\eta_1,0)$.
\item[$\bullet$] Assume now that $x=x_je_j + x_4e_4$ for some $j\in [3]$. One of the `trickiest' cases here is if $j=2$. Again, we distinguish the subcases $|x_4| \leq |x_2|$ and $|x_2|\leq |x_4|$ (with $\min(|x_2|,|x_4|)\leq \Theta_0$ by our assumptions). In the former we have
\begin{equation*}
x+(-\sign(x_2)|x_2|)(-\delta_1,1,-\eta_1,0) = (x_2\delta_1, 0, x_2\eta_1, x_4)\in {\rm int}(\B)
\end{equation*}
since $\abs{x_2\eta_1} < \abs{x_2\delta_1} < \frac{1-\Theta_0}{2}$ and $|x_4|\leq \Theta_0$, and thus the above displaced vector can be compared with the convex combination $\Theta_0e_4 + \frac{1-\Theta_0}{2}(e_1+e_3)\in \B$.

\smallskip

On the other hand, in the subcase where $|x_2|\leq |x_4|$, we will have
\begin{equation*}
x+ (-\sign(x_4)|x_4|)(-\eta_1,0,\delta_1,1) = (x_4\eta_1,x_2,-x_4\delta_1,0)
\end{equation*}
which is $\in {\rm int}(\B)$ by completely analogous reasoning.

\medskip

The remaining cases, where $x=x_je_j + x_4e_4$ with $j=1$ or $3$, are handled very similarly (and we illuminate $x$ using one of the first 12 directions of ${\cal F}_{\ref{prop:R^4-exactly-one-pair},\ref{prop:R^4-exactly-two-pairs},\delta_1, \eta_1}$).
\end{itemize}
\smallskip
\item[$\blacklozenge$] If $\abs{\cZ_x} = 1$, and we have $\cZ_x=\{4\}$, then we can use one of the first 8 directions to illuminate $x$ (since they contain all possible combinations of signs for the first three coordinates).
\begin{itemize}
\item[$\bullet$] Assume now that $\cZ_x=\{3\}$. If $|x_4|\leq \Theta_0$, then $\max(|x_1|,|x_2|) \geq 1-\Theta_0$ (because otherwise $x$ would not be a boundary point, since $\B$ contains $(1-\Theta_0)(e_1+e_2) + \Theta_0e_4$ and we can use Lemma \ref{lem:affine-set} with the section $\{\xi\in \B:\xi_4=\Theta_0\}$). If in addition $\min(|x_1|,|x_2|) < \frac{1-\Theta_0}{4}$, and we write $s_1$ for the index where the maximum is attained, then we can pick a direction $d$ from $\pm(1,\delta_1,\eta_1,0),\,\pm(\delta_1,-1,\eta_1,0)$ so that $m.c.(d) = s_1\in \{1,2\}$ and $d_{s_1}\cdot x_{s_1} < 0$. We will have that $x+|x_{s_1}|d\in \intr\B$ by our assumptions on $\delta_1$ and $\eta_1$.

\smallskip

If instead $|x_4|\leq \Theta_0$ and $\min(|x_1|,|x_2|)\geq \frac{1-\Theta_0}{4}$, then we pick $d$ from the same directions so that $d_s\cdot x_s < 0$ for both $s=1$ and $=2$. We compare $x+\frac{1-\Theta_0}{4}d$ with a convex combination of the form
\begin{equation*}
(1-\lambda)\big|\vec{x}\big| + \lambda\left(\frac{1+\Theta_0}{2}e_4 + \frac{1-\Theta_0}{2}e_3\right)
\end{equation*}
where $\lambda < \frac{1-\Theta_0}{4}\delta_1$. As long as $\eta_1 < 2\lambda < \frac{1-\Theta_0}{2}\delta_1$, we will have that $x+\frac{1-\Theta_0}{4}d\in \intr\B$.

\medskip

Finally, if $|x_4|> \Theta_0$, then $\max(|x_1|,|x_2|)\leq \Theta_0$. Therefore, we pick a direction $d\in \{\pm(-\eta_1,0,\delta_1,1)\}$ such that $d_4\cdot x_4 < 0$, and we can check that $x+|x_4|d\in \intr\B$ by comparing to the point $\frac{1+\Theta_0}{2}(e_1+e_2) + \frac{1-\Theta_0}{2}e_3\in \B$.
\item[$\bullet$] We argue completely analogously when $\cZ_x=\{1\}$, and we use the directions
\begin{equation*}
\pm(\delta_1,-1,\eta_1,0),\ \ \pm(-\eta_1,\delta_1,1,0),\ \ \pm(\eta_1,\delta_1,-1,0),\ \ \pm(-\eta_1,0,\delta_1,1)
\end{equation*}
to illuminate $x$.
\item[$\bullet$] Assume finally that $\cZ_x=\{2\}$. If $|x_4| \leq \Theta_0$, and at the same time $\sign(x_1)=-\sign(x_3)$, then we use one of the directions $\pm(-\eta_1,0,\delta_1,1), \,\pm(-\eta_1,0,\delta_1,-1)$ to illuminate $x$.

\smallskip

If instead $|x_4|\leq \Theta_0 < 1$ and $\sign(x_1)=\sign(x_3)$, then one of the directions $\pm(1,0,1,0)$ illuminates $x$ because of Corollary \ref{cor:affine-set}.

\medskip

On the other hand, if $|x_4| > \Theta_0$, then $|x_1|\leq \Theta_0$. Then we choose $d\in \bigl\{\pm(-\eta_1,0,\delta_1,1),$ $\pm(-\eta_1,0,\delta_1,-1)\bigr\}$ so that $d_s\cdot x_s < 0$ for $s\in \{3,4\}$. We will have that $x+|x_4|d\in \intr\B$ (where we distinguish subcases in our analysis based on whether $|x_3|<\frac{1-\Theta_0}{4}$ or not; in the latter subcase we compare $x+|x_4|d$ with a convex combination of the form $(1-\lambda)\big|\vec{x}\big| + \lambda e_1$ 
where $\lambda < \Theta_0\delta_1 < |x_4|\delta_1$, and observe that, as long as $\eta_1 < \lambda(1-\Theta_0) < \Theta_0(1-\Theta_0)\delta_1$, the desired conclusion will follow).
\end{itemize}
\item[$\blacklozenge$] Finally we deal with the cases where $\abs{\cZ_x}=0$. We can argue as in the very last subcase when $|x_4| >\Theta_0$: indeed, then $\max(|x_1|,|x_2|)\leq \Theta_0$, and thus we can pick $d\in \bigl\{\pm(-\eta_1,0,\delta_1,1),\,\pm(-\eta_1,0,\delta_1,-1)\bigr\}$ so that $d_s\cdot x_s < 0$ for $s\in \{3,4\}$ to illuminate $x$ (considering again the displaced vector $x+|x_4|d$ and distinguishing subcases based on whether $|x_3|<\frac{1-\Theta_0}{4}$ or not; the only change we have to make is that, in the latter subcase, we compare $x+|x_4|d$ to a convex combination of the form $(1-\lambda)\big|\vec{x}\big| + \lambda (e_1+e_2)$).

\medskip

On the other hand, if $|x_4| \leq \Theta_0 < 1$, then one of the first 8 directions in ${\cal F}_{\ref{prop:R^4-exactly-one-pair},\ref{prop:R^4-exactly-two-pairs},\delta_1, \eta_1}$ will illuminate $x$ by Corollary \ref{cor:affine-set}.
\end{itemize}
Gathering all the restrictions on $\eta_1$, we see that, as long as 
\begin{equation*}
\eta_1 < \frac{1-\Theta_0}{2}\delta_1
\end{equation*}
(which also implies that $\eta_1 < \Theta_0(1-\Theta_0)\delta_1$ since $\Theta_0\geq 1/2$),
the set ${\cal F}_{\ref{prop:R^4-exactly-one-pair},\ref{prop:R^4-exactly-two-pairs},\delta_1, \eta_1}$ will illuminate $\B$ in this second main case that we considered.
\end{proof}


\begin{proposition}\label{prop:R^4-exactly-three-pairs}
Suppose that for a given $\B \in \mathcal{U}^4$ there are exactly three pairs of indices $i_1,i_2\in [4]$ such that $\norm{e_{i_1}+e_{i_2}}_\B = 1$ (and at the same time there are NO triples of indices $j_1,j_2,j_3\in [4]$ such that $e_{j_1}+e_{j_2}+e_{j_3}\in \B$). Then at least one of the following two statements holds:
\smallskip\\
(i) there exist $\delta_1 > 0$ and $\eta_1=\eta_{\delta_1}>0$ so that $\B$ can be illuminated by some coordinate permutation of the set
\begin{align*}
{\cal F}_{\ref{prop:R^4-exactly-three-pairs},\ref{prop:R^4-exactly-four-pairs},\delta_1, \eta_1}:= \bigl\{\!&\pm(1,\delta_1,\eta_1,0),\ \pm(\delta_1,-1,-\eta_1,0),\ \pm(\delta_1,-\eta_1,1,0),\ \pm(\delta_1,\eta_1,-1,0),
\\
&\quad \pm(0,\pm(\eta_1,\delta_1),1),\ \pm(0,1,-1,0)\bigr\};
\end{align*}
(ii) there exist $\delta_2 > 0$ and $\eta_2=\eta_{\delta_2}>0$ 
so that $\B$ can be illuminated by some coordinate permutation of the set
\begin{align*}
{\cal F}_{\ref{prop:R^4-exactly-three-pairs}, {\rm alt}, \delta_2,\eta_2} := \bigl\{\!&\pm(\eta_2,1,\delta_2,0),\ \pm(-\eta_2,1,\delta_2,0),\ \pm(-\eta_2,1,-\delta_2,0),
\\
&\quad \pm(\eta_2,\eta_2,1,\delta_2),\ \pm(-\eta_2,-\eta_2,1,\delta_2),\ \pm(\eta_2,\eta_2,-\delta_2,1),\ \pm(-\eta_2,-\eta_2,-\delta_2,1)\bigr\}.
\end{align*}
\end{proposition}
\begin{proof}
We single out three non-equivalent cases, and observe that any other case here can be reduced to one of these three after a coordinate permutation:
\begin{description}
\item[Case 1:] $\B$ contains the points $e_1+e_2, e_1+e_3$ and $e_2+e_3$ (but does not contain the point $e_1+e_2+e_3$).
\item[Case 2:] $\B$ contains the points $e_1+e_2, e_1+e_3$ and $e_3+e_4$.
\item[Case 3:] $\B$ contains the points $e_1+e_2, e_1+e_3$ and $e_1+e_4$.
\end{description}
We will see that if either Case 1 or Case 2 holds, then ${\cal F}_{\ref{prop:R^4-exactly-three-pairs},\ref{prop:R^4-exactly-four-pairs},\delta_1, \eta_1}$ illuminates $\B$ (for some suitably chosen $\delta_1,\eta_1$), while if Case 3 holds, we may use ${\cal F}_{\ref{prop:R^4-exactly-three-pairs}, {\rm alt}, \delta_2,\eta_2}$ to illuminate $\B$.

\medskip

\emph{Proof for Case 1.} For every $i\in [3]$, set $\theta_{i,4}=\norm{e_i+e_4}_\B^{-1}$. By our assumptions for this main case, $\Theta_0:=\max\{\theta_{i,4}:i\in [3]\}\in (0,1)$.

\smallskip

We pick $\delta_1 < \frac{1-\Theta_0}{4}$, and $\eta_1< \frac{\delta_1}{2}$ (which we will restrict further by the end of the proof).

\smallskip

Note that the only \underline{extreme} points $x$ of $\B$ with $\abs{\cZ_x}=3$ are $\pm e_4$. To illuminate such a point $x$, we use the directions $\mp(0,\eta_1,\delta_1,1)$: e.g. $e_4 + (0,-\eta_1,-\delta_1,-1) = (0,-\eta_1,-\delta_1,0)\in \intr\B$ if we compare with the point $e_2+e_3$.

We now consider the other possibilities for $\abs{\cZ_x}$.

\begin{itemize}
\item[$\blacklozenge$] Assume that $\abs{\cZ_x}=2$, and consider first the (potentially extreme) points $\pm e_i+\pm e_j$ with $i,j\in [3]$. The ``trickiest'' case here are the points $\pm e_2+ \pm e_3$, for which we can pick a direction $d$ from $\pm(\delta_1,-1,-\eta_1,0),\,\pm(\delta_1,-\eta_1,1,0)$ so that $d_s\cdot x_s < 0$ for $s\in \{2,3\}$. E.g. $e_2-e_3 + (\delta_1,-\eta_1,1,0) = (\delta_1,1-\eta_1,0,0)\in \intr\B$, and this can be confirmed if we compare with the point $e_1+e_2$.

Note also that there are no \underline{other} extreme points of $\B$ with $\abs{\cZ_x}=2$ and $4\in \cZ_x$. 

Consider now a (potentially extreme) point $x\in\partial \B$ of the form $x_{i_1} e_{i_1} + x_4 e_4$, where $i_1\in [3]$ (write also $\{i_2,i_3\} = [3]\setminus\{i_1\}$). By our assumptions, $\min(|x_{i_1}|, |x_4|)\leq \Theta_0$. If $|x_{i_1}|\leq \Theta_0$, then we choose $d\in \{\pm(0,\eta_1,\delta_1,1)\}$ so that $d_4\cdot x_4 < 0$. Then $x+|x_4|d\in \intr\B$, which can be seen if we compare with one of the points
\begin{equation*}
(|x_{i_1}|+|x_4|\delta_1)e_{i_1} + |x_4|\delta_1 e_{i_2} \quad \hbox{or}\quad |x_{i_1}|e_{i_1} + |x_4|\delta_1(e_{i_2}+e_{i_3})
\end{equation*}
(where the values of $i_2, i_3$ from $[3]\setminus\{i_1\}$ are suitably chosen based on $x$). Note that both the above points are interior points of $\B$ by our restriction on $\delta_1$, and at least one of them has larger in absolute value corresponding coordinates to those of $x+|x_4|d$.

If instead $|x_4|\leq \Theta_0$, then we pick a direction $d$ from the first 8 in ${\cal F}_{\ref{prop:R^4-exactly-three-pairs},\ref{prop:R^4-exactly-four-pairs},\delta_1, \eta_1}$ so that $d_{i_1}\cdot x_{i_1} < 0$. In a similar manner to above, we can compare the point $x+|x_{i_1}|d$ to the point $|x_4|e_4 + |x_{i_1}|\delta_1(e_{i_2}+e_{i_3})$ to see that the former point (as well as the latter) is in $\intr\B$.
\smallskip
\item[$\blacklozenge$] Next, assume that $\abs{\cZ_x}=1$. If $\cZ_x=\{4\}$, the first 8 directions of ${\cal F}_{\ref{prop:R^4-exactly-three-pairs},\ref{prop:R^4-exactly-four-pairs},\delta_1, \eta_1}$ illuminate $x$.
\begin{itemize}
\item[$\bullet$] If instead $\cZ_x=\{1\}$, and we assume first that $|x_4|\leq \Theta_0$, then one of the directions from $\big\{\!\pm\!(0,\pm(\eta_1,\delta_1),1),\,\pm(0,1,-1,0)\bigr\}$ illuminates $x$ (we use one of the first 4 if $\sign(x_2)=\sign(x_3)$, and we use one of the last 2 if $\sign(x_2)=-\sign(x_3)$ while relying on Corollary \ref{cor:affine-set} as well).

\smallskip

When $|x_4| > \Theta_0$, we will instead have that $|x_2|\leq \Theta_0$, and thus we can pick a direction $d^\prime$ from $\pm(0,\pm(\eta_1,\delta_1),1)$ so that $d^\prime_s\cdot x_s < 0$ for $s\in \{3,4\}$. Then $x+|x_4|d^\prime\in \intr\B$, since we can compare it to the point $(\Theta_0+\eta_1)e_2 + (1-\delta_1)e_3$ which is also an interior point of $\B$.
\item[$\bullet$] If $\cZ_x=\{2\}$, and we assume first that $|x_4|\leq \Theta_0$, then we use the directions $\pm(1,\delta_1,\eta_1,0),\,\pm(\delta_1,-\eta_1,1,0),\,\pm(\delta_1,\eta_1,-1,0)$ to illuminate $x$ (and we consider subcases based on whether $\min(|x_1|,|x_3|) < \frac{1-\Theta_0}{4}$ or not; in the latter subcase, we pick $d\in \bigl\{\pm(\delta_1,-\eta_1,1,0),\,\pm(\delta_1,\eta_1,-1,0)\bigr\}$ so that $d_s\cdot x_s < 0$ for $s\in\{1,3\}$, and we compare $x+\frac{1-\Theta_0}{4}d$ with a convex combination of the form 
\begin{equation*}
(1-\lambda)\big|\vec{x}\big| + \lambda \left(\frac{1-\Theta_0}{2}e_2 + \frac{1+\Theta_0}{2}e_4\right)
\end{equation*}
where $\lambda < \frac{1-\Theta_0}{4}\delta_1$; as long as $\eta_1 < 2\lambda < \frac{1-\Theta_0}{2}\delta_1$, we can conclude that $x+\frac{1-\Theta_0}{4}d\in \intr\B$).

\medskip

If instead $|x_4| > \Theta_0$, then we have that $\max(|x_1|,|x_3|) \leq \Theta_0$. We thus pick $d^\prime\in \bigl\{\pm(0,\pm(\eta_1,\delta_1),1)\bigr\}$ so that $d^\prime_s\cdot x_s < 0$ for $s\in \{3,4\}$, and we check that $x+|x_4|d^\prime\in \intr\B$ by comparing this point to a point of the form $\Theta_0(e_1+e_3) + \frac{1-\Theta_0}{2}e_2$, which is also an interior point of $\B$ (here we also use the fact that $\Theta_0\geq 1/2$, and thus $|(x+|x_4|d^\prime)_3| \leq \max\bigl(|x_3|-|x_4|\delta_1,\,|x_4|\delta_1\bigr)\leq \max(|x_3|,\,\delta_1)\leq \Theta_0$).
\item[$\bullet$] Finally, if $\cZ_x=\{3\}$, then we argue similarly to the previous subcase, and we illuminate $x$ using one of the directions
\begin{equation*}
\pm(1,\delta_1,\eta_1,0),\ \ \pm(\delta_1,-1,-\eta_1,0),\ \ \pm(0,\pm(\eta_1,\delta_1),1)
\end{equation*}
(while distinguishing subcases based on whether $|x_4| > \Theta_0$ or not, and in the latter case, based on whether $\min(|x_1|,|x_2|) < \frac{1-\Theta_0}{4}$ or not).
\end{itemize}
\item[$\blacklozenge$] It remains to deal with the cases where $\cZ_x = \emptyset$. If $|x_4|\leq \Theta_0 < 1$, then we use the first 8 directions in ${\cal F}_{\ref{prop:R^4-exactly-three-pairs},\ref{prop:R^4-exactly-four-pairs},\delta_1, \eta_1}$ to illuminate $x$.

If instead $|x_4| > \Theta_0$, then, as before, we observe that $\max(|x_1|,|x_2|)\leq \Theta_0$. Thus, we can illuminate $x$ using a direction $d$ from $\pm(0,\pm(\eta_1,\delta_1),1)$ which satisfies $d_s\cdot x_s < 0$ for $s\in \{3,4\}$ (to show that the point $x+|x_4|d\in \intr\B$, we distinguish subcases based on whether $|x_3|\geq \frac{1-\Theta_0}{4}$ or not; in those cases that $|x_3|$ is `not too small', we compare the point $x+|x_4|d$ with a convex combination of the form $(1-\lambda)\big|\vec{x}\big| + \lambda(e_1+e_2)$ where $\lambda < \Theta_0\delta_1 < |x_4|\delta_1$, and note that, as long as $\eta_1 < \lambda(1-\Theta_0) < (1-\Theta_0)\Theta_0\delta_1$, the desired conclusion holds).
\end{itemize}

We can conclude that, as long as $\delta_1 < \frac{1-\Theta_0}{4}$ and $\eta_1 < \frac{1-\Theta_0}{2}\delta_1$, the set ${\cal F}_{\ref{prop:R^4-exactly-three-pairs},\ref{prop:R^4-exactly-four-pairs},\delta_1, \eta_1}$ will illuminate the body $\B$ which contains the points $e_1+e_2, \,e_1+e_3$ and $e_2+e_3$, but not the points $e_i+e_4$, $i\in [3]$, or the point $e_1+e_2+e_3$.

\bigskip

\emph{Proof for Case 2.} Analogously to the previous main case, we set $\theta_{i,j}=\norm{e_i+e_j}_\B^{-1}$ for all $1\leq i < j\leq 4$ such that $(i,j) \notin\bigl\{(1,2),\,(1,3),\,(3,4)\bigr\}$, and then set 
\begin{equation*}
\Theta_0:=\max\big\{\theta_{i,j}: 1\leq i< j\leq 4,\,(i,j)\notin\{(1,2),\,(1,3),\,(3,4)\}\bigr\}.
\end{equation*}
We will again pick $\delta_1 < \frac{1-\Theta_0}{4}$ and $\eta_1 < \frac{1-\Theta_0}{2}\delta_1$.

Let $x$ be an \underline{extreme} boundary point of $\B$. In this main case, the additional assumption that $x$ is extreme implies that $\abs{\cZ_x}\leq 2$.
\begin{itemize}
\item[$\blacklozenge$] Assume that $\abs{\cZ_x} = 2$, and consider first the (potentially extreme) points $\pm e_i\pm e_j$ with $(i,j) \in\bigl\{(1,2),\,(1,3),\,(3,4)\bigr\}$. The ``trickiest'' cases here are the points $\pm (e_1+e_2)$ and $\pm e_3 \pm e_4$. We have e.g. that $-(e_1+e_2) + (1,\delta_1,\eta_1,0) = (0,-1+\delta_1,\eta_1,0)\in \intr\B$ because $(1-\delta_1) + \eta_1 < 1$. Similarly, for $\pm e_3 \pm e_4$ we use the directions $\pm(0,\pm(\eta_1,\delta_1),1)$, and we have e.g. that $-e_3 + e_4 +(0,\eta_1,\delta_1,-1) = (0,\eta_1,-1+\delta_1,0)\in \intr\B$ for the same reason as above.

\medskip

If $x$ is a different extreme point of $\B$ with $\abs{\cZ_x}=2$, then we must have $x=x_ie_i + x_je_j$ with $\{i,j\} = ([4]\setminus\cZ_x) \notin \bigl\{\{1,2\},\,\{1,3\},\,\{3,4\}\bigr\}$. But then $\min(|x_i|,|x_j|)\leq \Theta_0$. WLOG suppose that $|x_i|= \max(|x_i|,|x_j|)$, and pick a direction $d$ from the first 12 in ${\cal F}_{\ref{prop:R^4-exactly-three-pairs},\ref{prop:R^4-exactly-four-pairs},\delta_1, \eta_1}$ so that $d_i\cdot x_i < 0$. Then $x+|x_i|d\in \intr\B$, which can be seen in the same manner as in the previous main case, by comparing to coordinate permutations of the points
\begin{equation*}
(\Theta_0 + \delta_1)e_1 + \delta_1e_2 \quad \hbox{and}\quad \Theta_0e_1 + \delta_1(e_2+e_3)
\end{equation*}
(all coordinate permutations of these points are in $\B$, and are interior points because of our restriction on $\delta_1$).
\item[$\blacklozenge$] Next, assume that $\abs{\cZ_x} = 1$. We argue exactly as in the previous main case when $\cZ_x=\{4\}$ or when $\cZ_x=\{3\}$.

We also argue as in the previous main case when $\cZ_x=\{1\}$ or $=\{2\}$, and we additionally suppose that $|x_4| \leq \Theta_0$.
\begin{itemize}
\item[$\bullet$] If $\cZ_x=\{1\}$ and $|x_4| > \Theta_0$, we have that $|x_2| \leq \Theta_0$. As in the previous main case, we pick a direction $d^\prime$ from $\pm(0,\pm(\eta_1,\delta_1),1)$ so that $d^\prime_s\cdot x_s < 0$ for $s\in \{3,4\}$, but now we distinguish cases according to whether $|x_3|< \frac{1-\Theta_0}{4}$ or not. In both subcases, we consider the displaced vector $x+|x_4|d^\prime$. Moreover, when $|x_3|\geq \frac{1-\Theta_0}{4}$, we compare $x+|x_4|d^\prime$ to a convex combination of the form $(1-\lambda)\big|\vec{x}\big| + \lambda e_2$, where $\lambda < \Theta_0\delta_1$; similarly to before, as long as $\eta_1 < (1-\Theta_0)\lambda < (1-\Theta_0)\Theta_0\delta_1$, we obtain that $x+|x_4|d^\prime\in \intr\B$ in this subcase too.
\item[$\bullet$] If instead $\cZ_x=\{2\}$ and $|x_4| > \Theta_0$, we similarly have that $|x_1| \leq \Theta_0$ (but unlike the previous main case, we cannot claim anymore that $|x_3|\leq \Theta_0$). Still, as earlier, we pick a direction $d^\prime\in \bigl\{\pm(0,\pm(\eta_1,\delta_1),1)\bigr\}$ so that $d^\prime_s\cdot x_s < 0$ for $s\in \{3,4\}$, but now we distinguish cases based on whether $|x_3|\leq \Theta_0$ or not. In the former case, we continue as we did before, while, in the cases where $|x_3| > \Theta_0$, we compare $x+|x_4|d^\prime$ with a convex combination of the form
\begin{equation*}
(1-\lambda^\prime)\big|\vec{x}\big| + \lambda^\prime (e_1+e_2)
\end{equation*}
where $\lambda^\prime < \Theta_0\delta_1$, and obtain the desired conclusion as long as $\eta_1 < \lambda^\prime < \Theta_0\delta_1$.
\end{itemize}
\smallskip
\item[$\blacklozenge$] Finally, when $\cZ_x=\emptyset$, we argue exactly as in the previous main case.
\end{itemize}
We are done with the proof of Case 2 as well.

\bigskip

\emph{Proof for Case 3.} Just as in the previous main cases, for any pair $(i,j)$ for which $e_i+e_j\notin\B$, we set $\theta_{i,j}:=\norm{e_i+e_j}_\B^{-1}$, and then we define
\begin{equation*}
\Theta_0 : = \max\{\theta_{i,j}: 2\leq i < j\leq 4\}.
\end{equation*}
We pick $\delta_2 < \frac{1-\Theta_0}{4}$ and $\eta_2 < (1-\Theta_0)\delta_2$, and, under these restrictions, we show that ${\cal F}_{\ref{prop:R^4-exactly-three-pairs}, {\rm alt}, \delta_2,\eta_2}$ illuminates the convex body $\B$ of Case 3. We will be relying on the following

\medskip

{\bf Key Observation for Case 3.} Let $2\leq i < j\leq 4$, let $\mu_i, \mu_j\in (0,1)$ be such that $\mu_i+\mu_j < 1$, and let $\lambda_0\in (0,1)$. Then the point
\begin{equation*}
\lambda_0e_1 + \mu_i e_i + \mu_je_j
\end{equation*} 
is an interior point of $\B$, since it can be written as a (non-trivial) convex combination of the points $\lambda_0 e_1 + e_i,\,\lambda_0 e_1 + e_j\in \B$ and the interior point $\lambda_0 e_1$:
\begin{equation*}
\lambda_0e_1 + \mu_i e_i + \mu_je_j = \mu_i \bigl(\lambda_0 e_1 + e_i\bigr) \,+\,\mu_j\bigl(\lambda_0 e_1 + e_j\bigr)  + (1-\mu_i-\mu_j) \lambda_0 e_1.
\end{equation*}

\medskip

Consider now an extreme point $x\in \partial \B$. Observe that there are no such $x$ with $\abs{\cZ_x}=3$, thus we consider the remaining possibilities.
\begin{itemize}
\item[$\blacklozenge$] Assume that $\abs{\cZ_x}=2$, and consider first the (potentially extreme) points $\pm e_1 \pm e_j$, where $j\in \{2,3,4\}$. We pick a direction $d$ such that $m.c(d) = j$, and such that $d_s\cdot x_s < 0$ for $s\in \{1,j\}$. Then $x+|x_j|d = x+d$ satisfies:
\begin{itemize}
\item[$\cdot$] $(x+d)_j=0$,
\item[$\cdot$] $|(x+d)_1| = 1-\eta_2 < 1$,
\item[$\cdot$] and $|(x+d)_{i_1}| + |(x+d)_{i_2}| \leq \delta_2 + \eta_2 < 2\delta_2 < 1$, where $\{i_1,i_2\} = [4]\setminus\{1,j\}$.
\end{itemize}
We obtain that $x+d\in \intr\B$ from the above ``key observation''.

\smallskip

Next, note that there are no other extreme points in $\B$ with $\abs{\cZ_x}=2$ and $1\notin \cZ_x$. Thus, assume now that $x= x_i e_i + x_j e_j$ with $i,j\in \{2,3,4\}$, $i\neq j$. WLOG assume that $\max(|x_i|,|x_j|) = |x_i|$, which in turn implies that $|x_j| \leq \Theta_0$. Thus, we pick a direction $d^\prime$ such that $m.c.(d^\prime) = i$ and $d^\prime_i\cdot x_i < 0$. We will have that $x+|x_i|d^\prime$ satisfies:
\begin{itemize}
\item[$\cdot$] $(x+|x_i|d^\prime)_i = 0$,
\item[$\cdot$] $|(x+|x_i|d^\prime)_1| = |x_i|\eta_2\leq \eta_2$,
\item[$\cdot$] $|(x+|x_i|d^\prime)_j| \leq \Theta_0 + \delta_2$,
\item[$\cdot$] and $|(x+|x_i|d^\prime)_s| \leq \delta_2$, where $\{s\} = \{2,3,4\}\setminus\{i,j\}$.
\end{itemize}
Since $|(x+|x_i|d^\prime)_j| + |(x+|x_i|d^\prime)_s| \leq \Theta_0 + 2\delta_2 < 1$, we conclude from the ``key observation'' that $x+|x_i|d^\prime\in \intr\B$.
\item[$\blacklozenge$] Next, assume that $\abs{\cZ_x}=0$. Suppose first that $|x_4|\leq \Theta_0$. Except for the cases where $\sign(x_1) = \sign(x_2) = - \sign(x_3)$, we can use the first 6 directions of ${\cal F}_{\ref{prop:R^4-exactly-three-pairs}, {\rm alt}, \delta_2,\eta_2}$ to illuminate $x$.

On the other hand, if $\sign(x_1) = \sign(x_2) = - \sign(x_3)$, then one of the directions $\pm(-\eta_2,-\eta_2,1,\delta_2),\,\pm(\eta_2,\eta_2,-\delta_2,1)$ illuminates $x$ (based also on what $\sign(x_4)$ is).

\medskip

Next, suppose that $|x_4| > \Theta_0$. In that case $|x_2| \leq \Theta_0$. If we first assume that $|x_3| < \frac{1-\Theta_0}{4}$, then we pick a direction $d$ from $\pm(\eta_2,\eta_2,-\delta_2,1),\,\pm(-\eta_2,-\eta_2,-\delta_2,1)$ so that $d_s\cdot x_s < 0$ for $s\in \{1,4\}$. We will have that $x+|x_4|d$ satisfies:
\begin{itemize}
\item[$\cdot$] $(x+|x_4|d)_4 =0$,
\item[$\cdot$] $|(x+|x_4|d)_1| \leq 1-|x_4|\eta_2 < 1-\Theta_0\eta_2$,
\item[$\cdot$] $|(x+|x_4|d)_2| \leq \Theta_0 + \eta_2$,
\item[$\cdot$] and $|(x+|x_4|d)_3| \leq \frac{1-\Theta_0}{4} + \delta_2$.
\end{itemize}
Thus $|(x+|x_4|d)_2| + |(x+|x_4|d)_3| \leq \Theta_0 + \frac{1-\Theta_0}{4} + 2\delta_2 < \Theta_0 + \frac{3}{4}(1-\Theta_0) < 1$, which implies that $x+|x_4|d\in \intr\B$ because of the ``key observation''.

\smallskip

If instead $|x_3|\geq \frac{1-\Theta_0}{4}$, then we pick a direction $d^\prime$ from the last 8 in ${\cal F}_{\ref{prop:R^4-exactly-three-pairs}, {\rm alt}, \delta_2,\eta_2}$ so that $d^\prime_s\cdot x_s < 0$ for $s\in \{1,3,4\}$. 
We can compare $x+\frac{1-\Theta_0}{4}d^\prime$ with a convex combination of the form
\begin{equation*}
(1-\lambda)\big|\vec{x}\big| + \lambda (e_1+e_2)
\end{equation*}
where $\lambda < \frac{1-\Theta_0}{4}\delta_2$. As long as $\eta_2 < (1-\Theta_0)\delta_2$, we obtain that $x+\frac{1-\Theta_0}{4}d^\prime\in \intr\B$.
\item[$\blacklozenge$] Finally, we suppose that $\abs{\cZ_x}=1$. Assume first that $\cZ_x=\{r\}$ with $r\in \{2,3,4\}$. Let us write $\{i,j\} = \{2,3,4\}\setminus\cZ_x= \{2,3,4\}\setminus\{r\}$, and WLOG let us assume that $\max(|x_i|,|x_j|) = |x_i|$. Then we will also have that $|x_j| = \min(|x_i|,|x_j|) \leq \Theta_0$.

We pick a direction $d$ so that $m.c.(d) = i$, and so that $d_s\cdot x_s < 0$ for $s\in \{1,i\}$. Then $x+|x_i|d$ will satisfy:
\begin{itemize}
\item[$\cdot$] $(x+|x_i|d)_i = 0$,
\item[$\cdot$] $|(x+|x_i|d)_1| \leq 1-|x_i|\eta_2 < 1$,
\item[$\cdot$] and $|(x+|x_i|d)_j| + |(x+|x_i|d)_r| \leq \Theta_0 + |x_i|\delta_2 + |x_i|\delta_2 < 1$, given our restrictions on $\delta_2$ and $\eta_2$.
\end{itemize}
Thus, by the ``key observation'', $x+|x_i|d\in \intr\B$.

\medskip

It remains to deal with the cases where $\cZ_x = \{1\}$. Here, we first consider the subcases where $|x_3| < \frac{1-\Theta_0}{4}$. Let $\{i,j\} = \{2,4\}$, and let us write $i$ for the index where $\max(|x_2|,|x_4|)$ is attained (if $|x_2| = |x_4|$, then set $i=2$). We then know that $|x_j| \leq \Theta_0$. Again, we pick a direction $d^\prime$ so that $m.c.(d^\prime) = i$ and $d^\prime_i\cdot x_i < 0$. Then, similarly to above, we can check that $x+|x_i|d^\prime\in \intr\B$ by the ``key observation'' (given that $|(x+|x_i|d^\prime)_1|\leq \eta_2$, and $|(x+|x_i|d^\prime)_j| + |(x+|x_i|d^\prime)_3| \leq \Theta_0 + |x_i|\delta_2 + \frac{1-\Theta_0}{4} + |x_i|\delta_2 < 1$).

\smallskip

We argue very similarly when $|x_3| \geq \frac{1-\Theta_0}{4}$ while at the same time $|x_i|= \max(|x_i|,|x_j|)= \max(|x_2|,|x_4|)< \frac{1-\Theta_0}{4}$: in those subcases, we pick a direction $d^\prime$ from $\pm(\eta_2,\eta_2,1,\delta_2)$ so that $d^\prime_3\cdot x_3 < 0$, and check in an analogous way that $x+|x_3|d^\prime\in \intr\B$.

\medskip

The last subcase to consider is when $\min(|x_i|,|x_3|) \geq \frac{1-\Theta_0}{4}$. Then we pick a direction $d^{\prime\prime}$ so that $\bigl\{|d^{\prime\prime}_i|,\,|d^{\prime\prime}_3|\bigr\} = \{1,\delta_2\}$, and so that $d^{\prime\prime}_s\cdot x_s < 0$ for $s\in \{i,3\}$. We compare the displaced vector $x+\frac{1-\Theta_0}{4}d^{\prime\prime}$ with a convex combination of the form
\begin{equation*}
(1-\lambda^\prime)\big|\vec{x}\big| + \lambda^\prime (e_1+e_j)
\end{equation*}
where $\lambda^\prime < \frac{1-\Theta_0}{4}\delta_2$: the former vector is guaranteed to be in $\intr\B$ as long as $\frac{1-\Theta_0}{4}\eta_2 < \lambda^\prime(1-\Theta_0) \leq \lambda^\prime (1-|x_j|)\, \Leftrightarrow\,\eta_2 < 4\lambda^\prime < (1-\Theta_0)\delta_2$.
\end{itemize}

\noindent This completes the proof in all main cases.
\end{proof}

\begin{remark}\label{rem:alt-set-Case12-3squares}
Taking into account that parameters which appear as later subscripts depend on previous parameters, and can be chosen much smaller if needed, we can now also verify, through a minor adaptation of the above argument, that, for Cases 1 and 2 of Proposition \ref{prop:R^4-exactly-three-pairs}, we can use the illuminating set 
\begin{align*}
{\cal F}_{\ref{prop:R^4-no-pairs},\ref{prop:R^4-exactly-five-pairs},\delta, \eta,\zeta}:= \bigl\{\!&\pm(1,\delta,\eta,0),\ \pm(\delta,-1,-\eta,0),\ \pm(\delta,-\eta,1,\zeta),\ \pm(-\delta,-\eta,1,\zeta),
\\
&\quad \pm(0,\pm(\eta,\delta),1),\ \pm(0,1,-\delta,\eta)\bigr\}
\end{align*}
instead of
\begin{align*}
{\cal F}_{\ref{prop:R^4-exactly-three-pairs}, \ref{prop:R^4-exactly-four-pairs},\delta_1, \eta_1}:= \bigl\{\!&\pm(1,\delta_1,\eta_1,0),\ \pm(\delta_1,-1,-\eta_1,0),\ \pm(\delta_1,-\eta_1,1,0),\ \pm(\delta_1,\eta_1,-1,0),
\\
&\quad \pm(0,\pm(\eta_1,\delta_1),1),\ \pm(0,1,-1,0)\bigr\}.
\end{align*}
This reduces further the number of non-equivalent illuminating sets that we need (we still chose to work with the latter set to keep the proof a little more transparent).
\end{remark}

Next we prove the case where $\B$ contains \underline{five} $2$-dimensional unit subcubes, because it is much more similar to the previous settings compared to the case where $\B$ contains \underline{four} such subcubes (which we will handle last).


\begin{proposition}\label{prop:R^4-exactly-five-pairs}
Suppose that for a given $\B \in \mathcal{U}^4$ there are exactly five pairs of indices $i_1,i_2\in [4]$ such that $\norm{e_{i_1}+e_{i_2}}_\B=1$ (and at the same time there are NO triples of indices $j_1,j_2,j_3\in [4]$ such that $e_{j_1}+e_{j_2}+e_{j_3}\in \B$). Then there exist $\delta>0$, $\eta = \eta_\delta>0$ and $\zeta=\zeta_{\delta,\eta}>0$ such that $\B$ can be illuminated by a coordinate permutation of the set
\begin{align*}
{\cal F}_{\ref{prop:R^4-no-pairs},\ref{prop:R^4-exactly-five-pairs},\delta, \eta,\zeta}= \bigl\{\!&\pm(1,\delta,\eta,0),\ \pm(\delta,-1,-\eta,0),\ \pm(\delta,-\eta,1,\zeta),\ \pm(-\delta,-\eta,1,\zeta),
\\
&\quad \pm(0,\pm(\eta,\delta),1),\ \pm(0,1,-\delta,\eta)\bigr\}.
\end{align*}
\end{proposition}
\begin{proof}
WLOG we can assume that the only pair of indices $i_1\neq i_2\in [4]$ for which $e_{i_1}+e_{i_2}\notin \B$ satisfies $\{i_1,i_2\}=\{1,4\}$. Set $\Theta_0=\|e_1+e_4\|_\B^{-1} < 1$. Also, for each $j\in [4]$ set 
\[\gamma_j=\|{\bm 1}-e_j\|_\B^{-1}.\]
By our assumptions, $\gamma_0:=\max_{j\in [4]}\gamma_j < 1$. Fix now some
\begin{equation*}
\delta < \frac{\min(1-\Theta_0,\,1-\gamma_0)}{4}
\end{equation*}
and some $\zeta < \eta/2 < \delta/4$ (which we will further restrict shortly).

\smallskip

Clearly there are no \underline{extreme} boundary points $x\in \B$ with $\abs{\cZ_x}=3$, thus we focus on the remaining cases.
\begin{itemize}
\item[$\blacklozenge$] Assume that $\abs{\cZ_x}=2$. Here most cases are similar, except for the case where $\cZ_x=\{2,3\}$ (or in other words, where $x=(x_1,0,0,x_4)$). From the remaining cases the only potentially extreme points are of the form $\pm e_i + \pm e_j$ where $\{i,j\}\neq \{1,4\}$. We illuminate these points using a direction $d$ as follows:
\begin{center}
\begin{tabular}{|c|c|}
\hline
Boundary point & Possible illuminating directions
\\
\hline\hline
$\phantom{\Bigl(}\pm e_1+\pm e_2\phantom{\Bigr)}$ &  $\pm(1,\delta,\eta,0),\ \pm(\delta,-1,-\eta,0)$
\\
\hline
$\phantom{\Bigl(}\pm e_1+\pm e_3\phantom{\Bigr)}$ & $\pm(\delta,-\eta,1,\zeta),\ \pm(-\delta,-\eta,1,\zeta)$
\\
\hline
$\phantom{\Bigl(}\pm e_2+\pm e_3\phantom{\Bigr)}$ & $\pm(\delta,-1,-\eta,0),\ \pm(\delta,-\eta,1,\zeta)$
\\
\hline
$\phantom{\Bigl(}\pm e_2+\pm e_4,\ \pm e_3+\pm e_4\phantom{\Bigr)}$ & $\pm(0,\eta,\delta,1), \ \pm(0,-\eta,-\delta,1)$
\\
\hline
\end{tabular}
\end{center} 
so that $d_s\cdot x_s< 0$ for $s\in \{i,j\}$. E.g.
\begin{equation*}
(e_2+e_3) + (\delta,-1,-\eta,0) = (\delta,\,0,\,1-\eta,\,0)\in \intr\B,
\end{equation*}
which follows simply from the facts that $e_1+e_3\in \B$ and that $\delta,\eta\in (0,1)$. Similarly
\begin{equation*}
(-e_2+e_3) + (-\delta,\eta,-1,-\zeta) = \bigl(-\delta,\,-(1-\eta),\,0,\,-\zeta\bigr)\in \intr\B,
\end{equation*}
which can be seen by comparing to the point $(1-\eta)(e_1+e_2) + \zeta e_4$, that is also an interior point of $\B$ since $1 -\eta + \zeta < 1$.

\medskip

Now assume that $x=(x_1,0,0,x_4)$. Then $\min(|x_1|,|x_4|)\leq \Theta_0$. Thus, if $i$ is the (smallest) index at which $\max(|x_1|,|x_4|)$ is attained, we can illuminate $x$ choosing a direction $d$ from $\pm(1,\delta,\eta,0),\ \pm(0,\eta,\delta,1)$ so that $m.c.(d) = i$ and $d_i\cdot x_i < 0$. We will have that $x+|x_i|d\in \intr\B$, which can be seen by comparing with the points $\frac{1+\Theta_0}{2}e_j + \frac{1-\Theta_0}{2}(e_2+e_3)\in \B$, $j\in \{1,4\}$. 
\item[$\blacklozenge$] Next, assume that $\abs{\cZ_x}=1$. 
\begin{itemize}
\item[$\bullet$] If $\cZ_x=\{4\}$, then we illuminate $x$ choosing from the directions
\begin{equation*}
\pm(1,\delta,\eta,0),\ \pm(\delta,-1,-\eta,0),\ \pm(\delta,-\eta,1,\zeta),\ \pm(-\delta,-\eta,1,\zeta).
\end{equation*}
This is straightforward to do in the cases that $\sign(x_2)=\sign(x_3)$, so we examine how to handle the remaining subcases here. 

Note that $\min(|x_1|,|x_2|,|x_3|) \leq \gamma_4 \leq \gamma_0$. If $\sign(x_2)=-\sign(x_3)$, but also $|x_3| \leq \gamma_0$, then we still pick a direction $d$ from the first 4 above so that $d_s\cdot x_s < 0$ for $s\in [2]$. Then, if $i=m.c.(d)\in [2]$, we will have that $x+|x_i|d \in \intr\B$ by comparing it to one of the points $e_1+e_3$ or $e_2+e_3$.

\medskip

If instead $|x_3|>\gamma_0$, then we pick a direction $d^\prime$ from $\pm(\delta,-\eta,1,\zeta),\,\pm(-\delta,-\eta,1,\zeta)$ so that $d^\prime_s\cdot x_s < 0$ for $s\in [3]$. We then check that $x+|x_3|d^\prime$ satisfies:
\begin{itemize}
\item[$\cdot$] $(x+|x_3|d^\prime)_3 = 0$,
\item[$\cdot$] $|(x+|x_3|d^\prime)_4| = |x_3|\zeta$,
\item[$\cdot$] and $|(x+|x_3|d^\prime)_1| \leq 1-|x_3|\delta\leq 1-|x_3|\eta$, and similarly $|(x+|x_3|d^\prime)_2| \leq 1-|x_3|\eta$.
\end{itemize}
Thus we can compare $x+|x_3|d^\prime$ with the point $(1-|x_3|\eta)(e_1+e_2) + |x_3|\zeta \,e_4$, with the latter point being an interior point of $\B$, since $1-|x_3|\eta + |x_3|\zeta < 1$.
\item[$\bullet$] If $\cZ_x=\{1\}$, we use the directions $\pm(0,\pm(\eta,\delta),1),\,\pm(0,1,-\delta,\eta)$ to illuminate $x$. Again, this will be straightforward when $\sign(x_2)=\sign(x_3)$, so we examine the remaining subcases.

Note that $\min(|x_2|,|x_3|,|x_4|) \leq \gamma_1 \leq \gamma_0$. If $|x_4|\leq \gamma_0$ and $\sign(x_2)=-\sign(x_3)$, then we pick the unique direction $d\in \{\pm(0,1,-\delta,\eta)\}$ which satisfies $d_s\cdot x_s < 0$ for $s\in \{2,3\}$. We will have that $x+|x_2|d\in \intr\B$, which can be seen by comparing to the point $e_3+e_4$.

If instead $|x_4| > \gamma_0$, then $\min(|x_2|,|x_3|) \leq \gamma_0$. Let $i\in \{2,3\}$ be the index at which $\max(|x_2|,|x_3|)$ is attained, and pick $d^\prime\in \{\pm(0,\pm(\eta,\delta),1)\}$ so that $d^\prime_s\cdot x_s < 0$ for $s\in \{i,4\}$. Then $x+ |x_4|d^\prime\in \intr\B$, which can be seen by comparing to the point $e_2+e_3$.
\item[$\bullet$] Now, assume that $\cZ_x = \{2\}$. Then $\min(|x_1|,|x_4|)\leq \Theta_0$. If we also have that $|x_3|\leq \frac{1-\Theta_0}{4}$, and if $i$ is the (smallest) index at which $\max(|x_1|,|x_4|)$, then we pick a direction $d$ from $\pm(1,\delta,\eta,0),\,\pm(0,\eta,\delta,1)$ so that $m.c.(d)=i$ and $d_i\cdot x_i < 0$. We have that $x+|x_i|d\in \intr\B$, which can be seen by comparing to the point $\frac{1+\Theta_0}{2}e_j + \frac{1-\Theta_0}{2} (e_2+e_3)$, $j\in \{1,4\}\setminus\{i\}$.

Suppose now that $|x_3| > \frac{1-\Theta_0}{4}$. Then we pick a direction $d^\prime$ from
\begin{equation*}
\pm(\delta,-\eta,1,\zeta),\ \ \pm(-\delta,-\eta,1,\zeta),\ \ \pm(0,\pm(\eta,\delta),1)
\end{equation*}
so that $\{|d^\prime_i|,\,|d^\prime_3|\} = \{1,\delta\}$ and so that $d^\prime_s\cdot x_s < 0$ for $s\in \{i,3\}$. We compare $x+\frac{1-\Theta_0}{4}d^\prime$ with a convex combination of the form $(1-\lambda)\big|\vec{x}\big| + \lambda (e_2+e_j)$, where $j\in \{1,4\}\setminus\{i\}$ and $\lambda < \frac{1-\Theta_0}{4}\delta$. As long as $\eta < \delta$ and $\zeta < 4\lambda < (1-\Theta_0)\delta$, we can conclude that $x+\frac{1-\Theta_0}{4}d^\prime\in \intr\B$.
\item[$\bullet$] Analogously we argue if $\cZ_x = \{3\}$, while picking a direction $d$ from
\begin{equation*}
\pm(1,\delta,\eta,0),\ \ \pm(\delta,-1,-\eta,0),\ \ \pm(0,\pm(\eta,\delta),1)
\end{equation*}
to illuminate $x$. For most subcases we can simply rely on the restrictions $\eta < \delta < \frac{1-\Theta_0}{4}$. 

\smallskip

In the subcases where $|x_2| > \frac{1-\Theta_0}{4}$ and $|x_4| > |x_1|$, we pick $d\in \bigl\{\pm(0,\pm(\eta,\delta),1)\bigr\}$ so that $d_s\cdot x_s < 0$ for $s\in \{2,4\}$. We will have that $x+|x_4|d\in \intr\B$ because we can compare this displaced vector to the vector
\begin{equation*}
(1-|x_4|\eta)e_2 \ +\ \frac{1+\Theta_0}{2}e_1\  +\  \frac{1-\Theta_0}{4}e_3
\end{equation*}
which we can show is an interior point of $\B$ as well, in a similar manner to how we proved the ``Key Observation for Case 3'' of Proposition \ref{prop:R^4-exactly-three-pairs}.
\end{itemize}
\item[$\blacklozenge$] Finally, we consider the cases where $\cZ_x = \emptyset$. If $|x_1| \leq \Theta_0$, then we can use the directions $\pm(0,\pm(\eta,\delta),1),\,\pm(0,1,-\delta,\eta)$ to illuminate $x$, except in the subcases where $\sign(x_2)=-\sign(x_3)=-\sign(x_4)$. In these latter subcases, we can instead use one of the directions $\pm(\delta,-\eta,1,\zeta),\,\pm(-\delta,-\eta,1,\zeta)$ (which we choose based also on what $\sign(x_1)$ is).

\medskip

If instead $|x_1| > \Theta_0$, then $|x_4| \leq \Theta_0 < 1$. Here we also consider subcases according to whether $\max(|x_2|,|x_3|) < \frac{1-\Theta_0}{4}$ or not. If this maximum is ``small'', we simply pick $d\in \{\pm(1,\delta,\eta,0)\}$ so that $d_1\cdot x_1 < 0$, and we compare the displaced vector $x+|x_1|d$ to the vector
\begin{equation*}
\frac{1-\Theta_0}{2}(e_2+e_3) + \frac{1+\Theta_0}{2}e_4\in \B.
\end{equation*} 
If instead $\max(|x_2|,|x_3|) \geq \frac{1-\Theta_0}{4}$, then we pick a direction $d^\prime$ from among all the first 8 in ${\cal F}_{\ref{prop:R^4-no-pairs},\ref{prop:R^4-exactly-five-pairs},\delta, \eta,\zeta}$. We consider further subcases according to whether it also holds that $\min(|x_2|,|x_3|)\geq \frac{1-\Theta_0}{4}$ or not: if the minimum is ``small'', and $i\in \{2,3\}$ is the index at which $\max(|x_2|,|x_3|)$ is attained, then we pick $d^\prime$ so that $\{|d^\prime_1|,\,|d^\prime_i|\} = \{1,\delta\}$ and so that $d^\prime_s\cdot x_s < 0$ for $s\in \{1,i\}$; on the other hand, if the minimum is `not too small', we pick $d^\prime$ so that $d^\prime_s\cdot x_s < 0$ for all $s\in [3]$. Then, in all subcases we can conclude that, as long as
\begin{equation*}
\zeta < \eta < (1-\Theta_0)\delta, \qquad \hbox{and also}\ \ \zeta < (1-\Theta_0)\eta,
\end{equation*}
the displaced vector $x+\frac{1-\Theta_0}{4}d^\prime$ (with $x$ being displaced in the appropriately chosen direction $d^\prime$, as explained above) will be an interior point of $\B$.
\end{itemize}  
The proof is complete.
\end{proof}

\bigskip


\begin{proposition}\label{prop:R^4-exactly-four-pairs}
Suppose that for a given $\B \in \mathcal{U}^4$, {\bf which is not an affine image of the cube}, there are exactly four pairs of indices $i_1,i_2\in [4]$ such that $e_{i_1}+e_{i_2}\in \B$, and at the same time there are NO triples of indices $j_1,j_2,j_3\in [4]$ such that $e_{j_1}+e_{j_2}+e_{j_3}\in \B$. 
Then at least one of the following two statements holds:
\smallskip\\
(i) there exist $\delta_1 > 0$ and $\eta_1=\eta_{\delta_1}>0$ so that $\B$ can be illuminated by some coordinate permutation of the set
\begin{align*}
{\cal F}_{\ref{prop:R^4-exactly-three-pairs},\ref{prop:R^4-exactly-four-pairs},\delta_1, \eta_1} = \bigl\{\!&\pm(1,\delta_1,\eta_1,0),\ \pm(\delta_1,-1,-\eta_1,0),\ \pm(\delta_1,-\eta_1,1,0),\ \pm(\delta_1,\eta_1,-1,0),
\\
&\quad \pm(0,\pm(\eta_1,\delta_1),1),\ \pm(0,1,-1,0)\bigr\};
\end{align*}
(ii) there exist $\delta_2 > 0$ and $\eta_2=\eta_{\delta_2}>0$ so that $\B$ can be illuminated by some coordinate permutation of the set
\begin{align*}
{\cal F}_{\ref{prop:R^4-exactly-four-pairs}, {\rm alt}, \delta_2,\eta_2} := \bigl\{\!&\pm(1,-\eta_2,-\delta_2,-\delta_2),\ \pm(-\eta_2,1,-\delta_2,-\delta_2),
\\&\pm(\delta_2,0,1,-\eta_2),\ \pm(\delta_2,0,-\eta_2,1),\  \pm(0,\delta_2,1,-\eta_2),\ \pm(0,\delta_2,-\eta_2,1)\bigr\}.
\end{align*}
\end{proposition}
\begin{proof}
Up to coordinate permutations, there are two main cases to consider:
\begin{itemize}
\item[Case 1:] $\B$ contains the points $e_1+e_2$, $e_1+e_3$, $e_2+e_3$ and $e_3+e_4$ (and does not contain the point $e_1+e_2+e_3$).
\item[Case 2:] $\B$ contains the points $e_1+e_3$, $e_1+e_4$ and $e_2+e_3$, $e_2+e_4$. Here we need to further observe that the convex hull of all coordinate reflections of these points is $CP_1^2\times CP_1^2$ which is an affine image of the 4-dimensional cube, therefore by our assumptions $\B$ must contain at least one more point $z_0$ which satisfies $|z_{0,1}|+|z_{0,2}|>1$ and/or $|z_{0,3}|+|z_{0,4}|>1$. We can check that this is equivalent to having $\beta_0:=\max\bigl\{\|e_1+e_2\|_\B^{-1}, \|e_3+e_4\|_\B^{-1}\bigr\}>\frac{1}{2}$. WLOG we will assume here that $\|e_3+e_4\|_\B^{-1}=\beta_0>\frac{1}{2}$.
\end{itemize}

\medskip

\emph{Proof for Case 1.} 
We set $\Theta_0 := \max\bigl\{\norm{e_1+e_4}_\B^{-1},\,\norm{e_2+e_4}_\B^{-1}\bigr\}$, and note that $\Theta_0 \in (0,1)$.
We pick $\delta_1< \frac{1-\Theta_0}{4}$ and $\eta_1 < \frac{1-\Theta_0}{2}\delta_1$, and we will show that ${\cal F}_{\ref{prop:R^4-exactly-three-pairs},\ref{prop:R^4-exactly-four-pairs},\delta_1, \eta_1}$ illuminates the convex body $\B$ that we consider here, which satisfies the assumptions of the 1st main case.

\smallskip 

Let $x$ be a boundary point of $\B$. As previously, we ignore boundary points which are guaranteed to not be extreme, so we do not consider cases where $\abs{\cZ_x}=3$.
\begin{itemize}
\item[$\blacklozenge$] Assume that $\abs{\cZ_x}=2$. We first consider points of the form $\pm e_i+ \pm e_j$, where $\{i,j\}\notin\bigl\{\{1,4\},\,\{2,4\}\bigr\}$. Given such a point, we pick a direction $d$ as follows:
\begin{center}
\begin{tabular}{|c|c|}
\hline
Boundary point & Possible illuminating directions
\\
\hline\hline
$\phantom{\Bigl(}\pm e_1+\pm e_2\phantom{\Bigr)}$ &  $\pm(1,\delta_1,\eta_1,0),\ \pm(\delta_1,-1,-\eta_1,0)$
\\
\hline
$\phantom{\Bigl(}\pm e_1+\pm e_3\phantom{\Bigr)}$ & $\pm(\delta_1,-\eta_1,1,0),\ \pm(\delta_1,\eta_1,-1,0)$
\\
\hline
$\phantom{\Bigl(}\pm e_2+\pm e_3\phantom{\Bigr)}$ & $\pm(\delta_1,-1,-\eta_1,0),\ \pm(\delta_1,-\eta_1,1,0)$
\\
\hline
$\phantom{\Bigl(}\pm e_3+\pm e_4\phantom{\Bigr)}$ & $\pm(0,\eta_1,\delta_1,1), \ \pm(0,-\eta_1,-\delta_1,1)$
\\
\hline
\end{tabular}
\end{center} 
so that $d_s\cdot x_s< 0$ for $s\in \{i,j\}$. We will have that $x+d\in \intr\B$ because it can be compared to one of these points again: if $1\leq i < j \leq 3$, then $y_x:= (x_i e_i + x_j e_j) + d$ (where $|x_i|=|x_j| = 1$) will satisfy $\norm{y_x}_\infty < 1$, $\abs{\cZ_{y_x}}=2$, and $4\in \cZ_{y_x}$, thus we can see that $y_x\in\intr\B$ by comparing it to one of the points $e_1+e_2,\,e_1+e_3,\,e_2+e_3$.

Similarly, if e.g. $x=e_3-e_4$, then $x+(0,-\eta_1,-\delta_1,1)= (0,-\eta_1,1-\delta_1,0)\in \intr\B$, which can be seen by comparing with the point $e_2+e_3$.

\medskip

Note now that no other point of $\B$ with support the same as one of the above points can be extreme, as they will be contained in the convex hull of the above points, so all these other points can be illuminated by the same directions. This leaves two more subcases to consider here.
\begin{itemize}
\item[$\bullet$] Suppose that $x= x_1e_1 + x_4e_4$. Then $\min(|x_1|,|x_4|)\leq \Theta_0$. If $|x_4|\leq \Theta_0$, then we illuminate $x$ using the unique direction $d\in \{\pm(1,\delta_1,\eta_1,0)\}$ satisfying $d_1\cdot x_1 < 0$. We will have that $x+|x_1|d\in \intr\B$, by comparing it to the point $\frac{1+\Theta_0}{2}e_4 + \frac{1-\Theta_0}{2}(e_2+e_3)$.

\smallskip

Analogously, if $|x_4| >\Theta_0$, we use the unique direction $d^\prime\in\{\pm(0,\eta_1,\delta_1,1)\}$ which satisfies $d^\prime_4 \cdot x_4 < 0$: we will have that $x+|x_4|d^\prime\in \intr\B$, as before.
\medskip
\item[$\bullet$] Finally, suppose that $x=x_2e_2 + x_4e_4$. In this subcase, pick the unique direction $d\in \bigl\{\pm(0,\pm(\eta_1,\delta_1),1)\bigr\}$ which satisfies $d_s\cdot x_s < 0$ for $s\in \{2,4\}$. Then $x+|x_4|d\in \intr\B$, since $(x+|x_4|d)_4 = (x+|x_4|d)_1 = 0$, while $|(x+|x_4|d)_2|\leq 1-|x_4|\eta_1 < 1$, and $|(x+|x_4|d)_3|= |x_4|\delta_1<1$.
\end{itemize}
\medskip
\item[$\blacklozenge$] Now assume that $\abs{\cZ_x}=1$. If $\cZ_x = \{4\}$, then the first 8 directions of ${\cal F}_{\ref{prop:R^4-exactly-three-pairs},\ref{prop:R^4-exactly-four-pairs},\delta_1, \eta_1}$ illuminate $x$. We now examine the remaining subcases here.
\begin{itemize}
\item[$\bullet$] If $\cZ_x=\{1\}$, then we illuminate $x$ using the directions $\pm(0,\pm(\eta_1,\delta_1),1),\ \pm(0,1,-1,0)$. Indeed, if $|x_2| > \Theta_0$, then necessarily $|x_4|\leq \Theta_0$, and thus we can use the first 4 directions here if $\sign(x_2)=\sign(x_3)$, otherwise we can rely on Corollary \ref{cor:affine-set} and illuminate $x$ using one of the directions $\pm(0,1,-1,0)$.

\smallskip

If instead $|x_2|\leq \Theta_0$, then we pick $d$ from $\pm(0,\pm(\eta_1,\delta_1),1)$ so that $d_s\cdot x_s < 0$ for $s\in \{3,4\}$. We will have that $x+|x_4|d\in \intr\B$, since $(x+|x_4|d)_4 = (x+|x_4|d)_1 = 0$, while $|(x+|x_4|d)_3|\leq 1-|x_4|\delta_1 < 1$, and $|(x+|x_4|d)_2|\leq \Theta_0 + |x_4|\eta_1<1$.
\item[$\bullet$] Next assume that $\cZ_x=\{2\}$. Then we use the directions
\begin{equation*}
\pm(1,\delta_1,\eta_1,0),\ \ \pm(\delta_1,-\eta_1,1,0),\ \ \pm(\delta_1,\eta_1,-1,0),\ \ \pm(0,\pm(\eta_1,\delta_1),1)
\end{equation*}
to illuminate $x$. When $|x_1|>\Theta_0$ and $|x_3| < \frac{1-\Theta_0}{4}$, we pick $d\in\{\pm(1,\delta_1,\eta_1,0)\}$ so that $d_1\cdot x_1 < 0$: we have that $x+|x_1|d\in \intr\B$, which can be seen by comparing to the point $\frac{1+\Theta_0}{2}e_4 + \frac{1-\Theta_0}{2}(e_2+e_3)$ again (recall that $|x_1| > \Theta_0$ implies that $|x_4|\leq \Theta_0$).

If instead $|x_1| > \Theta_0$ and $|x_3| \geq \frac{1-\Theta_0}{4}$, we pick $d^\prime$ from $\pm(\delta_1,-\eta_1,1,0),\, \pm(\delta_1,\eta_1,-1,0)$ so that $d^\prime_s\cdot x_s < 0$ for $s\in \{1,3\}$. We can then compare the displaced vector $x+\frac{1-\Theta_0}{4}d^\prime$ with a convex combination of the form
\begin{equation*}
(1-\lambda)\big|\vec{x}\big| + \lambda\left(\frac{1-\Theta_0}{2}e_2 + \frac{1+\Theta_0}{2}e_4\right)
\end{equation*}
where $\lambda <\frac{1-\Theta_0}{4}\delta_1$. Then, as long as $\eta_1 < 2\lambda < \frac{1-\Theta_0}{2}\delta_1$, we will obtain that $x+\frac{1-\Theta_0}{4}d^\prime\in \intr\B$.

\medskip

Finally, if $|x_1| \leq \Theta_0$, then we pick $d\in \{\pm(0,\pm(\eta_1,\delta_1),1)\}$ so that $d_s\cdot x_s < 0$ for $s\in \{3,4\}$. We will have that $x+|x_4|d\in \intr\B$, which can be seen by comparing with the point
\begin{equation*}
(1-|x_4|\delta_1)(e_1+e_3) + |x_4|\eta_1\,e_2
\end{equation*}
which is also an interior point of $\B$ since $1-|x_4|\delta_1 + |x_4|\eta_1 < 1$.
\item[$\bullet$] It remains to consider the subcases $\cZ_x=\{3\}$ here: we will now illuminate $x$ using the directions
\begin{equation*}
\pm(1,\delta_1,\eta_1,0),\ \ \pm(\delta_1,-1,-\eta_1,0),\ \ \pm(0,\eta_1,\delta_1,1).
\end{equation*}
Similarly to the previous subcase, we first assume that $\max(|x_1|,|x_2|) > \Theta_0$. Then we will have that $|x_4|\leq \Theta_0$. If we also have that $\min(|x_1|,|x_2|) < \frac{1-\Theta_0}{4}$, and we write $i$ for the index at which $\max(|x_1|,|x_2|)$ is attained, then we simply pick $d$ from the first 4 directions above so that $m.c.(d) = i$ and so that $d_i \cdot x_i < 0$.

\medskip

If instead $\min(|x_1|,|x_2|) \geq \frac{1-\Theta_0}{4}$, then we pick $d^\prime$ from the first 4 directions again, but this time so that $d^\prime_s\cdot x_s < 0$ for $s\in [2]$. Similarly to above, we consider the displaced vector $x+\frac{1-\Theta_0}{4}d^\prime$, and conclude that it is in $\intr\B$ as long as $\eta_1 < \frac{1-\Theta_0}{2}\delta_1$.

\bigskip

Finally, if $\max(|x_1|,|x_2|)\leq \Theta_0$, then we pick $d\in \{\pm(0,\eta_1,\delta_1,1)\}$ so that $d_4\cdot x_4 < 0$. We will have that $x+|x_4|d\in \intr\B$, which can be seen by comparing with the point $\frac{1+\Theta_0}{2}(e_1+e_2) +\frac{1-\Theta_0}{2}e_3$.
\end{itemize}
\item[$\blacklozenge$] We now suppose that $\cZ_x=\emptyset$. If $|x_4|\leq \Theta_0$, then, as in previous propositions and subcases, we illuminate $x$ using the first 8 directions of ${\cal F}_{\ref{prop:R^4-exactly-three-pairs},\ref{prop:R^4-exactly-four-pairs},\delta_1, \eta_1}$ (which capture all combinations of signs for the first three coordinates).

\medskip

Next assume that $|x_4| > \Theta_0$. Then $\max(|x_1|,|x_2|)\leq \Theta_0$. Hence, we can pick $d\in \{\pm(0,\pm(\eta_1,\delta_1),1)\}$ so that $d_s\cdot x_s < 0$ for $s\in \{3,4\}$. If $|x_3| < \frac{1-\Theta_0}{4}$, then we simply compare the displaced vector $x+|x_4|d$ with the vector $\frac{1+\Theta_0}{2}(e_1+e_2) + \frac{1-\Theta_0}{2}e_3$.

\smallskip

On the other hand, if $|x_3| \geq \frac{1-\Theta_0}{4}$, we consider the displaced vector $x+\frac{1-\Theta_0}{4}d$, and compare it with a convex combination of the form
\begin{equation*}
(1-\lambda)\big|\vec{x}\big| + \lambda(e_1+e_2)
\end{equation*}
where $\lambda < \frac{1-\Theta_0}{4}\delta_1$. We will have that $x+\frac{1-\Theta_0}{4}d\in \intr\B$, as long as $\frac{1-\Theta_0}{4}\eta_1 < \lambda(1-\Theta_0) \leq \lambda(1-|x_2|)$ $\Leftrightarrow\,\eta_1 < 4\lambda < (1-\Theta_0)\delta_1$, which is already guaranteed by the restrictions we imposed on $\eta_1$.
\end{itemize}
\noindent This completes the proof of Case 1.

\bigskip

\emph{Proof for Case 2.} Recall that we have set $\beta_0:=\max\bigl\{\|e_1+e_2\|_\B^{-1}, \|e_3+e_4\|_\B^{-1}\bigr\}$, and we know that $\beta_0 > \frac{1}{2}$ (since we assumed that $\B$ is NOT an affine image of the cube). Recall also that we supposed WLOG that $\beta_0 = \|e_3+e_4\|_\B^{-1}$. We will then show that ${\cal F}_{\ref{prop:R^4-exactly-four-pairs}, {\rm alt}, \delta_2,\eta_2}$ illuminates $\B$ for some suitably chosen $\delta_2, \eta_2$. In this main case, we need some preparatory/key observations first. 

\bigskip

{\bf Observation 1 for Case 2.} Since $\frac{1}{2}(e_3+e_4)\in \intr\B$, we get that, for every $\epsilon\in (0,1)\setminus\{\frac{1}{2}\}$, the point
\begin{equation*}
(1-\epsilon)e_3 + \epsilon e_4
\end{equation*}
is also an interior point of $\B$. Indeed, assume first that $\epsilon < \frac{1}{2}$, and set $\lambda= 2\epsilon$ (in which case $\lambda\in (0,1)\,\!$). Then 
\begin{equation*}
(1-\lambda) e_3 + \lambda\,\tfrac{1}{2}(e_3+e_4)\in \intr\B
\end{equation*}
because it is a non-trivial convex combination of points in $\B$ with one of them being interior. But
\begin{equation*}
(1-\lambda) e_3 + \lambda\,\tfrac{1}{2}(e_3+e_4) = \left(1-\tfrac{\lambda}{2}\right)e_3 + \tfrac{\lambda}{2}e_4 = (1-\epsilon)e_3 + \epsilon e_4.
\end{equation*}
Analogously we show the result if $\epsilon \in (\frac{1}{2},\,1)$, by considering convex combinations of $\frac{1}{2}(e_3+e_4)$ with $e_4$.

\bigskip

{\bf Observation 2 for Case 2.} For every $a,\epsilon\in (0,1)$, we have that the points
\begin{equation*}
(a,0,1-\epsilon,\epsilon) \quad \hbox{and}\quad (0,a,1-\epsilon,\epsilon)
\end{equation*}
are interior points of $\B$. Indeed, by the previous key observation we know that the point $(0,0,1-\epsilon,\epsilon)\in \intr\B$. At the same time $\B$ contains the point
\begin{equation*}
(1,0,1-\epsilon,\epsilon) = (1-\epsilon) (e_1+e_3) + \epsilon (e_1+e_4).
\end{equation*}
But then
\begin{equation*}
(a,0,1-\epsilon,\epsilon) = a(1,0,1-\epsilon,\epsilon) + (1-a) (0,0,1-\epsilon,\epsilon),
\end{equation*}
which shows that it is an interior point of $\B$. Similarly we check that $(0,a,1-\epsilon,\epsilon)\in \intr\B$.

\bigskip

For the rest of the proof we fix $\delta_2 < \frac{1-\beta_0}{4}$, and $\eta_2 < (1-\beta_0)\delta_2$. 
We are ready to illuminate the boundary points of $\B$, and as before, we only focus on \underline{potentially} extreme points $x\in \B$. By our current assumptions for $\B$, there are certainly no such points with $\abs{\cZ_x}=3$, so we move on with the remaining possibilities for $\abs{\cZ_x}$.
\begin{itemize}
\item[$\blacklozenge$] Assume that $\abs{\cZ_x}=2$. If $x= \pm e_1+ \pm e_3$, then we use the directions
\begin{equation*}
\pm(1,-\eta_2,-\delta_2,-\delta_2)\ \ \pm(\delta_2,0,1,-\eta_2)
\end{equation*}
to illuminate $x$. Indeed, e.g.
\begin{equation*}
e_1-e_3 + (-1,\eta_2,\delta_2,\delta_2) = (0,\,\eta_2,\,-1+\delta_2,\,\delta_2)
\end{equation*}
which is in $\intr\B$ by the 2nd key observation.

\smallskip

On the other hand, if e.g. $x= -e_1 -e_3$, then $x+(\delta_2,0,1,-\eta_2)= (-1+\delta_2,0,0,-\eta_2)$, which we can immediately confirm is an interior point of $\B$ by comparing it to $e_1+e_4$.

\medskip

In a very analogous manner, we can illuminate all the points $\pm e_1 +\pm e_4,\,$ $\pm e_2 +\pm e_3$ and $\pm e_2+ \pm e_4$, and then we will have also illuminated every other point in their convex hull.

\bigskip

Assume now that $x= x_1 e_1 + x_2 e_2$. Then $\min(|x_1|,|x_2|)\leq \beta_0$. If $i$ is the index at which $\max(|x_1|,|x_2|)$ is attained, and $\{j\}=\{1,2\}\setminus\{i\}$, then we pick a direction $d$ from $\pm(1,-\eta_2,-\delta_2,-\delta_2),\,\pm(-\eta_2,1,-\delta_2,-\delta_2)$ so that $m.c.(d) = i$ and $d_i \cdot x_i < 0$. We will have that $x+|x_i|d\in \intr\B$, which can be readily seen if we compare with the point $\left(\beta_0 + \frac{1-\beta_0}{4}\right)e_j + \frac{1-\beta_0}{4}(e_3+e_4)\in \intr\B$.

\medskip

Similarly, if $x=x_3e_3+x_4e_4$, then $\min(|x_3|,|x_4|)\leq \beta_0$. If $r$ is the index at which $\max(|x_3|,|x_4|)$ is attained, and $\{t\}=\{3,4\}\setminus\{r\}$, then we pick a direction $d^\prime$ from $\pm(\delta_2,0,1,-\eta_2),\,\pm(\delta_2,0,-\eta_2,1)$ so that $m.c.(d^\prime) = r$ and $d^\prime_r\cdot x_r < 0$. We will have that $x+|x_r|d^\prime\in\intr\B$, which can be seen by comparing with the point $\frac{1+\beta_0}{2}(e_1+e_t)\in \intr\B$.
\medskip
\item[$\blacklozenge$] Assume now that $\abs{\cZ_x}=1$. If $\cZ_x = \{2\}$, then we consider the following subcases:
\begin{itemize}
\item[$\bullet$] $\sign(x_3) = - \sign(x_4)$. Then one of the directions $\pm(\delta_2,0,1,-\eta_2),\,\pm(\delta_2,0,-\eta_2,1)$ illuminates $x$.
\item[$\bullet$] $\sign(x_3) = \sign(x_4)$. We also recall that $\min(|x_3|,|x_4|)\leq \beta_0$; write $i$ for the index at which $\max(|x_3|,|x_4|)$ is attained, and $j$ for the other index.

\medskip

If at the same time $\sign(x_1) = -\sign(x_3)=-\sign(x_4)$, then the unique direction $d\in \{\pm(1,-\eta_2,-\delta_2,-\delta_2)\}$ satisfying $d_s\cdot x_s < 0$ for $s\in \{1,3,4\}$ illuminates $x$. Indeed, if $|x_j|= \min(|x_3|,|x_4|) < |x_1|\delta_2$, then $x+|x_1|d$ satisfies:
\begin{itemize}
\item[$\cdot$] $(x+|x_1|d)_1 = 0$,
\item[$\cdot$] $|(x+|x_1|d)_i|\leq \max\bigl(|x_i|-|x_1|\delta_2,\,|x_1|\delta_2-|x_i|\bigr)\leq \max\bigl(|x_i|-|x_1|\delta_2,\,|x_1|\delta_2\bigr)\leq 1-|x_1|\delta_2$,
\item[$\cdot$] $|(x+|x_1|d)_j|\leq \max\bigl(|x_j|-|x_1|\delta_2,\,|x_1|\delta_2-|x_j|\bigr)\leq |x_1|\delta_2$,
\item[$\cdot$] and $|(x+|x_1|d)_2|\leq |x_1|\eta_2 < 1$.
\end{itemize}
Thus $x+|x_1|d$ has smaller (in absolute value) corresponding coordinates compared to the vector \[|x_1|\eta_2 e_2 \,+\,(1-|x_1|\delta_2)e_i \,+\, |x_1|\delta_2 e_j.\] It remains to recall that the latter point is in $\intr\B$ because of the 2nd key observation.

\medskip

On the other hand, if $|x_j|= \min(|x_3|,|x_4|) \geq |x_1|\delta_2$, then we compare $x+|x_1|d$ with a point of the form
\begin{equation*}
(1-\lambda)\big|\vec{x}\big| + \lambda e_2
\end{equation*}
where $\lambda < |x_1|\delta_2$. Then, since we have assumed that $\eta_2 < (1-\beta_0)\delta_2 < \delta_2$, we can choose $\lambda$ so that $|x_1|\eta_2 < \lambda < |x_1|\delta_2$, which will then imply that $x+|x_1|d\in \intr\B$.

\bigskip

Next we consider the cases where $\sign(x_1)=\sign(x_3)=\sign(x_4)$. Then we pick $d^\prime$ from $\pm(\delta_2,0,1,-\eta_2),\,\pm(\delta_2,0,-\eta_2,1)$ so that $d^\prime_s\cdot x_s < 0$ for $s\in \{1,i\}$ (recall that we write $i\in \{3,4\}$ for the (smallest) index at which $\max(|x_3|,|x_4|)$ is attained, and $j$ for the remaining index). Then $x+|x_i|d^\prime$ satisfies:
\begin{itemize}
\item[$\cdot$] $(x+|x_i|d^\prime)_i=0 = (x+|x_i|d^\prime)_2$,
\item[$\cdot$] $|(x+|x_i|d^\prime)_1| \leq 1-|x_i|\delta_2 < 1$,
\item[$\cdot$] and $|(x+|x_i|d^\prime)_j| \leq \beta_0 + |x_i|\eta_2 < 1$.
\end{itemize}
It follows that $x+|x_i|d^\prime\in \intr\B$ since its non-zero coordinates are \underline{strictly} smaller than the corresponding coordinates of $e_1+e_j$.
\end{itemize}
\begin{itemize}
\item[$\bullet$] If $\cZ_x = \{1\}$, then we illuminate $x$ in a completely symmetric way compared to the previous subcase, by using one of the directions $\pm(-\eta_2,1,-\delta_2,-\delta_2),\,\pm(0,\delta_2,1,-\eta_2),$ $\pm(0,\delta_2,-\eta_2,1)$.
\item[$\bullet$] If $\cZ_x=\{3\}$, then we illuminate $x$ using one of the directions
\begin{equation*}
\pm(1,-\eta_2,-\delta_2,-\delta_2),\ \ \pm(-\eta_2,1,-\delta_2,-\delta_2),\ \ \pm(\delta_2,0,-\eta_2,1),\ \ \pm(0,\delta_2,-\eta_2,1).
\end{equation*}
We recall that $\min(|x_1|,|x_2|)\leq \beta_0$; let us write $i$ for the index at which $\max(|x_1|,|x_2|)$ is attained, and $j$ for the other index in $[2]$. 

\medskip

If $\sign(x_i)=-\sign(x_4)$, then we pick $d$ from $\pm(1,-\eta_2,-\delta_2,-\delta_2),\,\pm(-\eta_2,1,-\delta_2,-\delta_2)$ so that $m.c.(d) = i$ and so that $d_s\cdot x_s < 0$ for $s\in \{i,4\}$. Then $x+|x_i|d$ will have smaller (in absolute value) corresponding coordinates compared to the point
\begin{equation*}
(\beta_0+|x_i|\eta_2)e_j + |x_i|\delta_2 e_3 + (1-|x_i|\delta_2)e_4,
\end{equation*}
which is in $\intr\B$ by the 2nd key observation. Hence $x+|x_i|d\in \intr\B$ too.

\medskip

If instead $\sign(x_i)=\sign(x_4)$, then we pick $d^\prime$ from $\pm(\delta_2,0,-\eta_2,1),\ \ \pm(0,\delta_2,-\eta_2,1)$ so that $|d^\prime_i|=\delta_2$ and so that $d^\prime_s\cdot x_s < 0$ for $s\in \{i,4\}$. If $|x_i| < \frac{1-\beta_0}{4}$, then we compare $x+|x_4|d^\prime$ with the vector
\begin{equation*}
\frac{1-\beta_0}{4}\bigl(e_i+e_j+e_3\bigr)\in \intr\B.
\end{equation*} 
On the other hand, if $|x_i| \geq \frac{1-\beta_0}{4}$, then we compare $x+|x_4|d^\prime$ with a convex combination of the form
$(1-\lambda^\prime)\big|\vec{x}\big| + \lambda^\prime(e_j+e_3)$ where $\lambda^\prime < |x_4|\delta_2$. Our assumption that $\eta_2 < \delta_2$ implies that we can choose such a $\lambda^\prime$ so that $|x_4|\eta_2 < \lambda^\prime < |x_4|\delta$, which in turn implies that $x+|x_4|d^\prime\in \intr\B$.
\item[$\bullet$] We illuminate $x$ in a symmetric fashion when $\cZ_x=\{4\}$, by using one of the directions
\begin{equation*}
\pm(1,-\eta_2,-\delta_2,-\delta_2),\ \ \pm(-\eta_2,1,-\delta_2,-\delta_2),\ \ \pm(\delta_2,0,1,-\eta_2),\ \ \pm(0,\delta_2,1,-\eta_2).
\end{equation*}
\end{itemize}
\item[$\blacklozenge$] Finally, assume that $\cZ_x = \emptyset$. We know that $\min(|x_1|,|x_2|) \leq \beta_0$, and the same inequality holds true for $\min(|x_3|,|x_4|)$. Let us write $i$ for the index at which $\max(|x_1|,|x_2|)$ is attained, and $j$ for the other index in $[2]$. Similarly, let us write $r$ for the index at which $\max(|x_3|,|x_4|)$ is attained, and $t$ for the other index in $\{3,4\}$.
\begin{itemize}
\item[$\bullet$] Assume first that $\sign(x_3)=-\sign(x_4)$. Then, by Corollary \ref{cor:affine-set}, $x$ is illuminated by the unique direction $d$ among the last 8 directions in ${\cal F}_{\ref{prop:R^4-exactly-four-pairs}, {\rm alt}, \delta_2,\eta_2}$ which satisfies $d_i\neq 0$ and $d_s\cdot x_s < 0$ for $s\in \{i,3,4\}$.
\item[$\bullet$] Next, assume that $\sign(x_3)=\sign(x_4)$. 
\begin{itemize}
\item[--] If $\sign(x_i) = \sign(x_3)=\sign(x_4)$, then we pick $d^\prime$ from the last 8 directions of ${\cal F}_{\ref{prop:R^4-exactly-four-pairs}, {\rm alt}, \delta_2,\eta_2}$ so that $m.c.(d^\prime) = r$, $|d^\prime_i|=\delta_2$ and $d^\prime_s\cdot x_s < 0$ for $s\in \{i,r\}$. 
\begin{itemize}
\item[$\cdot$] If in addition $|x_i| < \frac{1-\beta_0}{4}$, then we simply compare the displaced vector $x+|x_r|d^\prime$ to the point
\begin{equation*}
\frac{1-\beta_0}{4}(e_i+e_j) + \frac{1+\beta_0}{2}e_t = \frac{1-\beta_0}{4}(e_1+e_2) + \frac{1+\beta_0}{2}e_t\in \B
\end{equation*}
to conclude that $x+|x_r|d^\prime\in \intr\B$ (we can do this because, by our assumptions here, $|(x+|x_r|d^\prime)_i| \leq \max(|x_i|-|x_r|\delta_2,\,|x_r|\delta_2-|x_i|)\leq \max(|x_i|,\,\delta_2) < \frac{1-\beta_0}{4}$ and $|x_j|\leq |x_i| < \frac{1-\beta_0}{4}$).
\item[$\cdot$] If instead $|x_i| \geq \frac{1-\beta_0}{4}$, then we compare $x+|x_r|d^\prime$ to a convex combination of the form
\begin{equation*}
(1-\lambda)\big|\vec{x}\big| + \lambda(e_j+e_t)
\end{equation*}
where $\lambda < |x_r|\delta_2$. As long as $|x_r|\eta_2 < \lambda(1-\beta_0)\leq \lambda(1-|x_t|)$ (which a suitably chosen $\lambda$ can satisfy given the restriction $\eta_2 <(1-\beta_0)\delta_2$), we will have that $x+|x_r|d^\prime\in \intr\B$.
\end{itemize} 
\medskip
\item[--] It remains to consider the cases where $\sign(x_i) = -\sign(x_3)=-\sign(x_4)$. Then we pick the unique direction $d$ from $\pm(1,-\eta_2,-\delta_2,-\delta_2),\,\pm(-\eta_2,1,-\delta_2,-\delta_2)$ so that $m.c.(d) = i$ and so that $d_s\cdot x_s < 0$ for $s\in \{i,3,4\}$. 
\begin{itemize}
\item[$\cdot$] If $|x_t|=\min(|x_3|,|x_4|) < |x_i|\delta_2$, then, similarly to one of the subcases of the settings where $\cZ_x=\{1\}$ or $\cZ_x=\{2\}$, we will have that the displaced vector $x+|x_i|d$ has smaller (in absolute value) corresponding coordinates compared to the vector
\begin{equation*}
\left(\beta_0+\tfrac{1-\beta_0}{4}\right)e_j + (1-|x_i|\delta_2)e_r + |x_i|\delta_2 e_t
\end{equation*}
which is an interior point of $\B$ itself, by the 2nd key observation. Thus $x+|x_i|d\in \intr\B$.
\item[$\cdot$] If $|x_t|=\min(|x_3|,|x_4|) \geq |x_i|\delta_2$, then we compare $x+|x_i|d$ with a convex combination of the form
\begin{equation*}
(1-\lambda)\big|\vec{x}\big| + \lambda e_j
\end{equation*}
where $\lambda < |x_i|\delta_2$. As long as $|x_i|\eta_2 < \lambda(1-\beta_0) \leq \lambda(1-|x_j|)$, which is again possible for some $\lambda\in (0,|x_i|\delta_2)$ because of the restriction $\eta_2 < (1-\beta_0)\delta_2$, we will have that $x+|x_i|d\in \intr\B$.
\end{itemize}
\end{itemize}
\end{itemize}
\end{itemize}
This completes the proof of Case 2 as well.
\end{proof}

\bigskip

We can also make a similar note to Remark \ref{rem:alt-set-Case12-3squares}.

\begin{remark}
By slightly adjusting the proof of Case 1 of this last proposition, we can also confirm that any set $\B\in \U^4$ which contains the points $e_1+e_2, e_1+e_3, e_2+e_3$ and $e_3+e_4$, but does not contain $e_1+e_4$ and $e_2+e_4$ (nor does it contain the `triple' $e_1+e_2+e_3$) can be illuminated by the set 
\begin{align*}
{\cal F}_{\ref{prop:R^4-no-pairs},\ref{prop:R^4-exactly-five-pairs},\delta, \eta,\zeta}:= \bigl\{\!&\pm(1,\delta,\eta,0),\ \pm(\delta,-1,-\eta,0),\ \pm(\delta,-\eta,1,\zeta),\ \pm(-\delta,-\eta,1,\zeta),
\\
&\quad \pm(0,\pm(\eta,\delta),1),\ \pm(0,1,-\delta,\eta)\bigr\}
\end{align*}
instead of the set
\begin{align*}
{\cal F}_{\ref{prop:R^4-exactly-three-pairs},\ref{prop:R^4-exactly-four-pairs},\delta_1, \eta_1}:= \bigl\{\!&\pm(1,\delta_1,\eta_1,0),\ \pm(\delta_1,-1,-\eta_1,0),\ \pm(\delta_1,-\eta_1,1,0),\ \pm(\delta_1,\eta_1,-1,0),
\\
&\quad \pm(0,\pm(\eta_1,\delta_1),1),\ \pm(0,1,-1,0)\bigr\}.
\end{align*}
\end{remark}


\medskip
\medskip

\noindent \textsc{Department of Mathematical and Statistical Sciences,
University of Alberta, CAB 632, Edmonton, AB, Canada T6G 2G1}

\smallskip

\noindent 
{\it E-mail addresses:} \texttt{wrsun@ualberta.ca, vritsiou@ualberta.ca}

\end{document}